\documentclass[a4paper, 11pt]{article}
\usepackage[utf8x]{inputenc}
\usepackage{graphicx,amsfonts,amsmath,amssymb,amsthm}
\usepackage[sort, numbers]{natbib}
\usepackage{dsfont}
\usepackage{setspace}
\usepackage[a4paper,top=2.5cm,bottom=2.5cm,left=2.5cm,right=2.5cm,marginparwidth=1.75cm]{geometry}
\usepackage{fullpage}
\usepackage{pgfplots}
\pgfplotsset{compat=1.9}
\usepackage{enumitem}
\usepackage{commands}
\usepackage{tikzSettings}
\usepackage{caption}
\usepackage{subcaption}
\usepackage{url}
\usepackage{thmtools}
\usepackage{thm-restate}
\usepackage{longtable}
\usepackage{booktabs}
\usepackage{multirow}
\usepackage{siunitx}

\DeclareMathOperator*{\argmax}{\arg\!\max}
\usetikzlibrary{patterns,automata,arrows,arrows.meta, calc, quotes, tikzmark,shapes.misc,positioning,3d,pgfplots.statistics, pgfplots.colorbrewer}
\pgfplotsset{compat = 1.15, cycle list/Set1-8}

\newcommand{\boxplot}[7]{%
	\filldraw[fill=#7, opacity=0.75,line width=0.2mm] let \n{boxxl}={#1-0.1}, \n{boxxr}={#1+0.1} in (axis cs:\n{boxxl},#3) rectangle (axis cs:\n{boxxr},#4);   
	\draw[line width=0.2mm, color=black] let \n{boxxl}={#1-0.1}, \n{boxxr}={#1+0.1} in (axis cs:\n{boxxl},#2) -- (axis cs:\n{boxxr},#2);             
	\draw[line width=0.2mm] (axis cs:#1,#4) -- (axis cs:#1,#6);                                                                           
	\draw[line width=0.2mm] let \n{whiskerl}={#1-0.025}, \n{whiskerr}={#1+0.025} in (axis cs:\n{whiskerl},#6) -- (axis cs:\n{whiskerr},#6);        
	\draw[line width=0.2mm] (axis cs:#1,#3) -- (axis cs:#1,#5);                                                                           
	\draw[line width=0.2mm] let \n{whiskerl}={#1-0.025}, \n{whiskerr}={#1+0.025} in (axis cs:\n{whiskerl},#5) -- (axis cs:\n{whiskerr},#5);        
	}

\definecolor{aqua}{rgb}{0.0, 1.0, 1.0}
\newcommand*\rot{\rotatebox{60}}
\DeclarePairedDelimiter\floor{\lfloor}{\rfloor}

\author{Victor Cohen$^1$,Axel Parmentier$^1$}
\title{Integer programming for weakly coupled stochastic dynamic programs with partial information}
\date{%
    $^1$Ecole des Ponts Paristech, CERMICS, F-77455 Marne-la-Vall\'ee, France\\%
    \today}

\providecommand{\keywords}[1]
{
  \small	
  \textbf{\textit{Keywords---}} #1
}
\begin{document}
\maketitle

\begin{abstract}
	This paper introduces algorithms for problems where a decision maker has to control a system composed of several components and has access to only partial information on the state of each component.
	Such problems are difficult because of the partial observations, and because of the curse of dimensionality that appears when the number of components increases.
	Partially observable Markov decision processes (POMDPs) have been introduced to deal with the first challenge, while weakly coupled stochastic dynamic programs address the second.
	Drawing from these two branches of the literature, we introduce the notion of weakly coupled POMDPs.
	The objective is to find a policy maximizing the total expected reward over a finite horizon.
	Our algorithms rely on two ingredients.
	The first, which can be used independently, is a mixed integer linear formulation for generic POMDPs that computes an optimal memoryless policy.
	The formulation is strengthened with valid cuts based on a probabilistic interpretation of the dependence between random variables, and its linear relaxation provide a practically tight upper bound on the value of an optimal history-dependent policies.
	The second is a collection of mathematical programming formulations and algorithms which provide tractable policies and upper bounds for weakly coupled POMDPs.
	Lagrangian relaxations, fluid approximations, and almost sure constraints relaxations enable to break the curse of dimensionality.
	We test our generic POMDPs formulations on benchmark instance forms the literature, and our weakly coupled POMDP algorithms on a maintenance problem. 
	Numerical experiments show the efficiency of our approach.
\end{abstract}



%

\section{Introduction}
\label{sec:intro}

Many real world situations involve the control of a stochastic system composed of $M$ components on which the decision maker has only partial informations.
Typically, at each time-step, a component $m$ in $[M]$ is in a state $s^m$ in a state space $\calX_S^m$.
The decision maker does not know the exact state $s^m$ of component $m$, but has access to a noisy observation $o^m$ providing only partial information on $s^m$.
Based on these observations, the decision maker chooses the actions $a^m$ that should be performed on each components of the system.
This decision typically involves allocating some limited resources between the different components.
Each component then transitions to the new state $s'^m$ with a given probability, and the decision maker receives a reward $\sum_{m} r^m(s^m,a^m,s'^m)$ where $r^m$ is the individual reward of component $m$. 
The goal of the decision maker is to find a policy $\bfdelta$ that prescribes which actions $(a^m)_{m \in [M]}$ to take given observations $(o^m)_{m \in [M]}$ in order to maximize the expected total reward over a finite horizon.

Typical applications are (1) predictive maintenance problems, (2) multi-armed restless bandits and their applications to clinical trials, (3) inventory problems with inventory records inaccuracy, and (4) nurse assignments problems. See Appendix~\ref{app:examples} for a detailed treatment of these applications.
As an illustration, consider the case of airplane predictive maintenance.
An airplane regularly undergoes maintenance at slots known in advance.
At the beginning of each maintenance slot, the decision maker must decide which equipments it will maintain given its limited maintenance resources (machines, skilled technicians, etc.). 
On recent generations of airplanes, sensors signals are recorded during flights,
and a score that gives a noisy evaluation of the equipment wear is deduced from these signals.
The reward is in term of costs avoided. Each maintenance has a given cost.
And if a component fails in between two maintenance, the plane cannot take-off, and its flights of the day must be canceled leading to huge client refunding costs.
Based on the scores observed on each components, its maintenance resources, and an evaluation of the failures risks, the decision maker chooses which equipments to maintain.

When the system has only a single component whose state is observed by the decision maker, the problem is naturally modeled as a \emph{Markov Decision Process} (MDP) and well solved using dynamic programming algorithms. 
Controlling a partially observed multi-equipments system is more difficult because (1) the decision maker has only access to partial information, and (2) the system is composed of several components.
Let us briefly recap how these difficulties have been addressed in the literature.

The problem of controlling a system with only partial information is naturally modeled as a
\emph{Partially Observable Markov Decision Process} (POMDP) over a finite horizon.
As detailed in the survey of~\citet{Cassandra2003ASO}, a wide range of applications have been modeled as POMDPs, among which maintenance problems~\citep{Eckles1968} or clinical decision making~\citep{Denton2018}. 
On such partially-observed systems, POMDP approaches typically provide lower cost policies than MDP approaches.
This performance comes from the statistical models used to describe reality: Markov chains for MDP, versus hidden Markov models (HMM) for POMDP.
When data is scarce and only small dimensional model can be learned, HMM provide a much needed additional flexibility. 

Finding an optimal policy of the POMDP problem over a finite horizon is known to be PSPACE-complete \citep{Tsitsiklis1987}.
State-of-the-art approaches
rely on the fact that a POMDP is equivalent to a continuous state MDP in the \emph{belief state} space \citep[Theorem 4]{Eckles1968}.
The belief state is the posterior probability distribution on the state space given all the past decisions and observations. 
Leveraging this result, several dynamic programming algorithms have been derived for POMDPs with finite \citep{Sondik1973} or infinite horizon \citep{Sondik1978}. 
More details about these kind of algorithms can be found for instance in the survey of \citet{POMDPMonahan1982}, or more recently the book of \citet{Krishnamurthy2016}.
While the exact algorithms become quickly intractable when the size of the state and observation spaces grow, some approximate algorithms have shown significant improvements. \citet{Hauskrecht2000} and \citet{Pineau2008} survey various approximations methods for POMDPs. 
For POMDP problems with a finite horizon, \citet{Walraven2019} recently proposed an effective algorithm with guarantees. \citet{Aras2007} propose a mixed integer programming approach to exactly solve POMDP problems over a finite horizon. 
However, solving such a program is computationally expensive even for small instances.
Part of this difficulty comes from the fact that an optimal POMDP policy depends on the all history of action and decisions.
Using a memoryless policy leads to a more tractable problem, still NP-hard \citep{Littman94memoryless} but no more PSPACE-hard.
And there is a broad class of systems where memoryless policies perform well \citep{Barto1983,Li2011}, even if it is not the case on some pathological cases \citet{Littman94memoryless}.
Furthermore, our problem suffers from the \emph{curse of dimensionality} (like MDPs) since the size of the spaces grow exponentially with the number of components.
MDPs are in theory well-solved by dynamic programming algorithms, but the curse of dimensionality makes such approaches impractical to control systems with several components.
Approximate dynamic programming algorithms \citep{Bertsekas2007,Farias01thelinear,Powell2011,Meuleau1998} provide generic methodologies to address this curse.
In another paradigm, reinforcement learning approaches \citep{Sutton1998} to build such algorithms is an active research area in the machine learning community.
In the operations research literature, the notion of weakly coupled dynamic programs \citep{Meuleau1998,hawkins2003langrangian,Adelman2008} and decomposable Markov decision processes \citep{bertsimas2016decomposable} have been introduced in order to catch the specific structure of multi-components, where component-specific actions on each components must be coordinated, or a single action affects all components, respectively. 
And \citet{Walraven2018} recently introduced a different POMDP problem with multiple components where the resource constraints enforce policies to induce a maximum expected resource consumption over the finite time horizon. 
Such models have been applied to stochastic inventory routing with limited vehicle capacity \citep{Kleywegt2002}, stochastic multi-product dispatch problems \citep{Papadaki2003}, scheduling problems \citep{Whittle1988}, resources allocation \citep{Gittins1979}, revenue management \citep{Topaloglu2009} among others \citep{hawkins2003langrangian}.
The specific structure of these problems can then be leveraged in mathematical programming formulations that use approximate value functions \citep{de2003linear} or approximate moments \citep{bertsimas2016decomposable} as variables, notably using the Lagrangian relaxation of non-anticipativity constraints \citep{Brown2010} or of the linking constraints \citep{Ye2014}.
All the mathematical programs in this paper can be formulated either using moments variables (or marginal probabilities) or value function variables.
Since they lead to better numerical results, we include in the paper the formulations with moment variables. 
The value-function formulations can be found in the PhD dissertation of the first author~\citep{VCohenThesis2020}.

When the decision maker has only access to a partial observation for each component, each subsystem is a POMDP and all the subsystems are linked by resource constraints on the actions taken on each components.
This structure requires to extend the notion of weakly coupled stochastic dynamic program to the notion of \emph{weakly coupled POMDP}, which has been introduced by
\citet{Parizi2019} on infinite horizon problems.
\citet{Parizi2019,Abbou2019} both propose approximate policies based on the MDP relaxations. 
Numerical experiments  on finite horizon problems in Appendix~\ref{sub:app_nums:implicit_policy} show that taking into account the fact that observations are partial improve the performance of such algorithms.
For multi-armed bandits, the optimization problem has been introduced by \citet{MeshramMG18} as \emph{restless partially observable multi-armed bandit} and new index policies have been proposed recently to solve it \citep{MeshramMG18,KazaMMM2019}. 
However, the algorithms and results proposed hold when the state spaces and observation spaces contain at most two elements, and the algorithms do not scale to larger spaces.

There is therefore a need for more efficient algorithms for weakly coupled POMDPs that address both the partially observable aspect and the curse of dimensionality.
To the best of our knowledge, mathematical programming techniques have not been exploited for the control of partially observed systems with multiple components. 
The purpose of this paper is to show that such techniques lead to practically efficient algorithms for such systems, both in the single and the multi-component cases.
Our first contributions are mathematical programming formulations for generic POMDPs.
\begin{enumerate}
	\item We propose an exact Non-Linear Program (NLP) formulation and an exact Mixed Integer Linear Program (MILP) formulation for POMDP problems with memoryless policies. 
	When solved with off-the-shelve solvers, they provide a practically efficient solution approach.
	Numerical experiments on instances from the literature show that, on several class of problems over a finite horizon among which maintenance problems, our MILP approach provides better solutions than a state-of-the-art POMDP algorithms such as SARSOP (which is not restricted to memoryless policies).
	\item We introduce an extended formulation with valid inequalities that improve the resolution of our MILP. Such inequalities come from a probabilistic interpretation of the dependence between random variables.
	Numerical experiments show their efficiency.
	\item We prove that the linear relaxation of our MILP is equivalent to the MDP relaxation of the POMDP, which corresponds to the case where the action depend on the current state. In addition, we show that the strengthened linear relaxation is a (practically tight) upper bound of the optimal value of POMDP with history-dependent policies. 
\end{enumerate}
We then show how, as for the fully observed cases, mathematical programming techniques can successfully exploit the specific structure of multi-components systems to address the curse of dimensionality that affect them. 
More precisely, we introduce several formulations and algorithms for weakly coupled POMDP with memoryless policies that leverage our formulations and valid inequalities for generic POMDPs.
\begin{enumerate}[resume]
	\item Our main contribution is a history-dependent policy for weakly coupled POMDPs, each value of which is defined through the resolution of a MILP approximation. 
	To the best of our knowledge, it is the first algorithm that can address large scale instances over finite horizon: When approaches in the literature considered instances with $10$ components each of them with $2$ states, we can provide policies with at most $10\%$ optimality gap on instances with $20$ components and $5$ states within a reasonable computation time.
	\item We introduce tractable upper bounds on the value of an optimal memoryless policy and on the value of an optimal history-dependent policy. 
	The first one is a Lagrangian relaxation bound computed thanks to a column generation algorithm, and the second comes from an extended linear programming formulation.
	These bounds are compute  enable to compute the optimality gaps mentioned in the previous point.
	\item We provide a shared hierarchy of lower and upper bounds on the value of optimal memoryless and history dependent policies, decomposable policies, the value our formulations, and the value of natural information relaxations.

	\item Detailed numerical experiments on maintenance and multi-armed bandits problems evaluate the different algorithms and bounds proposed.
\end{enumerate}

The paper is organized as follows.
Section~\ref{sec:problem} recalls the formal definition of a POMDP and introduces the notion of weakly coupled POMDP.
Section~\ref{sec:pomdp} introduces our mathematical programming formulations (NLP and MILP) for POMDP with memoryless policies, the valid inequalities and the theoretical study of the linear relaxations.
Section~\ref{sec:wkpomdp} introduces our heuristic, which is based on approximate integer formulation, and also describes different mathematical programming formulations, including the Lagrangian relaxation, that give useful bounds for weakly coupled POMDPs.
Finally, Section~\ref{sec:num} summarizes our numerical experiments.
The proofs of all the theorems in Sections~\ref{sec:pomdp} and~\ref{sec:wkpomdp} are respectively available in Appendices~\ref{app:pomdp} and~\ref{app:wkpomdp:proofs}.
\section{Weakly coupled POMDP}
\label{sec:problem}


\subsection{Background on POMDPs}
\label{sub:problem:background_pomdp}
 

A POMDP is a multi-stage stochastic optimization problem defined as follows.
It models on a horizon $T$ in $\bbZ_+$ the evolution of a system.
At each time $t$ in $[T]$, the system is in a random state $S_t$, which takes value in a finite state space $\calX_S$.
The system starts in state $s$ in $\calX_S$ with probability $p(s) := \bbP\left(S_1 = s \right)$.
At time $t$, the decision maker does not have access to $S_t$, but observes $O_t$, whose value belongs to a finite state space $\calX_{O}$. 
When the system is in state $S_t=s$, it emits an observation $O_t=o$ with probability $\bbP(O_t=o | S_t=s) := p(o|s)$.
Then, the decision maker takes an action $A_t$, which belongs to a finite space $\calX_A$. Given an action $A_t=a$, the system transits from state $S_t=s$ to state $S_{t+1} = s'$ with probability $\bbP\left(S_{t+1}=s' | S_t=s, A_t=a \right) := p(s'|s,a)$, and the decision maker receives the immediate reward $r(s,a,s')$, where the reward function is defined as a real valued function $r \colon \calX_S \times \calX_A \times \calX_S \rightarrow \bbR$, which we will also view as a vector $\bfr = \left(r(s,a,s')\right) \in \bbR^{\calX_S\times\calX_A\times\calX_S}$.
We denote by $\pfrak$ the vector of probabilities $\pfrak = \left(p(s), p(o|s), p(s'|s,a) \right)_{\substack{s,s' \in \calX_S, o\in \calX_O \\ a \in \calX_A}}$.
A POMDP is parametrized by the fifth tuple $\left(\calX_S,\calX_O,\calX_A, \pfrak,\bfr \right)$.


\paragraph{POMDP problem.}
Let $\left(\calX_S,\calX_O,\calX_A, \pfrak,\bfr \right)$ be a POMDP.
Given a finite horizon $T\in \bbZ_{+}$, the choices made by the decision maker are modeled using a policy $\bfdelta=(\delta^1, \ldots, \delta^T)$, where $\delta^t$ is the conditional probability distribution of taking action $A_t$ at time $t$ given the history of observations and actions $H_t = \left(O_{1},A_{1},\ldots,A_{t-1},O_t\right)$ in $\calX_{H}^{t}:= \left(\calX_O \times \calX_A \right)^{t-1} \times \calX_O$, i.e., $\delta_{a|h}^t := \bbP(A_t=a|H_t=h),$ 
for any $a$ in $\calX_A$ and $h$ in $\calX_H^{t}$.
We denote by $\Deltahis$ the set of policies
\begin{align*}
\displaystyle \Deltahis = \bigg\{ \bfdelta \in \bbR^{\sum_{t=1}^T\calX_H^t \times \calX_A} \colon \sum_{a \in \calX_A} \delta^t_{a|h} = 1 \ \mathrm{and} \ \delta^t_{a|h} \geq 0, \forall h \in \calX_H^t, a \in \calX_A, t \in [T] \bigg\}.
\end{align*}
In $\Delta_{\mathrm{his}}$, ``his'' refers to policies that take into account the history of observations and actions by opposition to memoryless policies which will be introduced below.
A policy $\bfdelta \in \Deltahis$ leads to the probability distribution $\bbP_{\bfdelta}$ on $\left(\calX_S \times \calX_O \times \calX_A \right)^T \times \calX_S$ such that
\begin{align*}
    \bbP_{\bfdelta} \left((S_t = s_t, O_t=o_t, A_t=a_t)_{1 \leq t \leq T}, S_{T+1} = s_{T+1}\right) = 
    p(s_1) \prod_{t=1}^T p(o_t | s_t) p(s_{t+1} \vert s_t, a_t ) \delta^t_{a_t|h_t}, 
\end{align*}
where $h_t = (s_1,o_1,\ldots,a_{t-1},o_t)$.
We denote by $\bbE_{\bfdelta}$ the expectation according to $\bbP_{\bfdelta}$.
The goal of the decision maker is to find a policy $\bfdelta$ in $\Deltahis$ maximizing the expected total reward over the finite horizon $T$. The POMDP problem is:
\begin{equation}\label{pb:POMDP_perfectRecall}
    \max_{\bfdelta \in \Deltahis} \bbE_{\bfdelta} \bigg[ \sum_{t=1}^{T}r(S_t,A_t,S_{t+1})\bigg] \tag{$\rm{P}_{\rm{his}}$}
\end{equation}
It is known that~\ref{pb:POMDP_perfectRecall} is PSPACE-hard \citep{Tsitsiklis1987}.

\paragraph{POMDP problem with memoryless policies.}
Let $\left(\calX_S,\calX_O,\calX_A, \pfrak,\bfr \right)$ be a POMDP. Given a finite horizon $T \in \bbZ_{+}$, 
A memoryless policy is a vector $\bfdelta=(\delta^1, \ldots, \delta^T)$, where $\delta^t$ is the conditional probability distribution at time $t$ of action $A_t$ given observation $O_t$, i.e., $\delta_{a|o}^t := \bbP(A_t=a|O_t=o)$ 
for any $a$ in $\calX_A$ and $o$ in $\calX_O$. 
Such policies are said memoryless because the choice of $A_t$ only depends on the current observation $O_t$, in contrast with the history of observations and actions $H_t$. 
We denote by $\Deltaml$ the set of memoryless policies, where ``ml'' refers to memoryless
\begin{align}\label{eq:problem:def_policy_set}
\displaystyle \Deltaml = \bigg\{ \bfdelta \in \bbR^{T \times \calX_A \times \calX_O } \colon \sum_{a \in \calX_A} \delta^t_{a|o} = 1 \ \mathrm{and} \ \delta^t_{a|o} \geq 0,\enskip \forall o \in \calX_O, \enskip a \in \calX_A \bigg\}.
\end{align}
With an abuse of notation, the definition of $\Deltaml$ ensures that $\Deltaml \subseteq \Deltahis$ in the sense that for any policy $\bfdelta$ in $\Deltaml$, we can define a policy $\tilde{\bfdelta}$ in $\Deltahis$ with the same values, i.e., such that $\tilde{\delta}_{a|h}^t := \delta_{a|o}^t$ for any $a$ in $\calX_A$, $o$ in $\calX_O$, $h$ in $\calX_H^t$ and $t$ in $[T]$.
This time, the policy $\bfdelta \in \Deltaml$ endows $\left(\calX_S \times \calX_O \times \calX_A \right)^T \times \calX_S$ with the probability distribution $\bbP_{\bfdelta}$ characterized by
\begin{align}\label{eq:problem:probability_distrib}
    \bbP_{\bfdelta} \left((S_t = s_t, O_t=o_t, A_t=a_t)_{1 \leq t \leq T}, S_{T+1} = s_{T+1}\right) = p(s_1) \prod_{t=1}^T p(o_t | s_t) p(s_{t+1} \vert s_t, a_t ) \delta^t_{a_t|o_t}, 
\end{align}
and the decision maker now seeks a memoryless policy $\bfdelta$ in $\Deltaml$ maximizing the expected total reward over the finite horizon $T$. The POMDP problem with memoryless policy is therefore
\begin{equation}\label{pb:POMDP}
    \max_{\bfdelta \in \Deltaml} \bbE_{\bfdelta} \bigg[ \sum_{t=1}^{T}r(S_t,A_t,S_{t+1})\bigg]. \tag{$\rm{P}_{\rm{ml}}$}
\end{equation}
Problem~\ref{pb:POMDP} is NP-hard \citep{Littman94memoryless}.
In Section~\ref{sec:num}, we provide numerical experiments showing that memoryless policies perform well on different kinds of problems modeled as POMDPs.

\subsection{Problem definition of weakly coupled POMDP} 
\label{sub:problem:wkpomdp}


A weakly coupled POMDP models a system composed of $M$ components, each of them evolving independently as a POMDP.
Let $S_t^m$ and $O_t^m$ be random variables that represent respectively the state and the observation of component $m$ at time $t$, and that belong respectively to the state space $\calX_S^m$ and the observation space $\calX_O^m$ of component $m$. Each component is assumed to evolve individually as a POMDP. 
We denote by $\pfrak^m$ and $\bfr^m$ respectively the probability distributions and the immediate reward functions of component $m$.
We denote by $\bfS_t = \big( S_t^1, \ldots, S_t^M \big)$ and $\bfO_t = \big( O_t^1, \ldots, O_t^M \big)$ the state and the observation of the full system at time $t$, which lie respectively in the \emph{state space} $\calX_S = \calX_S^1 \times \cdots \times \calX_S^M$ and the \emph{observation space} $\calX_O = \calX_O^1 \times \cdots \times \calX_O^M$. 
The spaces $\calX_S$ and $\calX_O$ represent the state space and the observation space of the full system.
We assume that the action space $\calX_A$ can be written
\begin{align}\label{eq:problem:def_form_action_space}
    \calX_A = \left\{\mathbf{a} \in \calX_A^1 \times \cdots \times \calX_A^M \colon \sum_{m =1}^M \mathbf{D}^m(a^m) \leq \mathbf{b} \right\},
\end{align}
where $\calX_A^m$ is the individual action space of component $m$, and $\mathbf{D}^m : \calX_A^m \rightarrow \bbR^q$ is a given function for each component $m$ in $[M]$, and $\mathbf{b} \in \bbR^q$ is a given vector for some finite integer $q$. 
We assume in the rest of the paper that the vector $\bfb$ is non-negative.
This is without loss of generality since an instance with a generic $\bfb$ can be turned into an equivalent one with non-negative $\bfb$: It suffices to set $\bfD'^m := \bfD^m-\frac{k}{M}$ for every $m$ in $[M]$, $\bfb' := \bfb - k$ where $k=\min_{i \in [q]} b_i$ and $b_i$ is the $i$-th coordinate of vector $\bfb$, and $\calX_A = \left\{\bfa \in \calX_A^1 \times \cdots \times \calX_A^M \colon \sum_{m=1}^M \bfD'^m(a^m) \leq \bfb' \right\}$. 

Each component is assumed to evolve independently, hence the joint probability of emission factorizes as 
\begin{equation}\label{eq:problem:factor_emission}
    \bbP(\bfO_t = \bfo | \bfS_t = \bfs) =\prod_{m=1}^M p^m(o^m | s^m),
\end{equation}
 and the joint probability of transition factorizes as
\begin{equation}\label{eq:problem:factor_transition}
\bbP(\bfS_{1} = \bfs) = \prod_{m=1}^M p^m(s^m) \quad \text{and} \quad \bbP(\bfS_{t+1} = \bfs' | \bfS_t=\bfs, \bfA_t=\bfa) = \prod_{m=1}^M p^m(s'^m |s^m, a^m),
\end{equation}
for all $t$ in $[T]$. In addition, the reward is assumed to decompose additively
\begin{equation}\label{eq:problem:factor_rewards}
    r(\bfs,\bfa,\bfs') = \sum_{m=1}^M r^m(s^m, a^m, s'^m).
\end{equation}
Hence, the weakly coupled POMDP problem with memoryless policies is the following:
\begin{equation}\label{pb:decPOMDP_wc}
      \max_{\bfdelta \in \Deltaml} \bbE_{\bfdelta} \bigg[ \sum_{t=1}^{T}r(\bfS_t,\bfA_t,\bfS_{t+1}) \bigg] \tag{{$\rm{P}_{\rm{ml}}^{\rm{wc}}$}}
\end{equation}
where the expectation is taken according to $\bbP_{\bfdelta}$ defined in~\eqref{eq:problem:probability_distrib}.
Similarly, we define the weakly coupled POMDP problem with history-dependent policies $\rm{P}_{\rm{his}}^{\rm{wc}}$ by replacing $\Deltaml$ with $\Deltahis$. 
Note that unless $\calX_A = \emptyset$ there always exists a feasible policy of~\ref{pb:decPOMDP_wc}. 
A weakly coupled POMDP is fully parametrized by $\left((\calX_S^m,\calX_O^m,\calX_A^m,\pfrak^m,\bfr^m,\bfD^m)_{m\in [M]}, \bfb\right)$.

Remark that a weakly coupled POMDP~\ref{pb:decPOMDP_wc} is a POMDP~\ref{pb:POMDP} with state space $\calX_S = \calX_S^1 \times \cdots \times \calX_S^M$, observation space $\calX_O = \calX_O^1 \times \cdots \times \calX_O^M$, the action space $\calX_A$ defined in~\eqref{eq:problem:def_form_action_space}. 


\begin{rem}\label{rem:problem:def_emission_proba}
	In the definition of POMDP, we could have considered a variant where the observation $O_t$ may depend on $A_{t-1}$ given $S_t$ and the emission probability distribution becomes $\bbP(O_t=o|A_{t-1}=a',S_t=s):= p(o|a',s)$. 
	All the mathematical programming formulations and theoretical results in this paper can be extended to this case. We choose to consider the case above to lighten the notation.
\end{rem}

\begin{rem}\label{rem:problem:decomposablePOMDP}
\citet{bertsimas2016decomposable} consider the notion of \emph{decomposable MDP} as an alternative to weakly coupled MDP.
The only difference is that the action space is generic and does not decompose along the components.
We can define a similar notion of decomposable POMDP.
It turns out that using a transformation similar to the one introduced by \citet[Sec 4.3]{bertsimas2016decomposable}, we can prove that the frameworks of weakly coupled and decomposable POMDPs are equivalent.
Indeed, given a generic action space $\calX_A$, it suffices to define the set of individual action spaces as $\calX_A^m = \calX_A$ for each component $m$ in $[M]$. 
For any $(a^1,\ldots,a^M) \in \calX_A^1 \times \cdots \times \calX_A^M$, we enforce the following linking constraints $a^m = a^{m+1}$ for all $m$ in $[M-1]$. Therefore, the action space can be written $\big\{\bfa \in \calX_A^1 \times \cdots \times \calX_A^M \colon a^m = a^{m+1}, \ \forall m \in [M-1] \big\}$, which has the requested form~\eqref{eq:problem:def_form_action_space}.
In this paper, we choose the weakly coupled POMDP framework because we want to exploit it in our algorithms the explicit structure it presupposes on the action space~\eqref{eq:problem:def_form_action_space}.

\end{rem}

\subsection{Example of application to a maintenance problem}
\label{sub:problem:example}

We now illustrate how the maintenance problem mentioned in the introduction can be casted as a weakly coupled POMDP.
As mentioned in the introduction, Appendix~\ref{app:examples} describes three more examples, and more POMDPs applications can be found for instance in the survey of~\citet{Cassandra2003ASO}.
The system is composed of $M$ components. 
The time discretization $t \in [T]$ corresponds to the different maintenance slots. The decision maker chooses which equipment to maintain at the beginning of each of these slots. 
We model the degradation of component $m$ using a state $S_t^m$, which belongs to a finite state space $\calX_S^m$ and is not observed by the decision maker. 
We assume that there is a \emph{failure state} $s^{m,F}$ in $\calX_S^m$ for each component $m$ in $[M]$, corresponding to its most critical degradation state. 
Component $m$ starts in state $s$ with probability $p^m(s)$. At each time $t$, component $m$ is in state $S_t^m = s$, and it emits an observation $O_t^m=o$ with probability $p^m(o|s)$. 
Then, the decision maker takes an action $\mathbf{A}_t$ in $\{0,1\}^M$ where $A_t^m$ is a binary variable equal to $1$ when component $m$ is maintained. 
At each maintenance slot, the decision maker can maintain at most $K$ components. Hence, we write the action space $\calX_A$ as follows
\begin{align}\label{eq:int_prog_weakPOMDP:problem:pred_maint:def_action_space}
\calX_A = \left\{\bfa=(a^1,\ldots,a^M) \in \{0,1\}^M \colon \sum_{m=1}^M a^m\leq K \right\}.
\end{align}
Therefore, $\calX_A$ contains only one scalar constraint ($q=1$) and satisfies~\eqref{eq:problem:def_form_action_space} by setting $D^m(a) = a$ for every $a \in \{0,1\}$ and $m \in [M]$, and $b = K$.
We assume that each component $m$ evolves independently from state $S_t^m =s$ to state $S_{t+1}^{m}=s'$ with probability $p^m(s'|s,a)$, and the decision maker receives reward $r^m(s, a, s')$. In addition, we assume that when a component is maintained, it behaves like a new one, i.e.,
$p^m(s'|s,1) =p^m(s')$,
for any $s,s'$ in $\calX_S^m$, and the conditional probabilities factorize as~\eqref{eq:problem:factor_emission} and~\eqref{eq:problem:factor_transition}. 
Each component has a maintenance cost $C_R^m$ and a failure cost $C_F^m$ at each component $m$. The individual immediate reward function can be written
$r^m(s,a,s') = -\mathds{1}_{s_F^m}(s') C_F^m - \mathds{1}_{1}(a) C_R^m$,
for any $s,s' \in \calX_S^m$ and $a \in \{0,1\}$.
We assume that the reward decomposes additively as~\eqref{eq:problem:factor_rewards}.
Given a finite horizon $T$, the predictive maintenance problem with capacity constraints consists in finding a policy in $\Deltaml$ that solves~\ref{pb:decPOMDP_wc} with $\calX_A$, $(\pfrak^m)_{m \in[M]}$ and $(\bfr^m)_{m\in [M]}$.

\section{Integer programming for POMDPs}
\label{sec:pomdp}

In this section, we provide an integer formulation that gives an optimal memoryless policy for~\ref{pb:POMDP} as well as valid inequalities. 
We are given a POMDP $(\calX_S, \calX_O, \calX_A, \pfrak, \bfr)$ and a finite horizon $T \in \bbZ_{+}$. We denote by $v_{\rm{ml}}^*$ the optimal value of \ref{pb:POMDP}.  

\subsection{An exact Nonlinear Program}
\label{sub:pomdp:NLP}

We introduce the following nonlinear program (NLP) with a collection of variables\\
$\bfmu = \left((\mu_s^1)_{s}, \left((\mu_{soa}^t)_{s,o,a},(\mu_{sas'}^t)_{s,a,s'}\right)_{t}\right)$, $\bfdelta = \left( \left(\delta_{a|o}^t\right)_{a,o} \right)_{t}$.

\begin{subequations}\label{pb:pomdp:NLP_pomdp}
\begin{alignat}{2}
\max_{\bfmu, \bfdelta}  \enskip & \sum_{t=1}^T \sum_{\substack{s,s' \in \calX_S \\ a \in \calX_A}} r(s,a,s') \mu_{sas'}^t  & \quad & \label{eq:pomdp:NLP_obj_function}\\
\mathrm{s.t.} \enskip
 & \sum_{o \in \calX_O} \mu_{soa}^t = \sum_{s' \in \calX_S} \mu_{sas'}^t &  \quad  \forall s \in \calX_S, a \in \calX_A, t \in [T] \label{eq:pomdp:NLP_consistency_sa}\\
 & \sum_{s \in \calX_S, a \in \calX_A} \mu_{sas'}^t = \sum_{o \in \calX_O, a \in \calX_A} \mu_{s'oa}^{t+1} &  \quad  \forall s' \in \calX_S, t \in [T] \label{eq:pomdp:NLP_consistency_s} \\
 & \mu_s^1 = \sum_{o \in \calX_O, a \in \calX_A} \mu_{soa}^{1} &  \quad  \forall s \in \calX_S \label{eq:pomdp:NLP_consistency_initial} \\
 & \mu_s^1 = p(s) &  \quad  \forall s \in \calX_S \label{eq:pomdp:NLP_initial2} \\
 & \mu_{sas'}^{t} = p(s'|s,a) \sum_{s'' \in \calX_S} \mu_{sas''}^t &  \quad  \forall s,s' \in \calX_S, a \in \calX_A, t \in [T] \label{eq:pomdp:NLP_indep_state}\\
 &\mu_{soa}^t = \delta^t_{a|o} p(o|s) \sum_{o' \in \calX_O,a' \in \calX_A} \mu_{so'a'}^t & \quad \forall s \in \calX_S, o \in \calX_O, a \in \calX_A, t \in [T] \label{eq:pomdp:NLP_indep_action} \\
 & \bfdelta \in \Deltaml, \bfmu \geq 0 \label{eq:pomdp:NLP_constraint_policy} 
\end{alignat}
\end{subequations}
Given a policy $\bfdelta \in \Deltaml$, we say that $\bfmu$ is the vector of \emph{moments} of the probability distribution $\bbP_{\bfdelta}$ induced by $\bfdelta$ when
\begin{subequations}\label{eq:pomdp:def_moments}
	\begin{alignat}{2}
		&\mu_s^1 = \bbP_{\bfdelta}(S_1=s), & & \quad \forall s \in \calX_S \\
		&\mu_{soa}^t = \bbP_{\bfdelta}(S_t=s, O_t=o, A_t=a), & & \quad \forall s \in \calX_S, o \in \calX_O, a \in \calX_A, \forall t \in [T] \\
		&\mu_{sas'}^t = \bbP_{\bfdelta}(S_t=s, A_t=a, S_{t+1}=s'), & & \quad \forall s,s' \in \calX_S, a \in \calX_A, \forall t \in [T]
	\end{alignat}
\end{subequations}
Thanks to the properties of probability distributions, such vector of moments~\eqref{eq:pomdp:def_moments} of $\bbP_{\bfdelta}$ satisfy the constraints of Problem~\eqref{pb:pomdp:NLP_pomdp}. Conversely, given a feasible solution of Problem~\eqref{pb:pomdp:NLP_pomdp}, Theorem~\ref{theo:pomdp:NLP_optimal_solution} ensures that $\bfmu$ is the vector of moments of $\bbP_{\bfdelta}$. 
We denote by $z^*$ the optimal value of Problem~\eqref{pb:pomdp:NLP_pomdp}.

\begin{theo}\label{theo:pomdp:NLP_optimal_solution}
    Let $(\bfmu, \bfdelta)$ be a feasible solution of NLP~\eqref{pb:pomdp:NLP_pomdp}. Then $\bfmu$ is the vector of moments of the probability distribution $\bbP_{\bfdelta}$ induced by $\bfdelta$, and $(\bfmu, \bfdelta)$ is an optimal solution of NLP~\eqref{pb:pomdp:NLP_pomdp} if and only if $\bfdelta$ is an optimal policy of ~\ref{pb:POMDP}. In particular, $v_{\rm{ml}}^* = z^*$.
\end{theo}

\begin{rem}\label{rem:pomdp:with_observation} \emph{Taking into account initial observations.}
	Suppose that the decision maker has access to an initial observation $\ovo$ in $\calX_O$.
	Hence, for any policy $\bfdelta$ in $\Deltaml$ we have $\bbP_{\bfdelta} (O_1=\ovo) = 1.$
	Taking into account the initial observation requires to slightly modify the constraints of Problem~\eqref{pb:pomdp:NLP_pomdp}:
	We replace constraints~\eqref{eq:pomdp:NLP_initial2} and \eqref{eq:pomdp:NLP_indep_action} in Problem~\eqref{pb:pomdp:NLP_pomdp} at time $t=1$ by
	\begin{subequations}
		\begin{alignat}{2}
			&\mu_s^1 = \bbP_{\bfdelta}(S_1=s|O_1=\ovo), & \quad \forall s \in \calX_S, \\
			&\mu_{soa}^1 =  \delta_{a|o}^1 \mathds{1}_{\ovo}(o) \sum_{o' \in \calX_O, a' \in \calX_A}\mu_{so'a'}^1, & \quad \forall s \in \calX_S, o \in \calX_O, a \in \calX_A \label{eq:pomdp:NLP_init_obs_indep},
		\end{alignat}
	\end{subequations}
	where the probability distribution $\bbP_{\bfdelta}(S_1|O_1)$ can be computed using Bayes formula.
	This remark will be useful in Section~\ref{sec:wkpomdp}.
\end{rem}

\subsection{Turning the nonlinear program into an MILP}
\label{sub:pomdp:MILP}

We define the set of \emph{deterministic memoryless policies} $\Deltaml^{\rm{d}}$ as
\begin{equation}\label{eq:pomdp:def_deterministic_policy}
\Deltaml^{\rm{d}} = \bigg\{ \bfdelta \in \Deltaml \colon \delta^t_{a|o} \in \{0,1\},\ \forall o \in \calX_O,\ \forall a \in \calX_A,\ \forall t \in [T] \bigg\}.
\end{equation}
The following proposition states that we can restrict our policy search in~\ref{pb:POMDP} to the set of deterministic memoryless policies.
\begin{prop}{\cite[Proposition 1]{Bagnell2004}}\label{prop:pomdp:det_policies}
	There always exists an optimal policy for~\ref{pb:POMDP} that is deterministic, i.e.,
	\begin{equation}\label{eq:pomdp:prop:det_policies}
    	\max_{\bfdelta \in \Deltaml} \bbE_{\bfdelta} \bigg[ \sum_{t=1}^{T} r(S_t,A_t, S_{t+1}) \bigg] = \max_{\bfdelta \in \Deltaml^{\rm{d}}} \bbE_{\bfdelta} \bigg[ \sum_{t=1}^{T} r(S_t,A_t, S_{t+1})\bigg].
	\end{equation}
\end{prop}

Theorem~\ref{theo:pomdp:NLP_optimal_solution} ensure that~\ref{pb:POMDP} and Problem~\eqref{pb:pomdp:NLP_pomdp} are equivalent, and in particular admit the same optimal solution in terms of $\bfdelta$. However, Problem~\eqref{pb:pomdp:NLP_pomdp} is hard to solve due to the nonlinear constraints~\eqref{eq:pomdp:NLP_indep_action}.
By Proposition~\ref{eq:pomdp:prop:det_policies}, we can add the integrality constraints of $\Deltaml^{\rm{d}}$ in~\eqref{pb:pomdp:NLP_pomdp}, and, by a classical result in integer programming, we can turn Problem~\eqref{pb:pomdp:NLP_pomdp} into an equivalent MILP: It suffices to replace constraint~\eqref{eq:pomdp:NLP_indep_action} by the following McCormick inequalities~\citep{Mccormick1976}.
\begin{subequations}\label{eq:pomdp:McCormick_linearization}
	\begin{alignat}{2}
		&\mu_{soa}^t \leq p(o|s) \sum_{o' \in \calX_O,a' \in \calX_A} \mu_{so'a'}^t & \quad \forall s \in \calX_S, o \in \calX_O, a \in \calX_A, t \in [T] \label{eq:pomdp:MILP_McCormick_1} \\
 		&\mu_{soa}^t \leq \delta^t_{a|o} & \quad  \forall s \in \calX_S, o \in \calX_O, a \in \calX_A, t \in [T] \label{eq:pomdp:MILP_McCormick_2} \\
 		&\mu_{soa}^t \geq p(o|s)\sum_{o' \in \calX_O,a' \in \calX_A} \mu_{so'a'}^t  + \delta^t_{a|o} - 1 & \quad \forall s \in \calX_S, o \in \calX_O, a \in \calX_A, t \in [T]. \label{eq:pomdp:MILP_McCormick_3}
	\end{alignat}
\end{subequations}
For convenience, we denote by $\mathrm{McCormick}\left(\bfmu, \bfdelta \right)$ the set of McCormick linear inequalities~\eqref{eq:pomdp:McCormick_linearization}. 
Thus, by using McCormick's linearization on constraints~\eqref{eq:pomdp:NLP_indep_action}, we get that~\ref{pb:POMDP} is equivalent to the following MILP: 
\begin{equation}
\begin{aligned}\label{pb:pomdp:MILP_pomdp}
\max_{\bfmu, \bfdelta}  \enskip & \sum_{t=1}^T \sum_{\substack{s,s' \in \calX_S \\ a \in \calX_A}} r(s,a,s') \mu_{sas'}^t  & \quad & \\
\mathrm{s.t.} \enskip
 & \bfmu \ \mathrm{satisfies}~\eqref{eq:pomdp:NLP_consistency_sa}-\eqref{eq:pomdp:NLP_indep_state} \\
 & \mathrm{McCormick}\big(\bfmu,\bfdelta\big) \\
 & \bfdelta \in \Deltaml^{\rm{d}},  \bfmu \geq 0. 
\end{aligned}
\end{equation}

For convenience, given a POMDP $\left(\calX_S,\calX_O,\calX_A,\pfrak,\bfr \right)$, we define respectively the feasible sets of Problem~\eqref{pb:pomdp:NLP_pomdp} and MILP~\eqref{pb:pomdp:MILP_pomdp} as $\calQ\left(T,\calX_S,\calX_O,\calX_A,\pfrak \right)$ and $\calQ^{\mathrm{d}}\left(T,\calX_S,\calX_O,\calX_A,\pfrak \right).$
We write respectively $\calQ$ and $\calQ^{\mathrm{d}}$ when $\left(T,\calX_S,\calX_O,\calX_A,\pfrak\right)$ is clear from the context.

\subsection{Valid inequalities}
\label{sub:pomdp:valid_cuts}

Before introducing our valid inequalities, we start by explaining why the linear relaxation of MILP~\eqref{pb:pomdp:MILP_pomdp} is not sufficient to define a feasible solution of Problem~\eqref{pb:pomdp:NLP_pomdp}.
It turns out that given a feasible solution $(\bfmu,\bfdelta)$ of the linear relaxation of MILP~\eqref{pb:pomdp:MILP_pomdp}, the vector $\bfmu$ is not necessarily the vector of moments of the probability distribution $\bbP_{\bfdelta}$ induced by $\bfdelta$. Indeed, when the coordinates of the vector $\bfdelta$ are continuous variables, the McCormick's constraints~\eqref{eq:pomdp:McCormick_linearization} are, in general, no longer equivalent to bilinear constraints~\eqref{eq:pomdp:NLP_indep_action}. Then, $(\bfmu,\bfdelta)$ is not necessarily a feasible solution of Problem~\eqref{pb:pomdp:NLP_pomdp} anymore, which implies that $\bfmu$ is not necessarily the vector of moments of the probability distribution $\bbP_{\bfdelta}$.
Intuitively, it means that the feasible set of the linear relaxation of MILP~\eqref{pb:pomdp:MILP_pomdp} is too large. Actually, we can reduce the feasible set of the linear relaxation of MILP~\eqref{pb:pomdp:MILP_pomdp} by adding valid inequalities.
To do so, we introduce new variables $\left((\mu_{s'a'soa}^t)_{s',a',s,o,a}\right)_{t}$ and the inequalities
\begin{subequations}\label{eq:pomdp:Valid_cuts_pomdp}
    \begin{alignat}{2}
    &\sum_{s'\in \calX_S, a' \in \calX_A} \mu_{s'a'soa}^t = \mu_{soa}^{t}, \quad &\forall s \in \calX_S, o \in \calX_O, a \in \calX_A, \label{eq:pomdp:Valid_cuts_pomdp_consistency1}\\
    &\sum_{a \in \calX_A} \mu_{s'a'soa}^t = p(o|s)\mu_{s'a's}^{t-1}, \quad &\forall s',s \in \calX_S, o \in \calX_O, a' \in \calX_A, \label{eq:pomdp:Valid_cuts_pomdp_consistency2}\\
    &\mu_{s'a'soa}^t = p(s|s',a',o)\sum_{\ovs \in \calX_S} \mu_{s'a'\ovs oa}^t, \quad  &\forall s',s \in \calX_S, o \in \calX_O, a',a \in \calX_A, \label{eq:pomdp:Valid_cuts_pomdp_main}
    \end{alignat}
\end{subequations} 
\noindent where we use the constants $$p(s|s',a',o) = \bbP(S_t=s | S_{t-1}=s',A_{t-1}=a',O_t=o),$$
for any $s,s' \in \calX_S$, $a'\in \calX_A$ and $o \in \calX_O$.
Note that $p(s|s',a',o')$ does not depend on the policy $\bfdelta$ and can be easily computed during a preprocessing using Bayes rules. Therefore, constraints in~\eqref{eq:pomdp:Valid_cuts_pomdp} are linear.
\begin{prop}\label{prop:pomdp:valid_cuts_pomdp}
    Inequalities \eqref{eq:pomdp:Valid_cuts_pomdp} are valid for MILP~\eqref{pb:pomdp:MILP_pomdp}, and there exists a solution $\bfmu$ of the linear relaxation of~\eqref{pb:pomdp:MILP_pomdp} that does not satisfy constraints~\eqref{eq:pomdp:Valid_cuts_pomdp}. 
\end{prop}


The MILP formulation obtained by adding inequalities~\eqref{eq:pomdp:Valid_cuts_pomdp} in MILP~\eqref{pb:pomdp:MILP_pomdp} is an extended formulation, and has many more constraints than the initial MILP~\eqref{pb:pomdp:MILP_pomdp}. 
Its linear relaxation therefore takes longer to solve.
Inequalities~\eqref{eq:pomdp:Valid_cuts_pomdp} strengthen the linear relaxation, and numerical experiments in Section \ref{sec:num} show that these inequalities enable to speed up the resolution of MILP~\eqref{pb:pomdp:MILP_pomdp}.

\paragraph{Probabilistic interpretation.}
	Given a feasible solution $(\bfmu,\bfdelta)$ of the linear relaxation of~\eqref{pb:pomdp:MILP_pomdp}, $\bfmu$ can still be interpreted as the vector of moments of a probability distribution $\bbQ_{\bfmu}$ over 
	$\left(\calX_S \times \calX_O \times \calX_A\right)^T \times \calX_S$.
	However, as it has been mentioned above, the vector $\bfmu$ does not necessarily correspond to the vector of moments of $\bbP_{\bfdelta}$, which is due to the fact that $(\bfmu,\bfdelta)$ does not necessarily satisfy the nonlinear constraints~\eqref{eq:pomdp:NLP_indep_action}.
	Besides, constraints~\eqref{eq:pomdp:NLP_indep_action} is equivalent to the property that, 
	\begin{equation}\label{eq:pomdp:strongIndep}
	\text{according to $\bbQ_{\bfmu}$, action $A_t$ is independent from state $S_t$ given observation $O_t$.} 
	\end{equation}
	Hence, given a feasible solution $(\bfmu,\bfdelta)$ of the linear relaxation of MILP~\eqref{pb:pomdp:MILP_pomdp}, the distribution $\bbQ_{\bfmu}$ does not necessarily satisfy the conditional independences~\eqref{eq:pomdp:strongIndep}.
	Remark that~\eqref{eq:pomdp:strongIndep} implies the weaker result that,
	\begin{equation}\label{eq:pomdp:weakIndep}
	    \text{according to $\bbQ_{\bfmu}$, $A_t$ is independent from $S_t$ given $O_t$, $A_{t-1}$ and $S_{t-1}$.}
	\end{equation}
	Proposition~\ref{prop:pomdp:valid_cuts_pomdp} says that the independences in~\eqref{eq:pomdp:weakIndep} are not satisfied in general by a feasible solution $(\bfmu,\bfdelta)$ of the linear relaxation of MILP~\eqref{pb:pomdp:MILP_pomdp}, but that we can enforce them using linear inequalities~\eqref{eq:pomdp:Valid_cuts_pomdp} on $(\bfmu,\bfdelta)$ in an extended formulation.

\subsection{Strengths of the relaxations}
\label{sub:pomdp:inf_relax}

When the decision maker directly observes the state of the system, the POMDP problem becomes a MDP problem and the resulting optimization problem is called the \emph{MDP approximation} \citep{Hauskrecht2000}. 
It turns out that the linear relaxation of MILP~\eqref{pb:pomdp:MILP_pomdp} is related to the MDP approximation.
Given a collection of variables $\bfmu = \left((\mu_s^1)_{s},(\mu_{sas'}^t)_{s,a,s'})_{t}\right)$ and the following linear program which is known to solve exactly a MDP (e.g. \citet{Epenoux1963}).
\begin{subequations}\label{pb:pomdp:LP_MDP}
\begin{alignat}{2}
\max_{\bfmu}  \enskip & \sum_{t=1}^T \sum_{\substack{s,s' \in \calX_S \\ a \in \calX_A}} r(s,a,s')\mu_{sas'}^t  & \quad &\\
\mathrm{s.t.} \enskip 
& \mu_{s}^1 = \sum_{a' \in \calX_A, s' \in \calX_S} \mu_{sa's'}^{1} & \forall s \in \calX_S \label{eq:pomdp:LP_MDP_consistent_action_initial}\\
& \sum_{s'\in \calX_S,a' \in \calX_A} \mu_{s'a's}^t = \sum_{a' \in \calX_A, s' \in \calX_S} \mu_{sa's'}^{t+1} & \forall s \in \calX_S, t \in [T] \label{eq:pomdp:LP_MDP_consistent_action}\\
& \mu_{s}^1 = p(s) & \forall s \in \calX_S \label{eq:pomdp:LP_MDP_initial} \\
& \mu_{sas'}^{t} = p(s'|s,a) \sum_{s'' \in \calX_S }\mu_{sas''}^t &  \forall s \in \calX_S, a \in \calX_A, t \in [T] \label{eq:pomdp:LP_MDP_consistent_state}
\end{alignat}
\end{subequations}
In Problem~\eqref{pb:pomdp:LP_MDP}, the variables $(\mu_s^1)_s$ and $(\mu_{sas'}^t)_{sas'}$ respectively represent the probability distribution of $S_1$ and $\left(S_t,A_t,S_{t+1}\right)$ for any $t$ in $[T]$.
Theorem~\ref{theo:pomdp:MDP_approx_equivalence} below states that the linear relaxation of MILP~\eqref{pb:pomdp:MILP_pomdp} is equivalent to the MDP approximation of~\ref{pb:POMDP}. 
We introduce the following quantities:
\begin{enumerate}[label={--}]
	\item $z_{\rm{R}}^*$: the optimal value of the linear relaxation of MILP~\eqref{pb:pomdp:MILP_pomdp}.
	\item $z_{\rm{R}^{\rm{c}}}^*$: the optimal value of the linear relaxation of MILP~\eqref{pb:pomdp:MILP_pomdp} with inequalities~\eqref{eq:pomdp:Valid_cuts_pomdp}.
	\item $v_{\rm{his}}^*$: the optimal value of~\ref{pb:POMDP_perfectRecall}.
	\item $v_{\rm{MDP}}^*$: the optimal value of linear program~\eqref{pb:pomdp:LP_MDP}, which is the optimal value of the MDP approximation.
\end{enumerate}
\begin{theo}\label{theo:pomdp:MDP_approx_equivalence}
	Let $(\bfmu,\bfdelta)$ be feasible solution of the linear relaxation of MILP~\eqref{pb:pomdp:MILP_pomdp}. Then $(\bfmu,\bfdelta)$ is an optimal solution of the linear relaxation of MILP~\eqref{pb:pomdp:MILP_pomdp} if and only if $\bfmu$ is an optimal solution of linear program~\eqref{pb:pomdp:LP_MDP}. In particular, $z_{\rm{R}}^* = v_{\rm{MDP}}^*$. 
    In addition, the following inequalities hold:
    \begin{equation}\label{eq:pomdp:inequality_information}
    	z^* \leq v_{\rm{his}}^* \leq z_{\rm{R}^{\rm{c}}}^* \leq z_{\rm{R}}^*.
    \end{equation}
\end{theo}

Inequality~\eqref{eq:pomdp:inequality_information} ensures that by solving MILP~\eqref{pb:pomdp:MILP_pomdp}, we obtain an integrality gap $z_{\rm{R}}^*-z^*$ that bounds the approximation error $v_{\rm{his}}^* - z^*$ due to the choice of a memoryless policy instead of a policy that depends on all history of observations and actions. In addition, Theorem~\ref{theo:pomdp:MDP_approx_equivalence} ensures that the integrality gap $z_{\rm{R}^{\rm{c}}}^*-z^*$ obtained using valid inequalities~\eqref{eq:pomdp:Valid_cuts_pomdp} gives a tighter bound on the approximation error.

\section{Integer programming for weakly coupled POMDPs}
\label{sec:wkpomdp}

We now focus on the problem~\ref{pb:decPOMDP_wc} of finding a high/maximum expected reward policy for a weakly coupled POMDP $\left((\calX_S^m,\calX_O^m,\calX_A^m,\pfrak^m,\bfr^m,\bfD^m)_{m\in[M]},\bfb\right)$.
Based on the results of Section~\ref{sec:pomdp}, a naive approach is to use the mathematical programs of the previous section
on the POMDP $\left(\calX_S,\calX_O,\calX_A',\pfrak',\bfr' \right),$ where $\calX_S = \calX_S^1 \times \cdots \times \calX_S^M,$ $\calX_O=\calX_O^1 \times \cdots \times \calX_O^M,$ and, $\calX_A'=\calX_A^1 \times \cdots \times \calX_A^M,$ $\pfrak',$ and $\bfr'$ are respectively defined by~\eqref{eq:problem:def_form_action_space}-\eqref{eq:problem:factor_rewards}.
More precisely, NLP~\eqref{pb:pomdp:NLP_pomdp} becomes
\begin{subequations}\label{pb:wkpomdp:POMDP_NLP}
 \begin{alignat}{2}
 \max_{\bfmu,\bfdelta} \enskip & \sum_{t=1}^T \sum_{\substack{\bfs,\bfs' \in \calX_S \\ \bfa \in \calX_A}} \mu_{\bfs \bfa \bfs'}^{t}\sum_{m=1}^M r^m(s,a,s') & \quad &\\
\mathrm{s.t.} \enskip 
 & \left(\bfmu,\bfdelta\right) \in \calQ \left(T,\calX_S, \calX_O, \calX_A',\pfrak'\right) & \label{eq:wkpomdp:POMDP_NLP:per_pomdp} \\
 & \left( \sum_{m=1}^M \bfD^m(a^m) - \bfb \right)\delta_{\bfa|\bfo}^t  \leq 0 & \forall \bfo \in \calX_O, \bfa \in \calX_A', t \in [T] \label{eq:wkpomdp:POMDP_NLP:link}
 \end{alignat}
\end{subequations}
where constraints~\eqref{eq:wkpomdp:POMDP_NLP:link} model the linking constraints $\sum_{m=1}^M \bfD^m(A_t^m) \leq \bfb$, which are not included in $\calQ \left(T,\calX_S, \calX_O, \calX_A',\pfrak'\right)$.
Indeed, given a solution $(\bfmu,\bfdelta)$ of NLP~\eqref{pb:wkpomdp:POMDP_NLP}, if $\bfa \in \calX_A' \backslash \calX_A$, then $\delta_{\bfa|\bfo} = 0$ for all $t$ in $[T]$ and $\bfo$ in $\calX_O$. This results ensures that $\bbP_{\bfdelta} \left(\left\{ \bfA_t \in \calX_A, \enskip \forall t \in [T] \right\} \right) = 1$, which shows that the linking constraint is satisfied almost surely.

NLP~\eqref{pb:wkpomdp:POMDP_NLP} is intractable for at least three reasons.
(A) The number of variables required to encode a policy $\bfdelta$ and its moments $\bfmu$ is exponential (the vector $(\mu_{\bfs \bfa \bfs'})_{\substack{\bfs,\bfs' \in \calX_S \\ a \in \calX_A}}$ has $|\calX_A| \prod_{m=1}^M |\calX_S^m| |\calX_O^m|$ coordinates) and the number of constraints required to ensure that $\bfmu$ corresponds to the moment of the distribution induced by $\bfdelta$ is exponential.
(B) The latter constraints are nonlinear.
And (C) the number of inequalities required to ensure that the linking constraints are satisfied is exponential.
If the second difficulty can be addressed using the approach developed in Section~\ref{sub:pomdp:MILP}, the two others are specific to weakly coupled POMDP.

\subsection{An approximate integer program}
\label{sub:wkpomdp:approx}


Consider the variables
$\bftau^m = \Big( (\tau_s^{1,m})_{s}, (\tau_{sas'}^{t,m})_{s,a,s'}, (\tau_{soa}^{t,m})_{s,o,a}, (\tau_a^{t,m})_{a} \Big)_{t \in [T]},$ $\bfdelta^m =(\delta^{t,m})_{t \in [T]}$ and the following MILP,
\begin{subequations}\label{pb:wkpomdp:decPOMDP_MILP}
 \begin{alignat}{2}
 \max_{\bftau,\bfdelta} \enskip & \sum_{t=1}^T \sum_{m=1}^M \sum_{\substack{s,s' \in \calX_S^m \\ a \in \calX_A^m}} r^m(s,a,s')\tau_{sas'}^{t,m} & \quad &\\
\mathrm{s.t.} \enskip 
 & \left(\bftau^m,\bfdelta^m\right) \in \calQ^{\mathrm{d}}\left(T,\calX_S^m, \calX_O^m, \calX_A^m,\pfrak^m\right)   & \forall m \in [M] \label{eq:wkpomdp:decPOMDP_MILP_per_pomdp} \\
 &\sum_{s \in \calX_S^m, o \in \calX_A^m}\tau_{soa}^{t,m} = \tau_a^{t,m} & \forall a \in \calX_A^m, m \in [M], t \in [T] \label{eq:wkpomdp:decPOMDP_MILP_consistent_action} \\
 &\sum_{m =1}^M\sum_{a \in \calX_A^m} \bfD^m(a) \tau_a^{t,m} \leq \bfb & \forall t \in [T] \label{eq:wkpomdp:decPOMDP_MILP_linking_cons}
 \end{alignat}
\end{subequations}

In order to obtain~\eqref{pb:wkpomdp:decPOMDP_MILP}, we modify NLP~\eqref{pb:wkpomdp:POMDP_NLP} in three ways.
\begin{itemize}
	\item[\textup{(A)}] We reduce the number of variables and constraints required to encode a policy and its moments by considering ``local'' variables $(\bftau^m,\bfdelta^m)$ in $\calQ\left(T,\calX_S^m, \calX_O^m, \calX_A^m,\pfrak^m\right)$ for each component $m$ instead of $(\bfmu,\bfdelta)$.
	\item[\textup{(B)}] We consider deterministic policies $\bfdelta^m \in \Delta^{\mathrm{d},m}$ and linearize constraints using McCormick inequalities: We replace $\calQ$ by $\calQ^{\mathrm{d}}$.
	\item[\textup{(C)}] We replace the almost sure linking constraint $\sum_{m=1}^M \bfD^m(A_t^m) \leq \bfb,$ by the constraint in expectation $\sum_{m =1}^M \bbE_{\bfdelta^m} \left[ \bfD^m(A_t^m) \right] \leq \bfb$, which enables to reduce the number of inequalities~\eqref{eq:wkpomdp:POMDP_NLP:link} required to encode it.
\end{itemize}
The reader is already familiar with \textup{(B)}, which we used in Section~\ref{sub:pomdp:MILP} to turn the NLP~\eqref{pb:pomdp:NLP_pomdp} into the MILP~\eqref{pb:pomdp:MILP_pomdp}. Let us now focus \textup{(A)} and \textup{(C)}.

(A) Moving from $\bfmu$ and $\bfdelta$ to $\bftau^m$ and $\bfdelta^m =(\delta^{t,m})_{t \in [T]}$ enables to obtain a MILP with $O\bigg(\sum_{m=1}\vert \calX_S^m\vert \vert \calX_A^m \vert\Big(\vert \calX_O^m\vert + \vert \calX_S^m\vert \Big)\bigg)$ variables and constraints when NLP~\eqref{pb:wkpomdp:POMDP_NLP} had $O\bigg(\prod_{m=1}\vert \calX_S^m\vert \vert \calX_A^m \vert\Big(\prod_{m=1}^M\vert \calX_O^m\vert + \prod_{m=1}^M \vert \calX_S^m\vert \Big)\bigg)$ variables and constraints.
Furthermore, by Theorem~\ref{theo:pomdp:NLP_optimal_solution}, constraints~\eqref{eq:wkpomdp:decPOMDP_MILP_per_pomdp} and ~\eqref{eq:wkpomdp:decPOMDP_MILP_consistent_action} ensure that the variables $\bftau^m$ can still be interpreted as the vector of moments of the probability distribution $\bbP_{\bfdelta^m}$ induced by the deterministic policy $\bfdelta^m$ for each component $m$ in $[M].$
However, given a feasible solution $(\bftau^m,\bfdelta^m)_{m\in [M]}$ there is no guarantee that there exists a policy $\bfdelta$ of~\ref{pb:decPOMDP_wc} such that the variables $(\bftau^m)_{m \in [M]}$ represent the moments of the distribution $\bbP_{\bfdelta}$ on the whole system.
It is the reason why, following the conventions in the graphical model literature \citep{wainwright2008graphical}, we denote by $\bftau$ the approximate vector of moments (or pseudomarginals) instead of $\bfmu$, which we reserve for exact vector of moments.

(C) If we leverage the interpretation of $(\bfdelta^m)$ as ``local'' policies for each component, the linking constraints $\sum_{m=1}^M \bfD^m(A_t^m) \leq \bfb$ are satisfied almost surely by the collection of constraints
\begin{equation}\label{eq:wkpomdp:decPOMDP_link_almost_surely}
\sum_{m =1}^M \delta^{t,m}_{a^m|o^m}\bfD^m(a^m) \leq \bfb, \quad \forall \bfo \in \calX_O, \enskip \bfa \in \calX_A', \enskip t \in [T].
\end{equation}
Indeed, since we use deterministic ``local'' policies for each component, the latter constraint ensures that $\sum_{m \colon \delta^{t,m}_{a^m|o^m} = 1} \bfD^m(a^m) \leq \bfb$.
However, the collection of constraints~\eqref{eq:wkpomdp:decPOMDP_link_almost_surely} has still an exponential number of elements.
When we replace the almost sure constraint $\sum_{m=1}^M \bfD^m(A_t^m) \leq \bfb$ by the (weaker) constraint in expectation $\bbE_{\bfdelta}\left[\sum_{m=1}^M \bfD^m(A_t^m)\right] \leq \bfb$, the linearity of expectation enables to decouple the different components and obtain the family of constraint~\eqref{pb:pomdp:MILP_pomdp}, which has a polynomial number of elements.

In addition to its tractability, MILP~\eqref{pb:wkpomdp:decPOMDP_MILP} has the advantage when it is solved as a subroutine in the algorithms of Sections~\ref{sub:wkpomdp:first_policy} and~\ref{sub:wkpomdp:matheuristic} that compute practically efficient policies for large scale weakly coupled POMDP.
But its optimal value is neither a lower bound nor an upper bound on $v_{\rm{ml}}^*$ (see Appendix~\ref{app:counter_example} for examples of this fact). 
We will see in the next section that it shares lower and upper bounds with~\ref{pb:decPOMDP_wc}.

As it was the case for MILP~\eqref{pb:pomdp:MILP_pomdp}, we can strengthen MILP~\eqref{pb:wkpomdp:decPOMDP_MILP} with the valid inequalities of Section~\ref{sub:pomdp:valid_cuts}, which are now reformulated as follows, and whose analysis is deferred to Section~\ref{sub:wkpomdp:valid_ineq}.
\begin{subequations}\label{eq:wkpomdp:dec_Valid_cuts}
    \begin{alignat}{2}
    &\sum_{s'\in \calX_S^m, a' \in \calX_A^m} \tau_{s'a'soa}^{t,m} = \tau_{soa}^{t,m},& \forall s \in \calX_S^m, o \in \calX_O^m, a \in \calX_A^m, \label{eq:wkpomdp:dec_Valid_cuts_consistency1} \\
    &\sum_{a \in \calX_A^m} \tau_{s'a'soa}^{t,m} = p^m(o|s)p^m(s|s',a')\tau_{s'a'}^{t-1,m}, & \forall s,s' \in \calX_S^m, o \in \calX_O^m, a'  \in \calX_A^m, \label{eq:wkpomdp:dec_Valid_cuts_consistency2} \\
    &\tau_{s'a'soa}^{t,m} = p^m(s|s',a',o)\sum_{\overline{s} \in \calX_S^m} \tau_{s'a'\overline{s} oa}^{t,m}, \quad \quad &\forall s,s' \in \calX_S^m, o \in \calX_O^m,a,a'  \in \calX_A^m. \label{eq:wkpomdp:dec_Valid_cuts_main}
    \end{alignat}
\end{subequations}

\subsection{Strengths of the linear relaxation}
\label{sub:wkpomdp:strengths_linear_relaxation}

While in Section~\ref{sub:pomdp:inf_relax} we showed that the linear relaxation of MILP~\eqref{pb:pomdp:MILP_pomdp} is equivalent to the MDP approximation, one may ask the question: How do we relate the linear relaxation of MILP~\eqref{pb:wkpomdp:decPOMDP_MILP} with the MDP approximation of a weakly coupled POMDP?
As stated the theorem below, we are able to link the value of the linear relaxation of MILP~\eqref{pb:wkpomdp:decPOMDP_MILP} (with and without valid inequalities~\eqref{eq:wkpomdp:dec_Valid_cuts}) with the optimal value $v_{\rmml}^*$ and $v_{\rmhis}^*$.
We denote respectively by $z_{\rmIP}$, $z_{\rmRc}$ and $z_{\rmR}$ the optimal values of MILP~\eqref{pb:wkpomdp:decPOMDP_MILP} and its linear relaxation with and without valid inequalities~\eqref{eq:wkpomdp:dec_Valid_cuts}.

\begin{theo}\label{theo:wkpomdp:ineq_decPOMDP}
	The linear relaxation of MILP~\eqref{pb:wkpomdp:decPOMDP_MILP} is a relaxation of the MPD approximation of the weakly coupled POMDP. Furthermore, the inequalities $v_{\rmMDP}^* \leq z_{\rmR}$ and $v_{\rmhis}^* \leq z_{\rmRc} \leq z_{\rmR}$ hold.
\end{theo}

It turns out the linear relaxation of MILP~\eqref{pb:wkpomdp:decPOMDP_MILP} is equivalent to the fluid formulation of \citet{bertsimas2016decomposable}, which is a relaxation of the MDP approximation of a weakly coupled POMDP. We show this result in Appendix~\ref{app:linksBertsimas}.

\subsection{An upper bound and a lower bound}
\label{sub:wkpomdp:lb_ub}

In this section we introduce bounds $z_{\rmLB}$, $z_{\rmUB}$, and $z_{\rmLR}$ that enable to quantify the distance between the optimal values $v_{\rm{ml}}^*$ and $z_{\rmIP}$ of~\ref{pb:decPOMDP_wc} and MILP~\eqref{pb:wkpomdp:decPOMDP_MILP}. More precisely, we prove 
\begin{equation*}
\begin{array}{rcl}
  z_{\rm{LB}} \leq& v_{\rm{ml}}^* &\leq z_{\rm{UB}}  \\
  z_{\rm{LB}} \leq&  z_{\rmIP} &\leq z_{\rm{UB}} 
\end{array} \quad \text{and} \quad z_{\rm{UB}} \leq z_{\rm{LR}} \leq z_{\rmRc} \leq z_{\rmR}
\end{equation*}
where $z_{\rm{R}}$ and $z_{\rm{R}^{\rm{c}}}$ denote the optimal values of the linear relaxation of MILP~\eqref{pb:wkpomdp:decPOMDP_MILP} and the linear relaxation of MILP~\eqref{pb:wkpomdp:decPOMDP_MILP} with valid inequalities~\eqref{eq:wkpomdp:dec_Valid_cuts}, respectively.
Bounds $z_{\rmLB}$, $z_{\rmUB}$, and $z_{\rmLR}$ are defined as optimal values of mathematical programs obtained by playing with approximations~(A), (B), and~(C). 
Since $z_{\rmLB}$ and~$z_{\rmUB}$ are hard to compute, their interest is mainly theoretical: They enable to bound the difference between $v_{\rm{ml}}^*$ and~$z_{\rmIP}$.
We also introduce the tractable upper bound~$z_{\rmLR}$ to evaluate the quality of the policy we use in the numerical experiments.

In this section, we need to compare MILP formulations that do not share the same set of variables. We therefore say that a problem P is a \emph{relaxation} of problem P' when given a feasible solution of P' we can build a feasible solution of P with the same value.

\subsubsection{The lower bound from decomposable policies}

A policy $\bfdelta$ in $\Delta$, defined on $\calX_A$, is \emph{decomposable} if there exist policies $\bfdelta^m$ in $\Delta^m$, defined on $\calX_A^m$, such that $\bfdelta = \prod_{m=1}^M \bfdelta^m.$
The main advantage of decomposable policies is that they decouple components:
given such a $\bfdelta$, the random variables of each component are independent under $\bbP_{\bfdelta}$.
We can therefore compute exact probabilities using local moments $\bftau^m$. 
An optimal policy can then be computed using the following MILP, whose optimal value we use as lower bound $z_{\rm{LB}}$.


\begin{subequations}\label{pb:wkpomdp:decPOMDP_MILP_LB}
 \begin{alignat}{2}
 z_{\rm{LB}}:= \max_{\bftau,\bfdelta} \enskip & \sum_{t=1}^T \sum_{m=1}^M \sum_{\substack{s,s' \in \calX_S^m \\ a \in \calX_A^m}} r_m(s,a,s')\tau_{sas'}^{t,m} & \quad &\\
\mathrm{s.t.} \enskip 
 &\left(\bftau^m,\bfdelta^m\right) \in \calQ^{\mathrm{d}}\left(T, \calX_S^m, \calX_O^m, \calX_A^m,\pfrak^m\right)  & \forall m \in [M] \label{eq:wkpomdp:MILP_per_POMDP_Lb} \\
 &\sum_{m =1}^M \sum_{a \in \calX_A^m} \bfD^m(a)\delta^{t,m}_{a|o^m} \leq \bfb & \forall \bfo \in \calX_O, t \in [T] \label{eq:wkpomdp:decPOMDP_link_constraints_Lb}
 \end{alignat}
\end{subequations}
\begin{theo}\label{theo:wkpomdp:lower_bound_MILP}
	Let $(\bftau^m,\bfdelta^m)$ be an optimal solution of~\eqref{pb:wkpomdp:decPOMDP_MILP_LB}. Then $\bfdelta = \prod_{m=1}^M \bfdelta^m$ is an optimal deterministic decomposable policy for~\ref{pb:decPOMDP_wc}.

	Furthermore, MILP~\eqref{pb:wkpomdp:decPOMDP_MILP} is also a relaxation of MILP~\eqref{pb:wkpomdp:decPOMDP_MILP_LB}. In particular, $z_{\rmLB} \leq v_{\rm{ml}}^*$ and $z_{\rm{LB}} \leq z_{\rmIP}.$
\end{theo}
Remark that MILP~\eqref{pb:wkpomdp:decPOMDP_MILP_LB} is obtained from NLP~\eqref{pb:wkpomdp:POMDP_NLP} by using modifications~(A) and (B), but not (C).
Indeed, the almost sure constraint~\eqref{eq:wkpomdp:decPOMDP_link_constraints_Lb} is required to ensure that the decomposable policy we obtain is feasible.
Since there is an exponential number of such constraints, a line generation approach is required to solve this problem in practice.

\subsubsection{An upper bound through a nonlinear formulation}


At first sight, it may seem that MILP~\eqref{pb:wkpomdp:decPOMDP_MILP} is a relaxation of NLP~\eqref{pb:wkpomdp:POMDP_NLP}.
Indeed, modification (A) relaxes the problem: Given a solution $(\bfmu,\bfdelta)$ of NLP~\eqref{pb:wkpomdp:POMDP_NLP}, and defining $\bftau^m$ and $\bfdelta^m$ respectively as the marginal moments and marginal policy as follows
\begin{align*}
	&\left(\tau_s^1,\tau_{soa}^{t,m},\tau_{sas'}^{t,m}\right) = \left( \sum_{s^{-m}} \mu_{\bfs}^1,\sum_{\substack{s^{-m},o^{-m} \\ a^{-m}}} \mu_{\bfs\bfo\bfa}^t, \sum_{\substack{s^{-m},a^{-m} \\ s'^{-m}}} \mu_{\bfs\bfa\bfs'}^t \right), & \delta_{a|o}^{t,m} = \sum_{\substack{s^{-m},o^{-m} \\ a^{-m}}} \delta^t_{\bfa|\bfo} \prod_{m' \neq m} p^{m'}(o^{m'}|s^{m'})\tau_{s^{m'}}^{t,m'},
\end{align*}
where the sum over $s^{-m}$, $o^{-m}$ and $a^{-m}$ indicates respectively the sum over all the vector in $\prod_{m'\neq m} \calX_S^{m'}$, $\prod_{m'\neq m} \calX_O^{m'}$ and $\prod_{m'\neq m} \calX_A^{m'}$,
we obtain a solution that satisfies constraints of $\calQ\left(T,\calX_S^m,\calX_O^m,\calX_A^m,\pfrak^m\right)$.
And modification (C) further relaxes the problem by turning an almost sure constraint into a constraint in expectation.
However, if there is an optimal policy $\bfdelta$ that is deterministic, the moments $\bfdelta^m$ of such a policy defined above are not necessarily deterministic anymore.
This is the reason why MILP~\eqref{pb:wkpomdp:decPOMDP_MILP}, which forces the $\bfdelta^m$ to be integer, is not a relaxation of NLP~\eqref{pb:wkpomdp:POMDP_NLP}.
But the following non-linear program,

\begin{subequations}\label{pb:wkpomdp:decPOMDP_NLP}
 \begin{alignat}{2}
z_{\rmUB}:=\max_{\bftau,\bfdelta} \enskip & \sum_{t=1}^T \sum_{m=1}^M \sum_{\substack{s,s' \in \calX_S^m \\ a \in \calX_A^m}} r_m(s,a,s')\tau_{sas'}^{t,m} & \quad &\\
\mathrm{s.t.} \enskip 
 & \left(\bftau^m,\bfdelta^m\right) \in \calQ\left(T,\calX_S^m,\calX_O^m,\calX_A^m,\pfrak^m \right) & \quad\forall m \in [M] \label{eq:wkpomdp:NLP_per_POMDP_v1} \\
  &\sum_{s \in \calX_S^m, o \in \calX_A^m}\tau_{soa}^{t,m} = \tau_a^{t,m} & \forall a \in \calX_A^m, m \in [M], t \in [T] \label{eq:wkpomdp:decPOMDP_NLP_consistent_action} \\
 &\sum_{m=1}^M\sum_{a \in \calX_A^m} \bfD^m(a) \tau_a^{t,m} \leq \bfb & \forall t \in [T] \label{eq:wkpomdp:decPOMDP_NLP_linking_cons}
 \end{alignat}
\end{subequations}
which is obtained by using only modifications (A) and (C), is a relaxation.
We denote by $z_{\rmUB}$ the optimal value of the nonlinear program~\eqref{pb:wkpomdp:decPOMDP_NLP}.
\begin{theo}\label{theo:wkpomdp:upper_bound_NLP}
	The nonlinear program~\eqref{pb:wkpomdp:decPOMDP_NLP} is a relaxation of~\ref{pb:decPOMDP_wc} and MILP~\eqref{pb:wkpomdp:decPOMDP_MILP}. In particular, $v_{\rm{ml}}^* \leq z_{\rm{UB}}$ and $z_{\rmIP} \leq z_{\rm{UB}}.$
\end{theo}
Once again, variables $\bftau^m$ can be interpreted as the vector of moments of the probability distribution $\bbP_{\bfdelta^m}$ on component $m$ but there is no guarantee that it defines a joint probability distribution over the whole system. 
Problem~\eqref{pb:wkpomdp:decPOMDP_NLP} is a Quadratically Constrained Quadratic Program (QCQP) due to constraints~\eqref{eq:wkpomdp:NLP_per_POMDP_v1} and is in general non-convex. 
Recent QCQP solvers such as \texttt{Gurobi 9.0} \citep{gurobi} are able to solve small instances of Problem~\eqref{pb:wkpomdp:decPOMDP_NLP} to optimality, but they cannot address large instances due to the limits of a Spatial Branch-and-Bound \citep{Liberti2008}.

\subsubsection{A tractable upper bound through Lagrangian relaxation}
\label{sub:wkpomdp:lagrangian_relaxations}

The main computational difficulty in~\eqref{pb:wkpomdp:decPOMDP_NLP} comes from the nonlinear constraints $\tau_{soa}^{t,m} = \delta_{a|o}^{t,m}\sum_{\substack{o' \in \calX_O^m \\ a' \in \calX_A^m}} \tau_{soa'}^{t,m}$, which we cannot linearize using McCormick inequalities as in Section~\ref{sub:pomdp:MILP} because we do not assume that $\bfdelta^m$ is integer anymore.
However, if we perform a Lagrangian relaxation where we relax linking constraint~\eqref{eq:wkpomdp:decPOMDP_NLP_linking_cons}, then
we obtain an subproblem for each component that is a POMDP, which we can now reformulate as a MILP using the approach of Section~\ref{sub:pomdp:MILP}.
Let $z_{\rmLR}$ denote the value of this Lagrangian relaxation of NLP~\eqref{pb:wkpomdp:decPOMDP_NLP}.

\begin{restatable}{prop}{lagrangian}\label{prop:wkpomdp:lagrangian_relax}
	The value of the Lagrangian relaxations of MILP~\eqref{pb:wkpomdp:decPOMDP_MILP} and NLP~\eqref{pb:wkpomdp:decPOMDP_NLP} are equal and the following inequalities hold:
	\begin{align}\label{eq:wkpomdp:weak_duality}
		z_{\rmUB} \leq z_{\rmLR} \leq z_{\rmRc} \leq z_{\rmR}
	\end{align}
\end{restatable}

Proposition~\ref{prop:wkpomdp:lagrangian_relax} ensures that $z_{\rmLR}$ can be computed using the Lagrangian relaxation of MILP~\eqref{pb:wkpomdp:decPOMDP_MILP}.
Geoffrion's theorem \citep{Geoffrion1974} ensures that we can compute $z_{\rmLR}$ by solving the Dantzig-Wolfe reformulation of MILP~\eqref{pb:wkpomdp:decPOMDP_MILP} within a column generation approach.
The proof of Proposition~\ref{prop:wkpomdp:lagrangian_relax} and the full details of the column generation algorithm are available in Appendix~\ref{app:ColGen}.
Note that Lagrangian relaxation has already been used in the literature on weakly coupled stochastic dynamic programs to compute upper bounds \citep{Adelman2008,Ye2014,hawkins2003langrangian}.




\subsubsection{Benefits and drawbacks of the formulations}
\label{sub:wkpomdp:lb_ub:size_formulations}




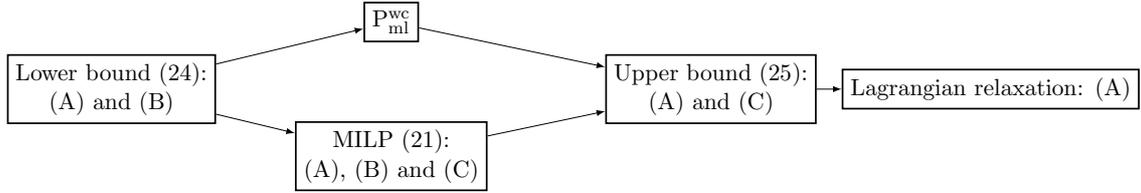
\begin{figure}
	\begin{center}
		\resizebox{15cm}{!}{
		\begin{tikzpicture}
			\def\l{1.3}
            \def\h{1.1}
            \node[draw,thick] (decPOMDP) at (-1.5*\l,0.0*\h) [align=center]{\ref{pb:decPOMDP_wc}};
            \node[draw,thick] (decPOMDPlb) at (-5.0*\l,-1.0*\h) [align=center]{Lower bound~\eqref{pb:wkpomdp:decPOMDP_MILP_LB}:\\
            (A) and (B)};
            \node[draw,thick] (decPOMDPmilp) at (-1.5*\l,-2.0*\h) [align=center]{MILP~\eqref{pb:wkpomdp:decPOMDP_MILP}:\\ 
            (A), (B) and (C)};
			\node[draw,thick] (decPOMDPnlp) at (2.5*\l,-1.0*\h) [align=center]{Upper bound~\eqref{pb:wkpomdp:decPOMDP_NLP}:\\ 
            (A) and (C)};
            \node[draw,thick] (lagrange) at (6.0*\l,-1.0*\h) [align=center]{Lagrangian relaxation: (A)};

            \draw[arc] (decPOMDPlb) -- (decPOMDP);
            \draw[arc] (decPOMDPlb) -- (decPOMDPmilp);
            \draw[arc] (decPOMDP) --  (decPOMDPnlp);
            \draw[arc] (decPOMDPmilp) --  (decPOMDPnlp);
            \draw[arc] (decPOMDPnlp) --  (lagrange);
	\end{tikzpicture}}
	\end{center}
	\caption{The relaxations of Section~\ref{sub:wkpomdp:lb_ub}. An arrow from Problem $X$ to a Problem $Y$ indicates that $Y$ is a relaxation of $X$ in the sense we defined at the beginning of this section. In each block, we indicate which approximations we use to obtain the formulations.}
	\label{fig:wkpomdp:interpret_bound}
\end{figure}

Figure~\ref{fig:wkpomdp:interpret_bound} summarizes the links between the feasible sets of the different formulations they have been established in the theorems.
Table~\ref{tab:wkpomdp:comparison_size_formulations} highlights the benefits and the drawbacks of the different formulations. 
It reports the behavior of the formulations regarding several criteria formulated as questions: are the numbers of variables (Pol. variables) and constraints (Pol. constraints) polynomial? Does the formulation have linking constraints between the components (Link. constraints)? Is the formulation linear (Linearity)? Are there integer variables in the formulation (Int. variables)? 
Does the formulation provide a feasible policy (Feas. pol.)? 
Is the optimal value an upper bound or a lower bound regarding $v_{\mathrm{ml}}^*$? Is the formulation tractable regarding the size of the instances (small with $\vert \calX_S \vert \leq 10^2$, medium with $10^2 \leq \vert \calX_S \vert \leq 10^4$ and large with $\vert \calX_S \vert \geq 10^4$)?
The tractability criteria should be understood as an advice on the formulation to choose and the scale order is only one indicator among others. 

\begin{table}[h!]
  \begin{center}
    \resizebox{16cm}{!}{
    \begin{tabular}{lccccccccccc}
      Formulations & \rot{Pol. variables} & \rot{Pol. constraints} & \rot{Link. constraints} & \rot{Linearity} & \rot{Int. variables} & \rot{Feas. policy} & \rot{Lower bound} & \rot{Upper bound} & Small & \multicolumn{1}{b{1.5cm}}{Tractability \newline Medium} & Large \\
      \hline
      NLP~\eqref{pb:wkpomdp:POMDP_NLP} &  &  & \checkmark & & & \checkmark &$-$ & $-$ & & & \\
      MILP~\eqref{pb:wkpomdp:decPOMDP_MILP} & \checkmark & \checkmark & \checkmark & \checkmark & \checkmark & & & & \checkmark & \checkmark&  \\
      Lower bound~\eqref{pb:wkpomdp:decPOMDP_MILP_LB} & \checkmark &  & \checkmark & \checkmark & \checkmark & \checkmark & \checkmark &  & \checkmark &  &  \\ 
      Upper bound~\eqref{pb:wkpomdp:decPOMDP_NLP} & \checkmark & \checkmark & \checkmark & & & & & \checkmark & \checkmark & & \\
      Lagrangian Relaxation
      & \checkmark & \checkmark &  & \checkmark & \checkmark & & & \checkmark & \checkmark& \checkmark& \checkmark  \\
      Linear Relaxation of~\eqref{pb:wkpomdp:decPOMDP_MILP} & \checkmark & \checkmark & \checkmark & \checkmark &  & &  & \checkmark &\checkmark & \checkmark & \checkmark\\
      \hline
    \end{tabular}}
  \end{center}
  \caption{Comparison of the properties of the formulations.}
  \label{tab:wkpomdp:comparison_size_formulations}
\end{table}


\subsection{Valid inequalities}
\label{sub:wkpomdp:valid_ineq}

Since Inequalities~\eqref{eq:pomdp:Valid_cuts_pomdp} are valid for MILP~\eqref{pb:pomdp:MILP_pomdp}, their local counterpart~\eqref{eq:wkpomdp:dec_Valid_cuts} are valid for all the MILPs introduced in this section.


\begin{prop}\label{prop:wkpomdp:decPOMDP_valid_cuts}
	Inequalities~\eqref{eq:wkpomdp:dec_Valid_cuts} are valid for MILP~\eqref{pb:wkpomdp:decPOMDP_MILP}, MILP~\eqref{pb:wkpomdp:decPOMDP_MILP_LB} and Problem~\eqref{pb:wkpomdp:decPOMDP_NLP}, and there exists a solution of the linear relaxation of~\eqref{pb:wkpomdp:decPOMDP_MILP} that does not satisfy constraints~\eqref{eq:wkpomdp:dec_Valid_cuts}.
\end{prop}

As in the generic POMDP case, inequalities~\eqref{eq:wkpomdp:dec_Valid_cuts} help the resolution of MILP~\eqref{pb:wkpomdp:decPOMDP_MILP} in practice. 
However, since the extended formulation obtained by adding inequalities~\eqref{eq:wkpomdp:dec_Valid_cuts} in MILP~\eqref{pb:wkpomdp:decPOMDP_MILP} has a large number of variables and constraints when the number of components is large ($M \geq 15$), the linear relaxation takes longer to solve. 
\subsection{Deducing a history-dependent policy from MILP~\eqref{pb:wkpomdp:decPOMDP_MILP}}
\label{sub:wkpomdp:first_policy}


In MILP~\eqref{pb:wkpomdp:decPOMDP_MILP}, we consider ``local'' policies $\bfdelta^m$ on each component $m$ in $[M].$
However, in general, given a vector of ``local'' policies $(\bfdelta^m)_{m\in[M]},$ there is no guarantee that there exists a policy $\bfdelta$ that coincides with $\bfdelta^m$ for every components $m$ in $[M].$
In this section we describe how we can use MILP~\eqref{pb:wkpomdp:decPOMDP_MILP} to deduce a history-dependent policy for weakly coupled POMDPs.

Consider a history $\bfh = (\bfo_1,\bfa_1,\ldots,\bfo_{t-1},\bfa_{t-1}, \bfo_t)$ available at time $t$.
Conditionally to $\bfh$, the vectors of state component $\left(S_{t'}^m\right)_{1\leq t'\leq t}$ for all $m$ in $[M]$ become independent, i.e.,
$$\bbP_{\bfdelta} \left(\bfS_t =\bfs \vert \bfH_t=\bfh \right) = \prod_{m=1}^M \overbrace{\bbP_{\bfdelta}\left(S_t^m=s^m \vert H_t^m=h^m \right)}^{p^m(s^m|h^m)}.$$
In the POMDP literature, the probability distribution $\bbP_{\bfdelta}\left(S_t^m\vert H_t^m \right)$ is called the \emph{belief state} of component $m$.
We can use the \emph{belief state update} (see Appendix.\ref{sub:app_wkpomdp_heuristic:belief_state}) on each of the components to compute the belief state $p^m(s^m|h^m)$.
We introduce the following algorithm:
\begin{algorithm}[H]
\caption{History-dependent policy $\rm{Act}_{\rmIP,T}^t(\bfh)$}
\label{alg:wkpomdp:heuristic_individual}
\begin{algorithmic}[1]
\STATE \textbf{Input} An history of observations and actions $\bfh \in (\calX_O \times \calX_A)^{t-1}\times \calX_O$.
\STATE \textbf{Output} An action $\bfa \in \calX_A.$
\STATE Compute the belief state $p^m(s|h^m)$ according to the belief state update (see Appendix.~\ref{sub:app_wkpomdp_heuristic:belief_state}) for every state $s$ in $\calX_S^m$ and every component $m$.
\STATE Remove constraints and variables indexed by $t'<t$ in MILP~\eqref{pb:wkpomdp:decPOMDP_MILP} and solve the resulting problem with horizon $T - t$, initial probability distributions $\left(p^m(s|h^m)\right)_{s \in \calX_S^m}$ for every component $m$ in $[M]$ and initial observation $\bfo_t$ (see Remark~\ref{rem:pomdp:with_observation}) to obtain an optimal solution $(\bftau^m,\bfdelta^m)_{m\in [M]}.$
\label{alg:wkpomdp:modify_constraints}
\STATE Return $\bfa = (a^{1},\ldots,a^{M})$ for which $\delta^{t,m}_{a^{m}|o^m} = 1$ for all $m$ in $[M].$ \label{alg:wkpomdp:take_action}
\end{algorithmic}
\end{algorithm}

Then we define the \emph{implicit} policy $\bfdelta^{\rmIP}$ as follows:
\begin{align*}
	\delta_{\bfa|\bfh}^{\rmIP,t}= \begin{cases}
							1, & \text{if}\ \bfa=\mathrm{Act}_{T}^{\rmIP,t}(\bfh) \\
      						0, & \text{otherwise}
						\end{cases}, & \quad \forall \bfh \in \calX_H^t, \enskip \bfa \in \calX_A, \enskip  t \in [T],
\end{align*}
where implicit means that each value of the policy is obtained by solving a mathematical programming formulation. 
It is not clear at first sight that policy $\bfdelta^{\rmIP}$ is a feasible policy in $\rm{P}_{\rm{his}}^{\rm{wc}}$ because it is not immediate to see that the action returned by Algorithm~\ref{alg:wkpomdp:heuristic_individual} belongs to $\calX_A.$
The theorem below ensures that the implicit policy $\bfdelta^{\rmIP}$ is a feasible policy of $\rm{P}_{\rm{his}}^{\rm{wc}}$, and that the belief updates can only improve the total expected reward. 
We denote by $\nu_{\rmIP}$ the total expected reward induced by policy $\bfdelta^{\rmIP}$.

\begin{theo}\label{theo:wkpomdp:heuristic}
	The implicit policy $\bfdelta^{\rmIP}$ is a feasible policy of $\rm{P}_{\rm{his}}^{\rm{wc}}$ and the inequality $z_{\rmIP} \leq \nu_{\rmIP} \leq v_{\rmhis}^*$ holds.
\end{theo}

It turns out that we can also use lower bound~\eqref{pb:wkpomdp:decPOMDP_MILP_LB} instead of MILP~\eqref{pb:wkpomdp:decPOMDP_MILP} at step~\ref{alg:wkpomdp:modify_constraints}. We denote by $\rm{Act}_{T}^{\rmLB,t}(\bfh)$ such an algorithm and $\bfdelta^{\rmLB}$ the induced policy.
Thanks to Theorem~\ref{theo:wkpomdp:lower_bound_MILP}, Theorem~\ref{theo:wkpomdp:heuristic} remains true when we use the induced implicit policy with $\rm{Act}_{T}^{\rmLB,t}(\bfh)$.
However, we cannot use our other formulations because the resulting actions would not necessary belong to $\calX_A$.
Appendix~\ref{app:LPpolicy} introduces a variant of Algorithm~\eqref{alg:wkpomdp:heuristic_individual} which can be used with these formulations, as well as numerical experiments showing that the resulting policies are not as efficient as $\bfdelta^{\rmIP}$.

\subsection{Rolling horizon matheuristic}
\label{sub:wkpomdp:matheuristic}


When the horizon $T$ is long it is computationally interesting to embed the implicit policy $\bfdelta^{\rmIP}$ in a rolling horizon heuristic, which consists in repeatedly solving an optimization problem with a smaller horizon $T_{\rmr} < T$ at each time step and to take action at the current time. 
This type of rolling horizon heuristic is commonly used for multistage stochastic optimization problems in operations research~\citep{Alistair2005}.
This approach is also called \emph{Model Predictive Control}~\citep{Bertsekas} in the optimal control literature.
See Appendix~\ref{sub:app_wkpomdp_heuristic:rolling_horizon} for more details.
\section{Numerical experiments}
\label{sec:num}

In this section we provide numerical experiments on the mathematical program formulations for POMDP and weakly coupled POMDP.
First, we illustrate the efficiency of our integer formulation~\eqref{pb:pomdp:MILP_pomdp} and the valid inequalities~\eqref{eq:pomdp:Valid_cuts_pomdp} for POMDP on random instances, and then we show the relevance of using our memoryless policies on different kinds of instances from the literature.
Second, we show the efficiency of using our rolling horizon matheuristic of Section~\ref{sub:wkpomdp:matheuristic} on a maintenance problem. Full details of the experiments and additional numerical results can be found in Appendix~\ref{app:nums}. In particular, it also includes numerical results on multi-armed bandit problems similar to the ones of~\citet{bertsimas2016decomposable}.
All the mathematical programs have been written in \texttt{Julia} \citep{bezanson2017julia} with the \texttt{JuMP} \citep{DunningHuchetteLubin2017} interface and solved using \texttt{Gurobi} 9.0. \citep{gurobi} with the default settings. 
Experiments have been run on a server with 192Gb of RAM and 32 cores at 3.30GHz.

\subsection{Generic POMDPs: Random instances}
\label{sub:num:random_instances}

In this section, we provide numerical experiments on generic POMDPs showing the efficiency of the valid inequalities~\eqref{eq:pomdp:Valid_cuts_pomdp} for our MILP~\eqref{pb:pomdp:MILP_pomdp}.
We solve instances of POMDP over different finite horizon $T \in \{10,20\}$
For each triplet $(\calX_S,\calX_O,\calX_A)$, we generate randomly $30$ instances.
We refer the reader to Appendix~\ref{app:nums} for more details about how we generate the instances.
The average results over the $30$ instances are reported in Table~\ref{tab:num:random_instances_results}.
The first four columns indicate the size of state space $\big\vert \calX_{S}\big\vert$, observation space $\big\vert \calX_{O}\big\vert$, action space $\big\vert \calX_{A}\big\vert$ and time horizon $T$. 
The fifth column indicates the mathematical program used to solve Problem~\eqref{pb:POMDP} with (strengthened) or without (basic)valid inequalities~\eqref{eq:pomdp:Valid_cuts_pomdp}. In the last three columns, we report the integrity gap $g_{\rm{i}}$, the final gap $g_{\rm{f}}$, the percentage of instances solved to optimality $\mathrm{Opt}$ and the computation time ($\mathrm{Time}$).

Table~\ref{tab:num:random_instances_results} shows that the valid inequalities introduced for MILP~\eqref{pb:pomdp:MILP_pomdp} are efficient.
Indeed, the integrity gap is always significantly lower with the strengthened formulation, which shows the tightening of linear relaxation and it greatly reduces the computation time.
In addition, Inequality~\ref{eq:pomdp:inequality_information} ensures that the integrality gaps of~\eqref{pb:pomdp:MILP_pomdp}  reported in Table~\ref{tab:num:random_instances_results} are also bounds of the relative gap between $v_{\mathrm{ml}}^*$ and $v_{\rm{his}}^*.$

\begin{table}
\begin{minipage}{0.5\textwidth}
  \begin{center}
    \resizebox{8.0cm}{!}{
    \begin{tabular}{|c|cc|c|SScc|}
        \hline
        \multirow{3}{*}{$(\vert\calX_S\vert,\vert\calX_O\vert,\vert\calX_A\vert)$} & \multirow{3}{*}{$T$} & \multirow{3}{*}{$\big\vert\Delta_{\mathrm{ml}}^{\mathrm{d}}\big\vert$} & \multirow{3}{*}{Formulation} & \multicolumn{4}{c}{\textbf{MILP~\eqref{pb:pomdp:MILP_pomdp}}} \vline \\ 
         & & & & {$g_{\mathrm{i}}$} & {$g_{\mathrm{f}}$} & {Opt} & {Time} \\ 
         & & & & {(\%)} & {(\%)} & {(\%)} & (s) \\
         \hline
        \text{$(3,3,3)$} & 10 & \text{${10^{14}}$} & Basic & 6.02 & Opt & 100 &  1.49 \\ \cline{5-8}
              &    &   & Strengthened & 1.70 & Opt & 100 & 0.62 \\ \cline{2-8}
                & 20 & \text{${10^{28}}$} & Basic & 6.04 & Opt & 100 & 664 \\ \cline{5-8}
                &    &                    & Strengthened & 1.52 & Opt & 100 & 466 \\ \cline{2-8} 
          \hline
        \text{$(3,4,4)$} & 10 & \text{${10^{24}}$} &Basic & 9.51 & 0.34 & 87 &  512 \\ \cline{5-8}
                &    &                   & Strengthened & 3.16 & 0.18 & 87 & 514.4 \\ \cline{2-8}
                & 20 & \text{${10^{48}}$} & Basic & 9.64 & 1.96 & 43 & 2221 \\ \cline{5-8}
                &    &                    & Strengthened & 2.86 & 1.13 & 61 & 1731\\ \cline{2-8} 
          \hline
        \text{$(3,5,5)$} & 10 & \text{${10^{34}}$} &Basic & 9.33 & 0.83 & 57 & 1591 \\ \cline{5-8}
                &    &                   & Strengthened & 2.35 & 0.38 & 70 & 1113 \\ \cline{2-8}
                & 20 & \text{${10^{69}}$} & Basic & 9.60 & 3.30 & 26 & 2702 \\ \cline{5-8}
                &    &                    & Strengthened & 2.14 & 1.14 & 52 & 1879 \\ \cline{2-8} 
          \hline
  \end{tabular}}
  \end{center}
\end{minipage}
\begin{minipage}{0.5\textwidth}
  \begin{center}
    \resizebox{8.0cm}{!}{
    \begin{tabular}{|c|cc|c|SScc|}
        \hline
        \multirow{3}{*}{$(\vert\calX_S\vert,\vert\calX_O\vert,\vert\calX_A\vert)$} & \multirow{3}{*}{$T$} & \multirow{3}{*}{$\big\vert\Delta_{\mathrm{ml}}^{\mathrm{d}}\big\vert$} & \multirow{3}{*}{Formulation} & \multicolumn{4}{c}{\textbf{MILP~\eqref{pb:pomdp:MILP_pomdp}}} \vline \\ 
         & & & & {$g_{\mathrm{i}}$} & {$g_{\mathrm{f}}$} & {Opt} & {Time} \\ 
         & & & & {(\%)} & {(\%)} & {(\%)} & (s) \\
         \hline
        \text{$(4,3,3)$} & 10 & \text{${10^{14}}$} &Basic & 7.39 & Opt & 100 & 26 \\ \cline{5-8}
                &    &                   & Strengthened& 2.28 & Opt & 100 & 9.16 \\ \cline{2-8}
                & 20 & \text{${10^{28}}$} & Basic & 6.01 & 1.01 & 60 & 1598 \\ \cline{5-8}
                &    &                    & Strengthened & 2.03 & 0.32 & 80 & 987 \\ \cline{2-8}
          \hline
        \text{$(4,4,4)$} & 10 & \text{${10^{24}}$} & Basic & 12.19 & 0.98 & 65 & 1477 \\ \cline{5-8}
                &    &                   &  Strengthened & 3.44 & 0.27 & 80 & 967 \\ \cline{2-8}
                & 20 & \text{${10^{48}}$} & Basic & 12.29 & 4.66 & 20 & 2900 \\ \cline{5-8}
                &    &                    & Strengthened & 3.05 & 1.48 & 30 & 2651 \\ \cline{2-8}
          \hline
        \text{$(4,5,5)$} & 10 & \text{${10^{34}}$} &Basic & 11.64 & 1.76 & 35 & 2427 \\ \cline{5-8}
                &    &                   & Strengthened & 3.09 & 0.62 & 65 & 1345 \\ \cline{2-8}
                & 20 & \text{${10^{69}}$} & Basic & 12.04 & 5.46 & 5 & 3413 \\ \cline{5-8}
                &    &                    & Strengthened & 3.20 & 1.67 & 32 & 2490\\ \cline{2-8} 
          \hline
  \end{tabular}}
  \end{center}
\end{minipage}
\caption{Efficiency of the valid inequalities~\eqref{eq:pomdp:Valid_cuts_pomdp} for MILP~\eqref{pb:pomdp:MILP_pomdp} on  random instances of generic POMDPs, with a time limit of $3600$s}
 \label{tab:num:random_instances_results}
 \end{table}

\subsection{Generic POMDP: Instances from the literature}
\label{sub:num:literature}

In this section, we evaluate the efficiency of MILP~\eqref{pb:pomdp:MILP_pomdp} on instances of POMDP drawn from the literature and we compare its performances with one of the state-of-the-art POMDP solver SARSOP of~\citet{Kurniawati08sarsop}. 
In particular, it shows how the memoryless policy provided by MILP~\eqref{pb:pomdp:MILP_pomdp} performs on instances from the literature.

In fact, SARSOP solver gives an approximate history-dependent policy for the discounted infinite horizon POMDP problem. 
To adapt this policy for the finite horizon POMDP problem, we proceed as~\citet{DujardinDC15}: We compute a policy using SARSOP solver with a discount factor $\gamma = 0.999$ and we compute the expected sum of rewards over the $T$ time steps by simulation of the history-dependent policy. We perform $10^4$ simulations to compute the expectation.
By definition, the objective value $z_{\mathrm{SARSOP}}$ obtained by using this policy is a lower bound of $v_{\mathrm{his}}^*.$
We run the SARSOP algorithm using the Julia library \texttt{POMDPs.jl} of \citet{EgorovSBWGK17}.

All the instances can be found at the link~\url{http://pomdp.org/examples/} and further descriptions of each instance are available in the indicated literature on the same website.
In particular, it contains two instances \texttt{bridge-repair} and \texttt{machine} that model maintenance problems.
For each instance, we report the objective functions $z^*$ and $z_{\rm{SARSOP}}$, and the relative gaps $g(z) = \frac{z_{\rmRc}^* - z}{z_{\rmRc}^*}$ for any $z$ belonging to $\{z^*, z_{\mathrm{SARSOP}}\}.$ 
The lower the value of $g(z)$, the closer the value of $z$ is to $v_{\rmhis}^*$.
All the results are reported in Table~\ref{tab:pomdp:benchmark_instances_results}.
The first column indicates the instance considered.
The three next columns indicate respectively the cardinality of $\calX_S,$ $\calX_O$ and $\calX_A$ of the instance.
The fourth column indicates the algorithm used. 
The last six columns indicate the total expected reward (Obj.) and the gap values for different finite horizon $T \in \{5,10,20\}.$

One may observe that in most cases the optimal value obtained with our MILP is close to the upper bound $z_{\rmRc}^*.$ Thanks to Theorem~\ref{theo:pomdp:MDP_approx_equivalence}, it says that memoryless policies perform well on finite horizon for these instances. In particular, the gap is noticeably small on the instance of maintenance problems. 
However, as mentioned in the introduction, one can observe that the memoryless policies fail on instances of navigation problems \cite{Littman94memoryless}. We observe this phenomenon on instances of navigation problems, where the goal is to find a target in a maze, and there are a large number of states relatively to a small number of observations. 
It is fairly natural: using a memoryless policy in a maze is misleading because if the decision maker meets a wall, he will act as it is the first time he meets a wall, and then will always take the same actions.
It seems that on these instances, the SARSOP policies work best on larger horizons, which is expected since the SARSOP policy is built for an infinite horizon problem.
The results in Table~\ref{tab:pomdp:benchmark_instances_results} support the remark of~\citet[Section 3]{Walraven2019} saying that the point-based algorithms for infinite discounted POMDP, such as SARSOP, produce policies that can be inefficient on finite horizon.

\begin{table}
    \begin{center}
    \resizebox{14cm}{!}{
    \begin{tabular}{|c|ccc|c|cccccc|}
        \hline
        \multirow{3}{*}{Instances} & \multicolumn{3}{c}{Size} \vline& \multirow{3}{*}{Algorithms} & \multicolumn{6}{c}{Horizon} \vline\\ 
        & {$\vert \calX_S \vert$} & {$\vert \calX_O \vert$} & {$\vert \calX_A \vert$}&  & \multicolumn{2}{c}{$T=5$} & \multicolumn{2}{c}{$T=10$} & \multicolumn{2}{c}{$T=20$} \vline\\ 
        & & & & & Obj. & Gap({$\%$}) & Obj. & Gap({$\%$}) & Obj. & Gap({$\%$}) \\
         \hline
        \texttt{1d.noisy} & 4 & 2 & 2 & MILP & \textbf{1.56} & \textbf{18.73} & \textbf{2.97} & \textbf{19.18} & \textbf{5.82} & \textbf{18.73}\\
                          & & & & SARSOP & 0.57 & 70.12 & 0.67 & 81.76 & 0.81 & 88.71 \\
        \texttt{4x5x2}$^{*}$ & 39 & 4 & 4 & MILP & \textbf{0.37} & \textbf{58.13} & \textbf{0.75} & \textbf{57.45} & \textbf{1.86} & \textbf{47.58}\\
                       & & & & SARSOP & 0.08 & 90.87 & 0.08 & 95.28 & 0.08 & 97.50 \\
        \texttt{aircraftID} & 12 & 5 & 6 & MILP & \textbf{5.69} & \textbf{0.00} & \textbf{10.10} & \textbf{0.00} & \textbf{19.76} & \textbf{0.00}\\
                     & & & & SARSOP & 3.39 & 40.46 & 7.63 & 24.46 & 17.32 & 12.41 \\
        \texttt{aloha.10} & 30 & 3 & 9 & MILP & 38.04 & 0.56 & 62.74 & 1.66 & 84.92 & 13.84\\
                          & & & & SARSOP & \textbf{38.15} & \textbf{0.25} & \textbf{63.74} & \textbf{0.20} &  \textbf{89.09} & \textbf{9.61} \\
        \texttt{cheng.D3-1} & 3 & 3 & 3 & MILP & \textbf{32.29} & \textbf{1.87} & \textbf{64.38} & \textbf{1.11} & \textbf{128.55} & \textbf{0.72}\\
                          & & & & SARSOP & 32.04 & 2.65 &  64.16 & 1.45 & 128.28 & 0.93 \\
        \texttt{cheng.D4-1} & 4 & 4 & 4 & MILP & \textbf{33.83} & \textbf{5.20} & \textbf{67.37} & \textbf{4.10} & \textbf{134.45} & \textbf{3.54}\\
                          & & & & SARSOP & 32.40 & 9.1 & 65.90 & 6.19 & 133.05 & 4.54 \\
        \texttt{cheng.D5-1} & 5 & 5 & 5 & MILP & \textbf{32.89} & \textbf{3.28} & \textbf{65.64} & \textbf{2.25} & \textbf{131.12} & \textbf{1.73}\\
                          & & & & SARSOP & 32.47 & 4.50 & 65.23 & 2.86 & 130.81 & 1.96 \\
        \texttt{learning.c3} & 24 & 3 & 12 & MILP & \textbf{1.63} & \textbf{45.3} & \textbf{2.20} & \textbf{26.76} & \textbf{2.33} & \textbf{22.48}\\
                          & & & & SARSOP & 0.33 & 88.89 & 0.33 & 89.00 & 0.34 & 88.67 \\
        \texttt{milos-aaai97}$^{*}$ & 20 & 8 & 6 & MILP & \textbf{26.83} & \textbf{10.28} & \textbf{53.41} & \textbf{36.06} & 92.09 & 55.06\\
                          & & & & SARSOP & 12.62 & 57.79 & 39.52 & 52.69 & \textbf{97.73} & \textbf{52.31} \\
        \texttt{network}$^{*}$ & 7 & 2 & 4 & MILP & 20.30 & 2.49 & 95.06 & 22.85 & 203.87 & 36.02\\
                          & & & & SARSOP & \textbf{20.88} & \textbf{0.00} & \textbf{95.78} & \textbf{22.26} & \textbf{245.88} & \textbf{22.98} \\
        \texttt{bridge-repair}$^{**}$ & 5 & 5 & 10 & MILP & \textbf{1992.77} & \textbf{0.15} & \textbf{7801.56} & \textbf{0.44} & \textbf{27937.93} & \textbf{0.13}\\
                          & & & & SARSOP & 1514.15 & 24.13 & 6832.99 & 12.80 & 26568.42 & 5.03 \\
        \texttt{query.s2} & 9 & 3 & 2 & MILP & \textbf{21.54} & \textbf{0.95} & \textbf{46.25} & \textbf{0.10} & \textbf{96.50} & \textbf{0.11}\\
                          & & & & SARSOP & 15.77 & 27.50 & 31.68 & 31.56 & 64.91& 30.66 \\
        \texttt{machine}$^{**}$ & 256 & 16 & 4 & MILP & \textbf{4.90} & \textbf{0.00} & \textbf{9.50} & \textbf{0.81} & \textbf{17.98} & \textbf{0.05}\\
                          & & & & SARSOP & 4.90 & 0.00 & 9.35 & 2.38 & 15.69 & 12.79 \\
        \hline
  \end{tabular}
  } 
  \\
  \footnotesize{$^{*}$ Instances of navigation problem, $^{**}$ Instances of maintenance problem}
  \end{center}
 \caption{Performances of our memoryless policy on benchmark instances. The results written in bold indicate the best value obtained for each instance.}
 \label{tab:pomdp:benchmark_instances_results}
 \end{table}

\subsection{Weakly coupled POMDP: Performances of the matheuristic on a maintenance problem}
\label{sub:nums:implicit_policy}

The aim of this section is to show how close $\nu_{\rmIP}$ is to the optimal value $v_{\rmhis}^*$, and that policy $\bfdelta^{\rmIP}$ can be computed in a reasonable amount of time on large-scale instances of a practical problem.
We evaluate the performances of the history-dependent policy $\bfdelta^{\rmIP}$ by running Algorithm~\ref{alg:wkpomdp:matheuristic} on a maintenance problem taken from the literature.
Like \citet[Section 5.2]{Walraven2018}, we consider a road authority that performs maintenance on $M$ bridges, each of them evolving independently over a finite horizon $H$.
Each bridge is modeled as a POMDP \citep{Ellis1995} with $5$ possible states and observations, and the authority must chooses at most $K$ bridges to maintain at each decision time. 
As mentionned in Example~\ref{sub:problem:example}, this problem can be modeled as a weakly coupled POMDP.
In Appendix~\ref{sub:app_nums:implicit_policy} we describe the maintenance problem and the settings.
We consider instances with different values of $M$ in $\{3,4,5,10,15,20\}$.
We choose a maintenance capacity $K = \max(\floor*{\gamma \times M},1)$, where $\gamma$ is a scalar belonging $\left\{ 0.2, 0.4, 0.6, 0.8\right\}$ (when $M=3$, then $K$ belongs to $\{1,2\}$) and we solve the problem over the finite horizon $T=24$.
We evaluate the performances of matheuristic~\ref{alg:wkpomdp:matheuristic} with formulation $\rmIP$ for rolling horizon $T_{\rmr}$ in $\{2,5\}$.
For each instance $(M,K,(\pfrak^m)_{m\in [M]})$, we perform $10^3$ runs of matheuristic~\ref{alg:wkpomdp:matheuristic}. We compute the average total cost $\vert \nu_{\rmIP} \vert$, the average number of failures $\rmF{\rmIP}$ over the $10^3$ simulations.
We compare $\nu^{\rmIP}$ with the upper bound $z_{\rmRc}$ and the Lagrangian bound $z_{\rmLR}$ by evaluating the average gap $\rmG_{\rmIP}^{\rmRc} = \frac{z_{\rmRc} - \nu_{\rmIP}}{\vert z_{\rmRc} \vert}$ and $\rmG_{\rmIP}^{\rmLR} = \frac{z_{\rmLR} - \nu_{\rmIP}}{\vert z_{\rmLR} \vert}$. 
Thanks to Theorem~\ref{theo:wkpomdp:ineq_decPOMDP}, the value of $G_{\rmIP}^{\rmRc}$ indicates how far is $\nu_{\rmIP}$ from $v_{\rmhis}^*$ because $\nu_{\rmIP}\leq v_{\rmhis}^* \leq  z_{\rmRc}$.
The lower the value of $G_{\rmIP}^{\rmRc}$, the better is the performance of policy $\bfdelta^{\rmIP}$. 
In addition, for each simulation, we compute the average computation time in seconds of the underlying formulation over all steps of the simulation. We then consider the average value over all the $N$ simulations.
For the quantities $\vert \nu_{\rmIP} \vert$ and $\rmF_{\rmIP}$ we also report the standard deviations over all simulations.

Tables~\ref{tab:wkpomdp:maintenance} displays several results (see Appendix~\ref{sub:app_nums:implicit_policy} for more results). 
For all the mathematical programs, we set the computation time limit to $3600$ seconds and a final gap tolerance (\texttt{MIPGap} parameter in \texttt{Gurobi}) of $1 \%$, which is enough for the use of our matheuristic.
If the resolution has not terminated before this time limit, then we keep the best feasible solution at the end of the resolution.

\begin{table}
  \centering
    \resizebox{!}{5.0cm}{
    \begin{tabular}{|c|cc|ccccccc|}
        \hline
          \multirow{2}{*}{$M$} & \multirow{2}{*}{$\gamma$} & \multirow{2}{*}{$T_\rmr$} & $\vert \nu_{\rmIP} \vert$ & Std. err. & $\rmF_{\rmIP}$ & Std. err. & $\rmG_{\rmIP}^{\rmLR}$ & $\rmG_{\rmIP}^{\rmRc}$ & Time \\
              &  & & {($\times10^3$)} & {($\times10^3$)} &  &  & $(\%)$ & $(\%)$ & ($\si{s}$) \\
              \hline
       $10$   & 0.2  & 2 & 22.63 & 5.45 & 18.0 & 5.5 & -34.20 & 17.71 & 0.012 \\
                &    & 5 & 22.06 & 5.24 & 17.5 & 5.2 & \textbf{-35.85} & \textbf{14.74} & \textbf{0.384}\\\cline{2-10}
                & 0.4  & 2 & 19.19 & 3.38 & 11.5 & 3.3 & 1.32 & 3.07 & 0.012 \\
                &      & 5 & 18.91 & 3.26 & 10.6 & 3.2 & \textbf{-0.16} & \textbf{1.57} & \textbf{0.196} \\ \cline{2-10}
                & 0.6  & 2 & 19.10 & 3.28 & 10.8 & 3.1 & 0.92 & 2.60 & 0.011 \\
                &      & 5 & 18.81 & 3.06 & 9.7 & 2.9 & \textbf{-0.62} & \textbf{1.03} & \textbf{0.138} \\\cline{2-10}
                & 0.8 & 2 & 19.09 & 3.27 & 10.8 &3.1 & 0.86 & 2.53 & 0.011 \\
                &     & 5 & 18.82 & 3.08 & 9.6 & 2.9 & \textbf{-0.56} & \textbf{1.09} & \textbf{0.137} \\
                \hline
            15  & 0.2 & 2 & 31.54 & 6.03 & 25.0 & 6.0 & -22.56 & 12.73 & 0.017 \\
                &     & 5 & 30.89 & 5.88 & 24.0 & 5.9 & \textbf{-24.14} & \textbf{10.43} & \textbf{0.591} \\\cline{2-10}
                & 0.4  & 2 & 28.18 & 4.22 & 19.1 & 4.1 & 0.45 & 1.93 & 0.016 \\
                &      & 5 & 27.67 & 3.97 & 16.8 & 3.9 & \textbf{-1.39} & \textbf{0.06} & \textbf{0.232} \\ \cline{2-10}
                & 0.6  & 2 & 28.12 & 4.16 & 18.8 & 4.0 & 0.10 & 1.70 & 0.016 \\
                &      & 5 & 27.67 & 3.84 & 16.3 & 3.7 & \textbf{-1.51} & \textbf{0.05} & \textbf{0.225} \\\cline{2-10}
                & 0.8  & 2 & 28.12 & 4.16 & 18.8 & 4.0 & 0.11 & 1.71 & 0.015 \\
                &      & 5 & 27.65 & 3.86 & 16.2 & 3.7 & \textbf{-1.57} & \textbf{-0.01} &\textbf{0.226} \\ 
                \hline
              20  & 0.2& 2 & 45.06 & 7.01 & 35.9 & 7.1 & -20.37 & 8.67 & 0.022 \\
                &      & 5 & 44.28 & 6.90 & 35.1 & 6.9 & \textbf{-21.74} & \textbf{6.80} & \textbf{0.660} \\ \cline{2-10}
                & 0.4  & 2 & 41.18 & 4.72 & 23.0 & 4.7 & 0.35 & 1.50 & 0.020 \\
                &      & 5 & 40.83 & 4.66 & 22.7 & 4.7 & \textbf{-0.49} & \textbf{0.66} & \textbf{0.469} \\ \cline{2-10}
                & 0.6 & 2 & 40.96 & 4.25 & 18.4 & 4.1 & 0.32 & 1.22 & 0.020 \\
                &     & 5 & 40.72 & 4.19 & 17.9 & 4.0 & \textbf{-0.28} & \textbf{0.62} & \textbf{0.316} \\ \cline{2-10}
                & 0.8  & 2 & 40.96 & 4.24 & 18.3 & 4.0 & 0.29 & 1.21 & 0.019 \\
                &      & 5 & 40.76 & 4.09 & 17.8 & 3.9 & \textbf{-0.19} & \textbf{0.73} & \textbf{0.313} \\
            \hline
        \end{tabular}}
    \caption{Performances of the matheuritic on different rolling horizon $T_{\rmr} \in \{2,5\}$: Numerical values of $ \vert \nu_{\rmIP} \vert$, $\rmF_{\rmIP}$ (and the corresponding standard errors), $\rmG_{\rmIP}^{\rmLR}$ and $\rmG_{\rmIP}^{\rmRc}$ obtained on an instance $(M,\gamma)$ with $M \in \{10,15,20\}$ and $\gamma \in \{0.2, 0.4, 0.6, 0.8\}$. The values written in bold indicate the best performances of policy $\bfdelta^{\rmIP}$ regarding optimality and scalability (computation time).}
    \label{tab:wkpomdp:maintenance}
\end{table}

One may observe that for all instances, the matheuristic involving our MILP~\eqref{pb:wkpomdp:decPOMDP_MILP} delivers promising results even in the most challenging instance ($M=20$).
In particular, the values of $\rmG_{\rmIP}^{\rmRc}$ show that the policy $\bfdelta^{\rmIP}$ gives an optimality gap (in the set of history-dependent policies) of at most $10\%$ on the large-scale instance, which is satisfying regarding the complexity of the optimization problem ($\vert \calX_S \vert = \vert \calX_O \vert \approx 10^{14}$).
In Table~\ref{tab:wkpomdp:maintenance}, the negative values of $\rmG_{\rmIP}^{\rmRc}$ result from error approximations due to the Monte-Carlo simulations.
It can also be noted that the gap $\rmG_{\rmIP}^{\rmLR}$ is takes negative values for some instances, which shows that $\nu_{\rmIP}$ can take larger values than the Lagrangian relaxation for some instances. It highlights the benefit of using the belief state upadtes in the definition of $\bfdelta^{\rmIP}$.
In addition, even for the largest instances ($M = 15$ or $M=20$) and for $T=5$, the average time per action of $\rmAct_{T_\rmr}^{\rmIP,t}(\bfh_t)$ is on the order of $1.0$ second; this amount of time is still feasible even if the $24$ decision times are close together.  
In Appendix~\ref{sub:app_nums:implicit_policy}, the numerical results show that our policy involving MILP~\eqref{pb:wkpomdp:decPOMDP_MILP} gives better performances than other formulations $\{ \rmLB, \rmRc, \rmR \}$.
It would seem that using the Lagrangian relaxation within our matheuristic yields competitive results in terms of performances. However, it takes more time to compute the actions and the sampling leads to higher standard errors.

\section{Conclusion}
\label{sec:conclusion}

In this paper, we have considered several mathematical programming formulations for POMDPs with memoryless policies.
Valid cuts based on properties of independences strengthen these formulations.
And we have leveraged these formulations to build practically efficient history-dependent policies.
Furthermore, in order to break the curse of dimensionality that impedes the design of efficient policies for POMDPs modeling systems with several components we have introduced the notion of weakly coupled POMDPs, which generalizes the weakly coupled stochastic dynamic programs of~\citet{Adelman2008}.
And we have provided mathematical programming formulations and policies tailored for these weakly coupled POMDPs. Numerical experiments show their efficiency, notably on some maintenance problems.

Our history-dependent policies are designed to be very efficient on POMDPs where memoryless policies are rather efficient, but not on POMDPs where only tiny bits of informations are observed, such as maze problems, on which ones memoryless policies lead to poor decisions. 
Future directions include the design of efficient policies for such weakly coupled POMDPs with very few observations.

\section*{Acknowledgments}

We are grateful to Prof. Fr\'ed\'eric Meunier for his helpful remarks.
The authors gratefully acknowledge the financial support of the Operations Research and Machine Learning chaire between Ecole des Ponts Paristech and Air France.
%
%
%






\bibliographystyle{plainnat}
\bibliography{paperRO}

\newpage
\appendix

\section{Examples of weakly coupled POMDP applications}
\label{app:examples}

In this section we describe three examples of multi-stage stochastic optimization problems that can be formalized as a weakly coupled POMDP. 

\begin{ex}\label{ex:multi_armed_bandit}
    \emph{Multi-armed and Restless Bandit} problems are classical resource allocation problems where there are several arms, each of them evolving independently as a MDP, and at each time the decision maker has to activate a subset of arms so as to maximize its expected discounted reward over infinite horizon. We can consider the \emph{regular} multi-armed bandit problem, where only the activated arm states transit randomly and give an immediate reward, or the \emph{restless} multi-armed bandit problem, where all the arm states transit randomly and give an immediate reward.
    When the decision maker has only access to a partial observation on each arm instead of the arm state, the problem becomes a partially observable multi-armed bandit problem \citep{Krishnamurthy09}. 
    In this case, each arm evolves individually and independently as a POMDP. Such a problem enables to model practical applications such as clinical trials. In this setting, each component represents a medical treatment and activating a component corresponds to testing the treatment. The state of a medical treatment corresponds to its efficiency and the observation corresponds to a noisy measure of the efficiency of a medical treatment.

    We can formalize the partially observable multi-armed bandit problem within our weakly coupled POMDP framework. Let $M$ be the number of arms. At each time $t$, the decision maker has to activate $K < M$ arms. Since each component evolves as POMDP, we use the same notation to represent the state and the observation of Section~\ref{sub:problem:wkpomdp}.
    We define the individual action space $\calX_A^m =\{0,1\}$ of arm $m$ and the full action space is 
    $$\displaystyle \calX_A = \big\{\bfa \in \calX_A^1\times \cdots \times \calX_A^M \colon \sum_{m = 1}^M a^m = K\big\},$$
    which has the form~\eqref{eq:problem:def_form_action_space} by setting $\bfD^m(a)=\begin{pmatrix}a\\-a\end{pmatrix}$ for all $m$ in $[M]$ and $\bfb = \begin{pmatrix}K\\-K\end{pmatrix}$.
    In the case of the regular bandit problems, the immediate reward of component $m$ satisfies $r^m(s,0,s') = 0$, and the transition probabilities satisfy 
    $p^m(s'|s,0)$ equals $1$ if $s=s'$ and $0$ otherwise, for every $s,s' \in \calX_S^m$. 
    The goal of the decision maker is to find a policy $\bfdelta$ in $\Deltaml$ (or $\Deltahis$) maximizing the total expected discounted reward over infinite horizon $\bbE_{\bfdelta} \left[ \sum_{t=1}^{T}r\left(\bfS_t,\bfA_t,\bfS_{t+1}\right)\right]$,
    where $T$ is a finite horizon.
\end{ex}
\begin{ex}\label{ex:inventory_control}
    Consider a supplier that delivers a product to $M$ stores.
    At each time $t$, we denote by $S_t^m$ the inventory level of store $m$.
    Unfortunately, due to ''inventory records inaccuracy'' \citep{Mersereau2013} from various uncertainties, the supplier does not observe directly this inventory level. 
    He has instead only access to a noisy observation $O_t^m$ of the inventory level of store $m$.
    We assume that the inventory level of store $m$ has a known limited capacity $C^m$. Hence, we set $\calX_S^m:= \{0,\ldots,C^m\}$. Then $O_t^m=o$ is a noisy observation of the current inventory level, whose value belongs to $\calX_O^m:=\calX_S^m$ and is randomly emitted given a current state $S_t^m=s$ according to a known probability $p^m(o|s) = \bbP\left(O_t^m=o|S_t^m=s \right)$.
    At each time, the supplier has to decide the quantity to produce and to deliver automatically to each store.
    We denote by $A_t^m$ the quantity of product delivered to store $m$, which belongs to the individual action space $\calX_A^m := \calX_S^m$.
    The production has to satisfy resource constraints (raw materials, staff, etc.). Hence, the set of feasible actions has the form
    $$\calX_A := \left\{\bfa \in \calX_A^1\times \cdots \times \calX_A^M \colon \sum_{m=1}^M h^m a^m \leq H \right\},$$ 
    where $h^m$ is the given number of resources used per unit produced and delivered for store $m$ and $H$ is the given available amount of resource. 
    This action space has the form~\eqref{eq:problem:def_form_action_space} by setting $\bfD^m(a^m) = h^m a^m$ for all $m$ in $[M]$ and $\bfb=H$.
    The quantity of products in store $m$ cannot exceed capacity $C^m$. Hence, the quantity $\max(S_t^m+A_t^m - C^m, 0)$ is wasted and it induces a waste cost.

    We denote by $D_t^m$ the random variable representing the demand at store $m$ between time $t$ and $t+1$. The vector of demand is exogenous and independent identically distributed in each store with a known probability distribution $\bbP_D^m$ for store $m$.
    The inventory level of store $m$ follows the  dynamic 
    $$S_{t+1}^m = \max \left(\min\left(S_t^m + A_t^m,C^m\right) - D_t^m, 0 \right),$$ 
    which gives the transition probability distribution $\bbP(S_{t+1}|S_{t},A_t).$
    Now we can define the immediate reward function
    $$r^m(s,a,s') = \mathrm{price}^m(s+a-s') - \mathrm{waste}^m \max\left(s+a - C^m, 0 \right) -\mathrm{shortage}^m\bbE_{\bbP_D^m} \left[\max\left(D^m - (s+a), 0\right) \right],$$
    where $\mathrm{price}^m$ is the selling price per unit, $\rm{waste}^m$ is the wastage cost per unit and $\mathrm{shortage}^m$ is the shortage cost per unit. 
    It leads us to model this problem as a weakly coupled POMDP. The goal of the supplier is to find a policy $\bfdelta$ in $\Deltaml$ (or $\Deltahis$) maximizing the total expected reward over a finite horizon $T$.
    This example has been introduced by \citet{Kleywegt2002} for fully observable inventory levels and \citet{Mersereau2013} justifies the relevance of the POMDP framework for the stochastic inventory control problem.
\end{ex}

\begin{ex}\label{ex:medical_center}
    Consider a nurse assignment problem for home health care.
    A medical center follows $M$ patients at home on a daily basis over a given period of time $T$.
    On day $t$, we denote by $S_t^m$ the health state of patient $m$, whose value belongs to a finite state space $\calX_S^m$. 
    The medical center does not directly observe the health state of each patient. 
    However, at each time $t$, the medical center has access to a partial observation $O_t^m$ corresponding to a signal sent by a machine which diagnoses patient $m$. 
    We assume that this signal is discrete and noisy. Hence, $\calX_O^m$ is a finite space and an observation $o$ is randomly emitted given a state $s \in \calX_S^m$ according to the probability $p^m(o|s)$. 
    At each time, the medical center has to assign nurses to patients. There are $K_1$ available nurses with skill $1$ and $K_2$ available nurses with skill $2$. 
    On day $t$, we denote by $A_t^m$ the action taken by the medical center on patient $m$, whose value belongs to $\calX_A^m = \{0,1,2,3\}$ and the following meaning.
    \begin{align*}
        A_t^m = \begin{cases}
            & 0 \quad \text{if no nurse is sent to patient $m$} \\            
            & 1 \quad \text{if a nurse with skill $1$ is sent to patient $m$} \\
            & 2 \quad \text{if a nurse with skill $2$ is sent to patient $m$} \\
            & 3 \quad \text{if two nurses, one with each skill, are sent to patient $m$} \\
        \end{cases}
    \end{align*}
    Depending on the skill of the nurses sent to patient, the health state of each patient evolves randomly according to a transition probability $p^m(s'|s,a)$, for any $s,s' \in \calX_S^m$, $a \in \calX_A^m$ and $m \in [M]$.
    Hence, the set of feasible actions is
    $$\calX_A = \bigg\{\bfa \in \calX_A^1 \times \cdots \times \calX_A^M \colon \sum_{m=1}^M \mathds{1}_{1}(a^m) + \mathds{1}_{3}(a^m) \leq K_1 \ \text{and} \ \sum_{m =1}^M \mathds{1}_{2}(a^m) + \mathds{1}_{3}(a^m) \leq K_2 \bigg\},$$
    which has the form~\eqref{eq:problem:def_form_action_space} by setting $\bfD^m(a)= \begin{pmatrix}\mathds{1}_{1}(a) + \mathds{1}_{3}(a)\\ \mathds{1}_{2}(a) + \mathds{1}_{3}(a)\end{pmatrix}$ for all $m$ in $[M]$, and $\bfb = \begin{pmatrix} K_1 \\ K_2 \end{pmatrix}$.
    Now we can define the immediate reward function
    $$r^m(s,a,s') = -\rm{cost}_1^m(\mathds{1}_{1}(a)+\mathds{1}_{3}(a)) - \rm{cost}_2^m(\mathds{1}_{2}(a)+\mathds{1}_{3}(a)) - \rm{emergency}^m \mathds{1}_{s_{\rm{critic}}^m}(s'),$$
    where $\rm{cost}_i^m$ is the cost induced by sending a nurse with skill $i \in \{1,2\}$ to patient $m$, $s_{\rm{critic}}^m$ is the critical health state of patient $m$ and $\rm{emergency}^m$ is the cost induced by an emergency because patient $m$ reaches its critical health state.
    It leads us to model this problem as a weakly coupled POMDP. The goal of the medical center is to find a policy $\bfdelta$ in $\Deltaml$ (or $\Deltahis$) maximizing its total expected reward over a finite horizon $T$. 
\end{ex}

\section{Proofs of Section~\ref{sec:pomdp}}
\label{app:pomdp}

\proof[Proof of Theorem~\ref{theo:pomdp:NLP_optimal_solution}]
	Let $(\bfmu, \bfdelta)$ be a feasible solution of Problem~\eqref{pb:pomdp:NLP_pomdp}. 
	We prove by induction on $t$ that $\mu_s^1 = \bbP_{\bfdelta}\big(S_1 = s\big)$, $\mu_{soa}^t=\bbP_{\bfdelta}\big(S_t = s, O_t=o, A_t=a \big)$ and $\mu_{sas'}^t =\bbP_{\bfdelta}\big(S_t = s, A_t=a, S_{t+1} = s' \big)$. 
	At time $t=1$, the statement is immediate. Suppose that it holds up to $t-1$. Constraints~\eqref{eq:pomdp:NLP_indep_action}, \eqref{eq:pomdp:NLP_consistency_s} and induction hypothesis ensure that 
	\begin{align*}
	\mu_{soa}^{t} = \delta_{a|o}^t p(o|s) \sum_{o',a'} \mu_{so'a'}^t = \delta_{a|o}^t p(o|s) \sum_{s',a'} \mu_{s'a's}^{t-1} &= \delta_{a|o}^t p(o|s)  \sum_{s',a'} \bbP_{\bfdelta}\left(S_{t-1}=s', A_{t-1} = a', S_t=s \right) \\
	& = \delta_{a|o}^t p(o|s) \bbP_{\bfdelta}\left(S_t=s \right) \\
	&= \bbP_{\bfdelta}\left(S_t=s, O_t=o, A_t=a \right),
	\end{align*}
	where the last equality comes from the conditional independences and the law of total probability. 
	Constraints~\eqref{eq:pomdp:NLP_consistency_sa},\eqref{eq:pomdp:NLP_indep_state} and the induction hypothesis ensure that:
	\begin{align*}
		\mu_{sas'}^t = p(s'|s,a)\sum_{\ovs} \mu_{sa\ovs}^t = p(s'|s,a) \sum_{o} \mu_{soa}^t &= p(s'|s,a) \sum_{o} \bbP_{\bfdelta} \left( S_t=s, O_t=o, A_t=a \right) \\
																							&= \bbP_{\bfdelta}(S_t=s, A_t=a,S_{t+1}=s')
	\end{align*}
	where the last equality comes from the conditional independences and the law of total probability.
	Consequently,
	$$\sum_{t=1}^T \sum_{\substack{s,s' \in \calX_S \\ a \in \calX_A}} r(s,a,s') \bbP_{\bfdelta}\left(S_t=s,A_t=a,S_{t+1}=s' \right) = \bbE_{\bfdelta} \bigg[ \sum_{t=1}^{T}r(S_t,A_t,S_{t+1})\bigg],$$
	which implies that $\bfdelta$ is optimal if and only if $(\bfmu,\bfdelta)$ is optimal for Problem~\eqref{pb:pomdp:NLP_pomdp} and $v_{\rm{ml}}^* = z^*$. It achieves the proof.
\qed

\proof[Proof of Proposition~\ref{prop:pomdp:valid_cuts_pomdp}]
	Let $(\bfmu,\bfdelta)$ be a feasible solution of Problem~\eqref{pb:pomdp:MILP_pomdp}. We define
	$$\mu_{s'a'soa}^t = \delta^t_{a|o}p(o|s)\mu_{s'a's}^{t-1}$$ for all $(s',a',s,o,a) \in \calX_S \times \calX_A \times \calX_S \times \calX_O \times \calX_A$, $t \in [T]$. These new variables satisfy constraints in \eqref{eq:pomdp:Valid_cuts_pomdp} :
	\begin{align*}
		\sum_{a \in \calX_A} \mu_{s'a'soa}^t 
		&= \left(\sum_{a \in \calX_A} \delta^t_{a|o}\right)p(o|s)\mu_{s'a's}^{t-1} = p(o|s)\mu_{s'a's}^{t-1}\\
		\sum_{a' \in \calX_A, s' \in \calX_S} \mu_{s'a'soa}^t &= \left(\sum_{a' \in \calX_A, s' \in \calX_S} \mu_{s'a's}^{t-1}\right) \delta^t_{a|o} p(o|s) = \delta^t_{a|o} p(o|s) \sum_{o' \in \calX_O,a' \in \calX_O} \mu_{so'a'}^{t} = \mu_{soa}^t
	\end{align*}

	\noindent The remaining constraint \eqref{eq:pomdp:Valid_cuts_pomdp_main} is obtained using the following observation :
	\begin{align*}
		\frac{\mu_{s'a'soa}^t}{\sum_{s'' \in \calX_S} \mu_{s'a's''oa}^t} = \frac{p(o|s)\mu_{s'a's}^{t-1}}{\sum_{s'' \in \calX_S} p(o|s'')\mu_{s'a's''}^{t-1}} =  \frac{p(o|s)p(s|s',a')\sum_{\ovs} \mu_{s'a'\ovs}^{t-1}}{\sum_{s'' \in \calX_S} p(o|s'')p(s''|s',a')\sum_{\ovs} \mu_{s'a'\ovs}^{t-1}}  =  \frac{\displaystyle p(o|s)p(s|s',a')}{\displaystyle \sum_{s'' \in \calX_S} p(o|s'')p(s''|s',a')}
	\end{align*}
	By setting $p(s|s',a',o) = \frac{\displaystyle p(o|s)p(s|s',a')}{\displaystyle \sum_{\overline{s} \in \calX_S} p(o|\overline{s})p(\overline{s}|s',a')}$, equality \eqref{eq:pomdp:Valid_cuts_pomdp_main} holds. If $\sum_{s'' \in \calX_S} \mu_{s'a's''oa}^t = 0$, then $\mu_{s'a'soa}^t=0$ and constraint~\eqref{eq:pomdp:Valid_cuts_pomdp_main} is satisfied.

	Now we prove that there exists a solution $\bfmu$ of the linear relaxation of MILP~\eqref{pb:pomdp:MILP_pomdp} that does not satisfy inequalities~\eqref{eq:pomdp:Valid_cuts_pomdp}. We define such a solution $(\bfmu,\bfdelta)$.
	We set $\mu_s^1 = p(s)$ for all $s$ in $\calX_S$, and for all $t$ in $[T]$:
	\begin{align}
		&\mu_{soa}^1 = \begin{cases}
								&p(o|s)\mu_s^1, \ \text{if} \ a = \phi(s) \label{eq:pomdp:proof_mu_soa}\\
								&0,\  \text{otherwise} 
							\end{cases}, & \text{if \ $t=1$}, \\
		&\mu_{sas'}^{t} = p(s'|s,a) \sum_{o \in\calX_O} \mu_{soa}^t & \\
		&\mu_{soa}^t = \begin{cases}
								&p(o|s)\sum_{s'\in \calX_S,a' \in \calX_A} \mu_{s'a's}^{t-1}, \ \text{if} \ a = \phi(s) \label{eq:pomdp:proof_mu_soa}\\
								&0,\  \text{otherwise} 
							\end{cases}, & \text{if \ $t\geq 2$}, \\
		&\delta_{a|o}^t = \begin{cases}
								&\frac{\sum_{s \in \calX_S}\mu_{soa}^t}{\sum_{s \in \calX_S, a \in \calX_A} \mu_{soa}^{t}} \ \text{if} \ \sum_{s \in \calX_S, a \in \calX_A} \mu_{soa}^{t} \neq 0 \label{eq:pomdp:proof_delta}\\
								&\mathds{1}_{\tilde{a}}(a),\  \text{otherwise} 
							\end{cases} &
	\end{align}
	where $\phi : \calX_S \rightarrow \calX_A$ is an arbitrary mapping and $\tilde{a}$ is an arbitrary element in $\calX_A$. We prove that $\bfmu$ is a feasible solution of the linear relaxation of MILP~\eqref{pb:pomdp:MILP_pomdp}. 

	First, it is easy to see constraints~\eqref{eq:pomdp:NLP_initial2}-\eqref{eq:pomdp:NLP_indep_state} are satisfied. It remains to prove that constraints ~\eqref{eq:pomdp:MILP_McCormick_1}, ~\eqref{eq:pomdp:MILP_McCormick_2}, ~\eqref{eq:pomdp:MILP_McCormick_3} are satisfied.
	First, ~\eqref{eq:pomdp:MILP_McCormick_1} holds because
	\begin{align*}
		\mu_{soa}^t \leq \max\left(0,p(o|s) \sum_{s' \in \calX_S, a' \in \calX_A} \mu_{s'a's}^{t-1}\right) \leq p(o|s) \sum_{s' \in \calX_S, a' \in \calX_A} \mu_{s'a's}^{t-1},
	\end{align*}
	Second, ~\eqref{eq:pomdp:MILP_McCormick_2} holds because
	\begin{align*}
		\mu_{soa}^t \leq \sum_{s' \in \calX_S}\mu_{s'oa}^t = \delta_{a|o}^t \sum_{s' \in \calX_S}p(o|s')\sum_{s'' \in \calX_S, a'' \in \calX_A} \mu_{s''a''s'}^{t-1} \leq \delta_{a|o}^t, 
	\end{align*}
	where we used definition~\eqref{eq:pomdp:proof_delta} from the first to second line.
	Third, \eqref{eq:pomdp:MILP_McCormick_3} holds because
	\begin{align*}
		\mu_{soa}^t - p(o|s)\sum_{\substack{s' \in \calX_S\\ a' \in \calX_A}} \mu_{s'a's}^{t-1} & \geq \sum_{s'' \in \calX_S} \overbrace{\mu_{s''oa}^t - p(o|s'')\sum_{s' \in \calX_S, a' \in \calX_A} \mu_{s'a's''}^{t-1}}^{\leq 0} \\
		&= \sum_{s'',s' \in \calX_S, a' \in \calX_A} p(o|s'')\mu_{s'a's''}^{t-1}(\delta_{a|o}^t - 1) \\
		&\geq \delta_{a|o}^t - 1,
	\end{align*}
	which yields ~\eqref{eq:pomdp:MILP_McCormick_3}.
	Therefore, $(\bfmu,\bfdelta)$ is a solution of the linear relaxation of MILP~\eqref{pb:pomdp:MILP_pomdp}.

	Now, we prove that such a solution does not satisfy inequalities~\eqref{eq:pomdp:Valid_cuts_pomdp}. We define the new variables:
	\begin{align*}
	  	\mu_{s'a'soa}^t = \begin{cases}
	  						& \mu_{s'a's}^{t-1}\frac{\mu_{soa}^t}{\sum_{o'\in \calX_O,a' \in \calX_A}\mu_{so'a'}^t} \ \text{if} \ \sum_{o'\in \calX_O,a' \in \calX_A}\mu_{so'a'}^t \neq 0 \\
	  						& 0 \ \text{otherwise}
	  						\end{cases}
	\end{align*}
	Hence, $\bfmu$ satisfies constraints~\eqref{eq:pomdp:Valid_cuts_pomdp_consistency1} and \eqref{eq:pomdp:Valid_cuts_pomdp_consistency2}. However, constraint~\eqref{eq:pomdp:Valid_cuts_pomdp_main} is not satisfied in general. 
	Indeed, since the mapping $\phi$ is arbitrary, we can set $\phi$ such that $p(s|s',a',o) >0$ and $\mu_{s'a'soa}^t = 0$.
	Therefore, there exists a solution $\bfmu$ of the linear relaxation of MILP~\eqref{pb:pomdp:MILP_pomdp} that does not satisfy inequalities~\eqref{eq:pomdp:Valid_cuts_pomdp}. It achieves the proof.
\qed

\proof[Proof of Theorem~\ref{theo:pomdp:MDP_approx_equivalence}]
	We first prove the equivalence between the linear relaxation of our MILP~\eqref{pb:pomdp:MILP_pomdp} and its MDP approximation.
	Note that the two objective functions are the same. Hence, we only need to prove that we can construct a feasible solution from a problem to another.
	
	Let $(\bfmu,\bfdelta)$ be a feasible solution of the linear relaxation of Problem~\eqref{pb:pomdp:MILP_pomdp}. 
	Constraints~\eqref{eq:pomdp:NLP_consistency_sa}- ~\eqref{eq:pomdp:NLP_indep_state} ensure that $(\mu_s^1,\mu_{sas'}^t)_{t \in [T]}$ is a feasible solution of Problem~\eqref{pb:pomdp:LP_MDP}.

	Let $\bfmu$ be a feasible solution of Problem~\eqref{pb:pomdp:LP_MDP}. It suffices to define variables $\delta_{a|o}^t$ and $\mu_{soa}^t$ for all $a$ in $\calX_A$, $o$ in $\calX_O$, $s$ in $\calX_S$, and $t$ in $[T]$. We define these variables using ~\eqref{eq:pomdp:proof_mu_soa} and~\eqref{eq:pomdp:proof_delta}. In the proof of Proposition~\ref{prop:pomdp:valid_cuts_pomdp}, we proved that $(\bfmu,\bfdelta)$ is a feasible solution of the linear relaxation of MILP~\eqref{pb:pomdp:MILP_pomdp}.  
	Consequently, the equivalence holds and $z_{\rm{R}}^* = v_{\rm{MDP}}^*$.

	Now we prove that inequalities~\eqref{eq:pomdp:inequality_information} hold.
	Note that Proposition~\eqref{prop:pomdp:valid_cuts_pomdp} ensures that
	$$z^* \leq z_{\rm{R}^{\rm{c}}}^* \leq z_{\rm{R}}^*.$$
	It remains to prove the two following inequalities.
	\begin{align}
		&z^* \leq v_{\rm{his}}^* \label{eq:pomdp:ineq1}\\
		&v_{\rm{his}}^* \leq z_{\rm{R}^{\rm{c}}}^* \label{eq:pomdp:ineq2}
	\end{align}

	First, we prove Inequality~\eqref{eq:pomdp:ineq1}. By definition, we have $\Deltaml \subseteq \Deltahis$.
	Hence, we obtain $v_{\rm{ml}}^* \leq v_{\rm{his}}^*$. Using Theorem~\ref{theo:pomdp:NLP_optimal_solution}, we deduce that $z^* \leq v_{\rm{his}}^*$.
	Therefore the inequality $v_{\rm{ml}}^* \leq v_{\rm{his}}^* \leq z_{\rm{R}}^*$ holds.

	Now we prove Inequality~\eqref{eq:pomdp:ineq2}.
	The proof is based on a probabilistic interpretation of the valid inequalities~\eqref{eq:pomdp:Valid_cuts_pomdp}. It suffices to proves that for any policy $\bfdelta$ in $\Deltahis$, the probability distribution $\bbP_{\bfdelta}$ satisfies the weak conditional independences ~\eqref{eq:pomdp:weakIndep}.
	Let $\bfdelta \in \Deltahis$. The probability distribution $\bbP_{\bfdelta}$ over the random variables $(S_t,A_t,O_t)_{1\leq t \leq T}$ according to $\bfdelta$ is exactly
	\begin{align}\label{eq:pomdp:proof_distrib}
	\bbP_{\bfdelta} (\left(S_t=s_t,O_t=o_t,A_t=a_t \right)_{1\leq t \leq T}) &= \bbP_{\bfdelta}(S_1=s_1)\prod_{t=1}^T \bbP_{\bfdelta}(S_{t+1}=s_{t+1}|S_t=s_t,A_t=a_t) \notag \\
	&\bbP_{\bfdelta}(O_t=o_t|S_t=s_t) \delta^t_{a_t|h_t}
	\end{align}
	where $h_t = \{O_1=o_1,A_1=a_1,O_2=o_2,\ldots, O_t=o_t\}$ is the history of observations and actions. Note that the policy at time $t$ is the conditional probability $\delta^t_{a_t|h_t} = \bbP_{\bfdelta}(A_t=a_t|H_t=h_t)$.
	We define:
	\begin{align*}
		&\mu_{s}^1 = \bbP_{\bfdelta}(S_1=s)\\
		&\mu_{soa}^t = \bbP_{\bfdelta}(S_t=s,O_t=o,A_t=a)\\
		&\mu_{sas'}^t = \bbP_{\bfdelta}(S_t=s,A_t=a,S_{t+1} = s')\\
		&\mu_{s'a'soa}^t = \bbP_{\bfdelta}(S_{t-1}=s',A_{t-1}=a',S_t=s,O_t=o,A_t=a)
	\end{align*}
	We define the policy $\tilde{\bfdelta}$ using~\eqref{eq:pomdp:proof_delta}.
	It is easy to see that constraints of~\eqref{pb:pomdp:MILP_pomdp}
	are satisfied.
	Furthermore, we have $\tilde{\bfdelta} \in \Deltaml$.
	It remains to prove that inequalities~\eqref{eq:pomdp:Valid_cuts_pomdp} are satisfied. By definition of a probability distribution, we directly see that constraints~\eqref{eq:pomdp:Valid_cuts_pomdp_consistency1} are satisfied.
	We prove \eqref{eq:pomdp:Valid_cuts_pomdp_consistency2} and \eqref{eq:pomdp:Valid_cuts_pomdp_main}. We compute the left-hand side of \eqref{eq:pomdp:Valid_cuts_pomdp_consistency2}:
	\begin{align*}
		&\sum_{a \in \calX_A} \mu_{s'a'soa}^t = \sum_{a \in \calX_A} \bbP_{\bfdelta}(S_{t-1}=s',A_{t-1}=a',S_t=s,O_t=o,A_t=a)\\
		&= \sum_{a \in \calX_A} \sum_{\substack{s_1,\ldots,s_{t-2} \\ h_{t-1}}} \\
		&\bbP_{\bfdelta}((S_i=s_i,O_i=o_i,A_i=a_i)_{1\leq i\leq t-2},S_{t-1}=s',O_{t-1}=o',A_{t-1}=a',S_t=s,O_t=o,A_t=a)\\
		&= p(o|s)p(s|s',a')\sum_{\substack{s_1,\ldots,s_{t-2} \\ h_{t-1}}} \bbP_{\bfdelta}((S_i=s_i,O_i=o_i,A_i=a_i)_{1\leq i\leq t-2},S_{t-1}=s',O_{t-1}=o_{t-1},A_{t-1}=a') \\
		& \sum_{a \in \calX_A} \delta_{a|h_t} \\
		&= p(o|s)p(s|s',a')\sum_{\substack{s_1,\ldots,s_{t-2} \\ h_{t-1}}} \bbP_{\bfdelta}((S_i=s_i,O_i=o_i,A_i=a_i)_{1\leq i\leq t-2},S_{t-1}=s',O_{t-1}=o',A_{t-1}=a')\\
		&=  p(o|s)p(s|s',a') \bbP_{\bfdelta}(S_{t-1}=s',A_{t-1}=a')\\
		&=  p(o|s)\mu_{s'a's}^{t-1}
	\end{align*}
	where we used the definition of the probability distribution~\eqref{eq:pomdp:proof_distrib} at the third equation. Therefore, constraints~\eqref{eq:pomdp:Valid_cuts_pomdp_consistency2} are satisfied by $\bfmu$. 
	To prove that constraints~\eqref{eq:pomdp:Valid_cuts_pomdp_main} are satisfied, we prove that
	$$\bbP_{\bfdelta}(S_t=s_t| S_{t-1}=s_{t-1},A_{t-1}=a_{t-1},O_t=o_t,A_t=a_t) = \bbP_{\bfdelta}(S_t=s_t|  S_{t-1}=s_{t-1},A_{t-1}=a_{t-1},O_t=o_t)$$
	We compute $\bbP_{\bfdelta}(S_t=s_t| S_{t-1}=s',A_{t-1}=a',O_t=o,A_t=a)$:
	\begin{align*}
		&\bbP_{\bfdelta}(S_t=s_t| S_{t-1}=s_{t-1},A_{t-1}=a_{t-1},O_t=o_t,A_t=a_t) \\
		&= \frac{\bbP_{\bfdelta}(S_{t-1}=s_{t-1},A_{t-1}=a_{t-1},S_t=s_t,O_t=o_t,A_t=a_t)}{\bbP_{\bfdelta}(S_{t-1}=s_{t-1},A_{t-1}=a_{t-1},O_t=o_t,A_t=a_t)}\\
		&= \frac{\sum_{\substack{s_1,\ldots,s_{t-2} \\ h_{t-1}}} \bbP_{\bfdelta}((S_i=s_i,O_i=o_i,A_i=a_i)_{1\leq i\leq t})}{\sum_{\substack{s_1,\ldots,s_{t-2},s_t \\ h_{t-1}}} \bbP_{\bfdelta}((S_i=s_i,O_i=o_i,A_i=a_i)_{1\leq i\leq t})} \\
		&= \frac{\sum_{\substack{s_1,\ldots,s_{t-2} \\ h_{t-1}}} \delta_{a_t|h_t}^t p(o_t|s_t)p(s_{t}|s_{t-1},a_{t-1}) \bbP_{\bfdelta}((S_i=s_i,O_i=o_i,A_i=a_i)_{1\leq i\leq t-1})}{\sum_{\substack{s_1,\ldots,s_{t-2},s_t' \\ h_{t-1}}} \delta_{a_t|h_t}^t p(o_t|s_t)p(s_{t}|s_{t-1},a_{t-1}) \bbP_{\bfdelta}((S_i=s_i,O_i=o_i,A_i=a_i)_{1\leq i\leq t-1})}
	\end{align*}
	\begin{align*}
		&= \frac{p(o_t|s_t)p(s_{t}|s_{t-1},a_{t-1}) \sum_{\substack{s_1,\ldots,s_{t-2} \\ h_{t-1}}} \delta_{a_t|h_t}^t \bbP_{\bfdelta}((S_i=s_i,O_i=o_i,A_i=a_i)_{1\leq i\leq t-1})}{\sum_{s_t'} p(o_t|s_t')p(s_{t}'|s_{t-1},a_{t-1})\sum_{\substack{s_1,\ldots,s_{t-2} \\ h_{t-1}}} \delta_{a_t|h_t}^t \bbP_{\bfdelta}((S_i=s_i,O_i=o_i,A_i=a_i)_{1\leq i\leq t-1})} \\
		&= \frac{p(o_t|s_t)p(s_{t}|s_{t-1},a_{t-1})}{\sum_{s_t'} p(o_t|s_t')p(s_{t}'|s_{t-1},a_{t-1})}
	\end{align*}
	where the last line goes from the fact that the term $\delta_{a_t|h_t}^t \bbP_{\bfdelta}((S_i=s_i,O_i=o_i,A_i=a_i)_{1\leq i\leq t-1})$ does not depend on $s_t$.
	Hence, constraints~\eqref{eq:pomdp:Valid_cuts_pomdp_main} are satisfied by $\bfmu$. We deduce that $\bfmu$ is a feasible solution of MILP~\eqref{pb:pomdp:MILP_pomdp} satisfying the valid inequalities~\eqref{eq:pomdp:Valid_cuts_pomdp}.
	Therefore,
	\begin{align*}
		\bbE_{\bfdelta} \left[ \sum_{t=1}^T r(S_t,A_t,S_{t+1}) \right] &= \sum_{t=1}^T\sum_{s,a,s'}\bbP_{\bfdelta}(S_t=s,A_t=a,S_{t+1}=s')r(s,a,s') \leq z_{\rm{R}^{\rm{c}}}^*
	\end{align*} 
	By maximizing over $\bfdelta$ the left-hand side, we obtain $v_{\rm{his}}^* \leq z_{\rm{R}^{\rm{c}}}^*$. It achieves the proof.
\qed

\section{The rolling horizon matheuristic}
\label{app:wkpomdp:heuristic}

In this section, we provide more details on the belief state update used in Algorithm $\rmAct_{T}^{\rmp,t}(\bfh)$ (Section~\ref{sub:wkpomdp:first_policy}) and the rolling horizon matheuristic (Section~\ref{sub:wkpomdp:matheuristic})

\subsection{Belief state update}
\label{sub:app_wkpomdp_heuristic:belief_state}

Given a POMDP $\left(\calX_S,\calX_O,\calX_A,\pfrak,\bfr \right)$, at each time $t$, the belief state $(p(s|H_t))_{s \in \calX_S}$ is a sufficient statistic of the history of actions and observations $H_t$ \cite[Theorem 4]{Eckles1968}.
Given the action $a_t$ taken at time $t$ an the observation $o_{t+1}$ received at time $t+1$, the belief state can be easily updated over time according to the belief state update \citep[Eq. (1)]{Littman1994}:
\begin{align}\label{eq:belief_state_update}
	p(s_{t+1}|H_{t+1}) = p(s_{t+1}|H_{t}, a_t, o_{t+1}) = \sum_{s \in \calX_S} \frac{p(o_{t+1}|s)p(s|s_t, a_t)}{\sum_{s' \in \calX_S} p(o_{t+1}|s')p(s'|s_t, a_t)}p(s| H_t)
\end{align}

\subsection{Matheuristic with a rolling horizon}
\label{sub:app_wkpomdp_heuristic:rolling_horizon}

We denote by $T_{\rmr} < T$ denote such a rolling horizon. It leads to a matheuristic with formulation $\rmp$.

\begin{algorithm}[H]
\caption{Rolling horizon matheuristic with formulation $\rmp$.}
\label{alg:wkpomdp:matheuristic}
\begin{algorithmic}[1]
\STATE \textbf{Input}: $T$, $T_{\rm{r}}$, $\big((\calX_S^m,\calX_O^m,\calX_A^m,\pfrak^m,\bfr^m,\bfD^m)_{m\in[M]},\bfb\big)$
\FOR{$t=1,\ldots,T$}
\STATE Receive observation $\bfo_t$
\STATE Take action $\mathrm{Act}_{t+T_\rmr}^{\rmp,t}(\bfh_t)$ \label{alg:wkpomdp:matheuristic:take_action}
\ENDFOR
\end{algorithmic}
\end{algorithm}

Figure~\ref{fig:wkpomdp:scheme_implicite_policy} illustrates two consecutive iterations of Algorithm~\ref{alg:wkpomdp:matheuristic} with a rolling horizon $T_{\mathrm{r}}=5.$

\begin{figure}
	\begin{subfigure}{1.0\textwidth}
	\centering
	\resizebox{11cm}{!}{
    \begin{tikzpicture}
    
        \def\h{1.5}
        \def\l{1.5}
    
        \node (x1) at (-4*\l,0*\h) {};
        \node (x2) at (4*\l,0*\h) {};

        \draw[dashed] (-2.5*\l,0*\h) -- (-2.5*\l,1*\h);
        \draw[dashed] (-0.5*\l,0*\h) -- (-0.5*\l,1*\h);
        \draw[<->] (-2.5*\l,1*\h) -- node [text width=2.0cm,midway,above,align=center] {$T_{\mathrm{r}}=5$} (-0.5*\l,1*\h);
        \node[red] (txt) at (-2.5*\l,-0.5*\h) {Take action $\mathrm{Act}_{3+T_{\mathrm{r}}}^{\rmp,3}(\bfh_{3})$};
        \pattern[color=red,fill opacity=0.2] (-2.5*\l,0.0*\h)--(-2.5*\l,1*\h)--(-0.5*\l,1*\h)--(-0.5*\l,0*\h)--cycle;
        \pattern[pattern=north east lines] (-3.5*\l,0*\h)--(-3.5*\l,0.5*\h)--(-2.5*\l,0.5*\h)--(-2.5*\l,0*\h)--cycle;
        \node[black] (history) at (-3.0*\l,0.7*\h) {$\bfh_3$};

        \node (Time) at (4*\l,-0.2*\h) {Time};
        \draw[arc] (x1) -- (x2);
        \draw [black, fill] (-3.5*\l,0) circle [radius=0.08] node[below]{$t=1$};
        \draw [black, fill] (-3.0*\l,0) circle [radius=0.08];
        \draw [red, fill] (-2.5*\l,0) circle [radius=0.08] node[below]{$t=3$};

        \draw[black, fill] (-2.0*\l,0) circle [radius=0.08];
        \draw [black, fill] (-1.5*\l,0) circle [radius=0.08];
        \draw [black, fill] (-1.0*\l,0) circle [radius=0.08];
        \draw [black, fill] (-0.5*\l,0) circle [radius=0.08];
        \draw [black, fill] (0.0*\l,0) circle [radius=0.08];
        \draw [black, fill] (0.5*\l,0) circle [radius=0.08];
        \draw [black, fill] (1.0*\l,0) circle [radius=0.08];
        \draw [black, fill] (1.5*\l,0) circle [radius=0.08];
        \draw [black, fill] (2.0*\l,0) circle [radius=0.08];
        \draw [black, fill] (2.5*\l,0) circle [radius=0.08];
        \draw [black, fill] (3.0*\l,0) circle [radius=0.08] node[below]{$T$};


    \end{tikzpicture}}
    \caption{Iteration $t=3$ of Algorithm~\eqref{alg:wkpomdp:matheuristic}}
  	\label{fig:wkpomdp:decision_time_3} 
    \end{subfigure}
    \newline
    \begin{subfigure}{1.0\textwidth}
    \centering
    \resizebox{11cm}{!}{
    \begin{tikzpicture}
    
        \def\h{1.5}
        \def\l{1.5}
    
        \node (x1) at (-4*\l,0*\h) {};
        \node (x2) at (4*\l,0*\h) {};

        \draw[dashed] (-2.0*\l,0*\h) -- (-2.0*\l,1*\h);
        \draw[dashed] (0.0*\l,0*\h) -- (0.0*\l,1*\h);
        \draw[<->] (-2.0*\l,1*\h) -- node [text width=2.0cm,midway,above,align=center] {$T_{\mathrm{r}}=5$} (0.0*\l,1*\h);
        \node[red] (txt) at (-2.0*\l,-0.5*\h) {Take action $\mathrm{Act}_{4+T_{\mathrm{r}}}^{\rmp,4}(\bfh_{4})$};
        \pattern[color=red,fill opacity=0.2] (-2.0*\l,0.0*\h)--(-2.0*\l,1*\h)--(0.0*\l,1*\h)--(0.0*\l,0*\h)--cycle;
        \pattern[pattern=north east lines] (-3.5*\l,0*\h)--(-3.5*\l,0.5*\h)--(-2.0*\l,0.5*\h)--(-2.0*\l,0*\h)--cycle;
        \node[black] (history) at (-2.75*\l,0.7*\h) {$\bfh_4$};

        \node (Time) at (4*\l,-0.2*\h) {Time};
        \draw[arc] (x1) -- (x2);
        \draw [black, fill] (-3.5*\l,0) circle [radius=0.08] node[below]{$t=1$};
        \draw [black, fill] (-3.0*\l,0) circle [radius=0.08];
        \draw [black, fill] (-2.5*\l,0) circle [radius=0.08];

        \draw[red, fill] (-2.0*\l,0) circle [radius=0.08] node[below]{$t=4$};
        \draw [black, fill] (-1.5*\l,0) circle [radius=0.08];
        \draw [black, fill] (-1.0*\l,0) circle [radius=0.08];
        \draw [black, fill] (-0.5*\l,0) circle [radius=0.08];
        \draw [black, fill] (0.0*\l,0) circle [radius=0.08];
        \draw [black, fill] (0.5*\l,0) circle [radius=0.08];
        \draw [black, fill] (1.0*\l,0) circle [radius=0.08];
        \draw [black, fill] (1.5*\l,0) circle [radius=0.08];
        \draw [black, fill] (2.0*\l,0) circle [radius=0.08];
        \draw [black, fill] (2.5*\l,0) circle [radius=0.08];
        \draw [black, fill] (3.0*\l,0) circle [radius=0.08] node[below]{$T$};
    \end{tikzpicture}}
    \caption{Iteration $t=4$ of Algorithm~\ref{alg:wkpomdp:matheuristic}}
  	\label{fig:wkpomdp:decision_time_4}    
  	\end{subfigure}
    \caption{Scheme of the evaluation of our implicit policy $\bfdelta^{\rmIP}$ in Algorithm~\ref{alg:wkpomdp:matheuristic} at time $t=3$  and $t=4.$ The decision maker observes $\bfo_{3}$ and takes action $\mathrm{Act}_{8}^{\rmp,3}(\bfh_{3}).$ Then, the decision maker observes $\bfo_{4}$ and takes action $\mathrm{Act}_{9}^{\rmp,4}(\bfh_{4}).$  The black points indicate the time steps and the red point corresponds to the time when the decision is taken. The black hatched lines represent the past at the current time (red). The red square indicates the horizon taken into account in the optimization problem.}
    \label{fig:wkpomdp:scheme_implicite_policy}
\end{figure}

\section{Proofs of Section~\ref{sec:wkpomdp}}
\label{app:wkpomdp:proofs}

\proof[Proof of Theorem~\ref{theo:wkpomdp:ineq_decPOMDP}]
	First, we prove that the linear relaxation of MILP~\eqref{pb:wkpomdp:decPOMDP_MILP} is a relaxation of the MDP approximation. Let $\bfmu$ be a feasible solution of the linear program~\eqref{pb:pomdp:LP_MDP}.
	The variable $\mu_{\bfs\bfa\bfs'}^t$ is defined for $\bfa \in \calX_A = \{ \bfa \in \calX_A^1\times\cdots\calX_A^M \colon \sum_{m=1}^M \bfD^m(a^m) \leq \bfb \}$. We extend its definition on $\calX_A' = \calX_A^1\times\cdots\calX_A^M$ by setting $\mu_{\bfs\bfa\bfs'}^t:=0$ for $\bfa \in \calX_A' \backslash \calX_A$.
	Now we define the variables $(\tau_{s}^m,\tau_{soa}^{t,m},\tau_{sas'}^{t,m})_{s,o,a,t,m}$ as follows
	\begin{align*}
		&\tau_s^{1,m} = \sum_{s^{-m} \in \calX_S^{-m}} \mu_{\bfs}^1,& \quad &\tau_{sas'}^{t,m} = \sum_{\substack{s^{-m} \in \calX_S^{-m} \\ a^{-m} \in \calX_A^{-m} \\ s'^{-m} \in \calX_S^{-m}}} \mu_{\bfs\bfa\bfs'}^{t},
	\end{align*}
	and, for each component $m$, we define the variables $\delta_{a|o}^{t,m}$ and $\tau_{soa}^{t,m}$ using~\eqref{eq:pomdp:proof_mu_soa} and~\eqref{eq:pomdp:proof_delta}. We proved that this solution is a feasible solution of the linear relaxation on each component. 
	It reamains to prove that $\sum_{m=1}^M \sum_{a \in \calX_A^m} \tau_{a}^{t,m} \bfD^m(a) \leq \bfb$. It comes from the following computation:
	\begin{align*}
		\sum_{m=1}^M \sum_{a \in \calX_A^m} \tau_{a}^{t,m} \bfD^m(a) &= \sum_{m=1}^M \sum_{\bfa \in \calX_A'} \mu_{\bfa}^t\bfD^m(a^m) \\
		&= \sum_{\bfa \in \calX_A'} \mu_{\bfa}^t \overbrace{\sum_{m=1}^M\bfD^m(a^m)}^{\leq \bfb} \leq \bfb
	\end{align*}
	Therefore, $(\bftau^m,\bfdelta^m)_{m\in [M]}$ is a feasible solution of the linear relaxation of MILP~\eqref{pb:wkpomdp:decPOMDP_MILP}.
	In addition, the objective value induced by variable $\bfmu$ is 
	\begin{align*}
		\sum_{t=1}^T \sum_{\bfs,\bfs',\bfa} \mu_{\bfs\bfa\bfs'}^t\sum_{m=1}^M r^m(s^m,a^m,s'^m) &= \sum_{t=1}^T \sum_{m=1}^M \sum_{\bfs,\bfs',\bfa} \mu_{\bfs\bfa\bfs'}^t r^m(s^m,a^m,s'^m) \\
		&= \sum_{t=1}^T \sum_{m=1}^M \sum_{s^m,a^m,s'^m} r^m(s^m,a^m,s'^m) \leq z_{\rmR}.
	\end{align*}
	It shows that $v_{\rmMDP}^* \leq z_{\rmR}$.

	Now we prove that $v_{\rmhis}^* \leq z_{\rmRc}$. To prove this result, it suffices to observe that the linear relaxation of MILP~\eqref{pb:wkpomdp:decPOMDP_MILP} with valid inequalities~\eqref{eq:wkpomdp:dec_Valid_cuts} is a relaxation of the linear relaxation of MILP~\eqref{pb:pomdp:MILP_pomdp} with valid inequalities~\eqref{eq:pomdp:Valid_cuts_pomdp}. 
	This result appears using the previous technique that defines the ``local'' marginals $\tau_{sas'o'a}^{t,m} = \sum_{s^{-m}a^{-m}s'^{-m}o'^{-m}a^{-m}} \mu_{\bfs\bfa\bfs'\bfo'\bfa'}^t$ on each component $m \in [M]$.
	Therefore, we obtain that $z_{\rmRc}^* \leq z_{\rmRc}$.
	On other hand, Theorem~\ref{theo:pomdp:MDP_approx_equivalence} ensures that $v_{\rmhis}^* \leq z_{\rmRc}^*$, which achieves the proof.    
\qed

\proof[Proof of Theorem~\ref{theo:wkpomdp:lower_bound_MILP}]
	Let $(\bftau^m,\bfdelta^m)_{m\in [M]}$ be a feasible solution of MILP~\eqref{pb:wkpomdp:decPOMDP_MILP_LB}.
	We prove that $(\bftau^m,\bfdelta^m)_{m\in [M]}$ is a feasible solution of MILP~\eqref{pb:wkpomdp:decPOMDP_MILP} and~\ref{pb:decPOMDP_wc}.
	First, we show that $(\bftau^m,\bfdelta^m)_{m\in [M]}$ is a feasible solution of MILP~\eqref{pb:wkpomdp:decPOMDP_MILP}. 
	We define the variables $\tau_a^{t,m}$ for any $a \in \calX_A^m,$ $m \in [M],$ and $t \in  [T]$ such that $\sum_{s \in \calX_S^m, o \in \calX_O^m} \tau_{soa}^{t,m} = \tau_a^{t,m}.$ In addition, we introduce the variables $\tau_o^{t,m} = \sum_{s \in \calX_S^m, o \in \calX_O^m} \tau_{soa}^{t,m}$ for any $o \in \calX_O^m,$ any $m \in [M]$ and $t \in [T].$
	It suffices to show that inequality~\eqref{eq:wkpomdp:decPOMDP_MILP_linking_cons} holds. 
	We compute the left-hand side of~\eqref{eq:wkpomdp:decPOMDP_MILP_linking_cons}.
	\begin{align*}
		\sum_{m =1}^M \sum_{a \in \calX_A^m} \bfD^m(a)\tau_{a}^{t,m} = \sum_{m =1}^M\sum_{a \in \calX_A^m,o \in \calX_O^m} \bfD^m(a) \delta_{a|o}^{t,m} \tau_{o}^{t,m} = \sum_{m =1}^M \sum_{a \in \calX_A^m} \bfD^m(a^m) \bbE_{\tau^{t,m}} [\delta_{a^m|O_t^m}^{t,m}] \leq \bfb 
	\end{align*}
	The first equality is a consequence of the tightness of the McCormick constraints~\eqref{eq:pomdp:MILP_McCormick_1}-\eqref{eq:pomdp:MILP_McCormick_3}. The second equality comes from the fact that the variables $(\tau_{o}^{t,m})_{o \in \calX_O^m}$ define a probability distribution over $\calX_O^m.$
	Finally, the last inequality results from 
	$$\sum_{m=1}^M \sum_{a \in \calX_A^m} \bfD^m(a) \bbE_{\tau^{t,m}} [\delta_{a|O_t^m}^{t,m}] \leq \sum_{m =1}^M \sum_{a \in \calX_A^m} \bfD^m(a)\max_{o \in \calX_O^m}(\delta_{a|o}^{t,m}) \leq \bfb$$ 
	Therefore, the inequality holds~\eqref{eq:wkpomdp:decPOMDP_MILP_linking_cons} and MILP~\eqref{pb:wkpomdp:decPOMDP_MILP} is a relaxation of MILP~\eqref{pb:wkpomdp:decPOMDP_MILP_LB}.

	Second, we show that $(\bftau^m,\bfdelta^m)_{m\in [M]}$ is a feasible solution of~\ref{pb:decPOMDP_wc}. We define a policy over $\calX_A \times \calX_O$.
	\begin{align}
		\delta_{\bfa|\bfo}^t = \prod_{m=1}^M \delta_{a^m|o^m}^{t,m}
	\end{align}
	for all $\bfa \in \calX_A$, $\bfo \in \calX_O$ and $t \in [T]$. It suffices to prove that $\bfdelta$ belongs to $\Delta$. Let $\bfo \in \calX_O$ and $t \in [T]$.
	\begin{align*}
	 	\sum_{\bfa \in \calX_A} \delta_{\bfa | \bfo}^t &= \sum_{\bfa \in \calX_A} \prod_{m=1}^M \delta_{a^m|o^m}^{t,m} = \sum_{\bfa \in \calX_A} \prod_{m=1}^M \delta_{a^m|o^m}^{t,m} =  1
	 \end{align*}
	 The second equality comes from the fact that for any $\bfa \in \calX_A$ such that $\sum_{m =1}^M \bfD^m(a^m) > \bfb,$ $\prod_{m=1}^M \delta_{a^m|o^m}^{t,m} = 0$ because of Constraints~\eqref{eq:wkpomdp:decPOMDP_link_constraints_Lb}. Therefore, $\bfdelta$ is a feasible policy of~\ref{pb:decPOMDP_wc}.

	 Since $\bfdelta= \prod_{m=1}^M \bfdelta^m$, then all the component are independents and the marginal probabilities $\bftau^m$ are exact in the sense that they derive from policy $\bfdelta$. 
	 It follows that the objective functions are the same and the inequalities $z_{\rmLB} \leq v_{\rmml}^{*}$ and $z_{\rmLB} \leq z_{\rmIP}$ hold.
\qed

\proof[Proof of Theorem~\ref{theo:wkpomdp:upper_bound_NLP}]
	First, we prove that the nonlinear Program~\eqref{pb:wkpomdp:decPOMDP_NLP} is a relaxation of MILP~\eqref{pb:wkpomdp:decPOMDP_MILP}.
	Let $(\bftau^m,\bfdelta^m)_{m\in [M]}$ be a feasible solution of MILP~\eqref{pb:wkpomdp:decPOMDP_MILP}.
	We prove that $(\bftau^m,\bfdelta^m)_{m\in [M]}$ is a feasible solution of the nonlinear Program~\eqref{pb:wkpomdp:decPOMDP_NLP}.
	It suffices to prove that for all $m \in [M],$ $(\bftau^m,\bfdelta^m)$ satisfies constraints~\eqref{eq:pomdp:NLP_indep_action}. It comes from the tightness of the McCormick inequalities~\eqref{eq:pomdp:MILP_McCormick_1}-\eqref{eq:pomdp:MILP_McCormick_3} when the policy is deterministic.
	Hence, it is a relaxation with the same objective function. Therefore, the inequality $z_{\rmIP} \leq z_{\rmUB}$ holds.

	Second, we prove that the nonlinear Program~\eqref{pb:wkpomdp:decPOMDP_NLP} is a relaxation of~\ref{pb:decPOMDP_wc}.
	Let $\bfdelta$ be a feasible policy of~\ref{pb:decPOMDP_wc}.
	We want to define a solution of the non-linear program~\eqref{pb:wkpomdp:decPOMDP_NLP}. We extend the domain of $\bfdelta$ to $\calX_A$ by setting $\delta_{\bfa|\bfo}^t = 0$ when $\sum_{m \in [M]} \bfD^m(a^m) > \bfb,$ for all $\bfo \in \calX_O.$ It is easy to see that $\bfdelta$ is still a policy in $\calX_A.$
	Theorem~\ref{theo:pomdp:NLP_optimal_solution} ensures that there exist $\bfmu$ such that $(\bfmu,\bfdelta)$ is a feasible solution of MILP~\eqref{pb:pomdp:MILP_pomdp}.
	We define the variables $\bftau^m$ and $\bfdelta^m$ on component $m \in [M]$ by induction
	\begin{align*}
		&\tau_s^{1,m} = \sum_{s^{-m} \in \calX_S^{-m}} \mu_{\bfs}^1,& \quad &\tau_{soa}^{t,m} = \sum_{\substack{s^{-m} \in \calX_S^{-m} \\ o^{-m} \in \calX_O^{-m} \\ a^{-m} \in \calX_A^{-m}}} \mu_{\bfs\bfo\bfa}^t, \\
		&\tau_{sas'}^{t,m} = \sum_{\substack{s^{-m} \in \calX_S^{-m} \\ a^{-m} \in \calX_A^{-m} \\ s'^{-m} \in \calX_S^{-m}}} \mu_{\bfs\bfa\bfs'}^{t}, & \quad &
		\delta_{a|o}^{t,m} = \sum_{\substack{s^{-m} \in \calX_S^{-m} \\ o^{-m} \in \calX_O^{-m} \\ a^{-m} \in \calX_A^{-m}}} \delta^t_{\bfa|\bfo} \prod_{m' \neq m} p^{m'}(o^{m'}|s^{m'})\tau_{s^{m'}}^{t,m'}
	\end{align*}
	for all $s \in \calX_S^m$, $o \in \calX_O^m$, $a \in \calX_A^m$ and $t \in [T].$
	By definition of $\bftau^m,$ if $\bfdelta^m$ is in $\Delta^m,$ then the constraints~\eqref{eq:wkpomdp:NLP_per_POMDP_v1} are satisfied by $(\bftau^m,\bfdelta^m)_{m\in [M]}.$
	We prove that $\bfdelta^m$ is in $\Delta_{\mathrm{ml}}^m.$
	\begin{align*}
		\sum_{a \in \calX_A^m} \delta_{a|o}^{t,m} &= \sum_{a \in \calX_A^m} \sum_{\substack{s^{-m} \in \calX_S^{-m} \\ o^{-m} \in \calX_O^{-m} \\ a^{-m} \in  \calX_A^{-m}}} \delta^t_{\bfa|\bfo} \prod_{m' \neq m} p^{m'}(o^{m'}|s^{m'})\tau_{s^{m'}}^{t,m'} = \sum_{\substack{s^{-m} \in \calX_S^{-m} \\ o^{-m} \in \calX_O^{-m}}} \prod_{m' \neq m} p^{m'}(o^{m'}|s^{m'})\tau_{s^{m'}}^{t,m'} =1
	\end{align*}
	for all $o \in \calX_O^m,$ $m \in [M]$ and $t \in [T].$ The last equality comes from the fact that by induction we have that $\sum_{s \in \calX_S^m} \tau_s^{t,m} = 1.$
	Therefore, $\bfdelta^m \in \Deltaml^m.$

	It remains to prove that constraints~\eqref{eq:wkpomdp:decPOMDP_NLP_linking_cons} are satisfied by $(\bftau^m,\bfdelta^m)_{m\in [M]}.$ We compute the left-hand side of constraint~\eqref{eq:wkpomdp:decPOMDP_NLP_linking_cons}.
	\begin{align*}
		\sum_{m=1}^M \sum_{a \in \calX_A^m} \bfD^m(a) \tau_{a}^{t,m} 
		&= \sum_{m=1}^M \sum_{a \in \calX_A^m} \bfD^m(a)\bigg(\sum_{\substack{\bfs \in \calX_S, \bfo \in \calX_O \\ \bfa' \in \calX_A:a'^m=a^m}} \delta^t_{\bfa'|\bfo} \prod_{m' =1}^M p^{m'}(o^{m'}|s^{m'})\tau_{s^{m'}}^{t,m'}\bigg)  \\
		&= \sum_{\substack{\bfs \in \calX_S, \bfo \in \calX_O \\ \bfa' \in \calX_A}} \delta^t_{\bfa'|\bfo} \prod_{m'=1}^M p^{m'}(o^{m'}|s^{m'})\tau_{s^{m'}}^{t,m'} \sum_{m =1}^M \sum_{a \in \calX_A^m} \bfD^m(a^m)\\
		&= \overbrace{\sum_{\substack{\bfs \in \calX_S, \bfo \in \calX_O \\ \bfa \in \calX_A}} \delta^t_{\bfa|\bfo} \prod_{m'=1}^M p^{m'}(o^{m'}|s^{m'})\tau_{s^{m'}}^{t,m'}}^{=1} \underbrace{\sum_{m=1}^M \bfD^m(a^m)}_{\leq \bfb}\\
		&\leq \bfb
	\end{align*}
	Therefore, constraints~\eqref{eq:wkpomdp:decPOMDP_NLP_linking_cons} are satisfied. Consequently, the nonlinear program~\eqref{pb:wkpomdp:decPOMDP_NLP} is a relaxation of~\ref{pb:decPOMDP_wc}.
	In addition, the objective functions are equal. We deduce that $v_{\rmml}^* \leq z_{\rmUB}.$
\qed

\proof[Proof of Proposition~\ref{prop:wkpomdp:decPOMDP_valid_cuts}]
	Proposition~\ref{prop:pomdp:valid_cuts_pomdp} ensures that inequalities~\eqref{eq:wkpomdp:dec_Valid_cuts} are valid on each component. Hence, these inequalities are valid for MILP~\eqref{pb:wkpomdp:decPOMDP_MILP}, MILP~\eqref{pb:wkpomdp:decPOMDP_MILP_LB} and Problem~\eqref{pb:wkpomdp:decPOMDP_NLP}.
	Proposition~\ref{prop:pomdp:valid_cuts_pomdp} also ensures that there are solutions of the linear relaxation of~\eqref{pb:pomdp:MILP_pomdp} that do not satisfy constraints~\eqref{eq:pomdp:Valid_cuts_pomdp} on each component.
\qed

\proof[Proof of Theorem~\ref{theo:wkpomdp:heuristic}]
	First, we prove that at each time $t \in [T]$, for every observation $\bfh\in\calX_H^t$ the element $\rm{Act}^t(\bfh)$ belongs to $\calX_A,$ i.e., $\sum_{m=1}^M \bfD^m\left(\rm{Act}^{t,m}(\bfh)\right) \leq \bfb.$
	Let $(\bftau^m,\bfdelta^m)_{m \in [M]}$ be a feasible solution of MILP~\eqref{pb:wkpomdp:decPOMDP_MILP} at step~\ref{alg:wkpomdp:modify_constraints}.
	Since $O_t^m= o_t^m$ almost surely, $\tau_{soa}^{t,m}$ is equal to $0$ when $o \neq o_t^m$.
	Hence, we obtain
	\begin{align*}
		\tau_a^{t,m} &= \sum_{s \in \calX_S^m, o \in \calX_O^m} \tau_{soa}^{t,m} = \sum_{s \in \calX_S^m} \tau_{s o_t^{m}a}^{t,m} = \delta_{a|o_t^m}^{t,m}
	\end{align*}
	It ensures that $\tau_{a}^{t,m} \in \{0,1\}$ for any $a \in \calX_A^m$ and $m \in [M]$. 
	Let $\bfa^*$ be the action taken at step~\ref{alg:wkpomdp:take_action}. Therefore, $\tau_{a}^{t,m} = 1$ when $a = a^{m*}$ and $0$ otherwise.
	Now we compute the linking constraint of~\eqref{eq:problem:def_form_action_space}.
	\begin{align*}
		\sum_{m =1}^M \bfD^m(a^{m*}) = \sum_{m =1}^M \bfD^m(a^{m*})\tau_{a^{m*}}^{t,m} =  \sum_{m=1}^M\sum_{a \in \calX_A^m} \bfD^m(a)\tau_{a}^{t,m} \leq \bfb
	\end{align*}
	The last inequality comes from the fact that $(\bftau^m,\bfdelta^m)_{m \in [M]}$ satisfies constraint~\eqref{eq:wkpomdp:decPOMDP_MILP_linking_cons}.

	Now we prove the inequalities.
	The inequality $\nu_{\rmIP}\leq v_{\rmhis}^*$ holds because $\bfdelta^{\rmIP} \in \Delta_{\rmhis}$. 
	It remains to show that $z_{\rmIP} \leq \nu_{\rmIP}$. We do it using a backward induction.
	Let $(\bftau^{*m},\bfdelta^{*m})_{m \in [M]}$ be an optimal solution of Problem~\eqref{pb:wkpomdp:decPOMDP_MILP}.
	We denote by $\rmP_t(\bfh_t)$ the feasible set of the optimization problem solved in $\rmAct_T^{\rmIP,t}(\bfh_t)$, for every $t$ in $[T]$.
	We consider the following induction hypothesis at time $t$: 
	$$\max_{(\bftau^m,\bfdelta^m)_{m\in [M]} \in \rmP_tb(\bfh_t)}\sum_{m=1}^M\bbE_{\bfdelta^{m}}\left[\sum_{t'=t}^T r^m(S_t^m,A_t^m,S_{t+1}^m) | H_{t}^m=h_{t}^m \right] \leq \bbE_{\bfdelta^{\rmIP}}\left[ \sum_{m=1}^M\sum_{t'=t}^T r^m(S_t^m,A_t^m,S_{t+1}^m) | \bfH_{t}=\bfh_{t} \right]$$
	If $t=T$, then consider left-hand side is exactly equal to the right-hand side.
	\begin{align*}
		&\max_{(\bftau^m,\bfdelta^m)_{m\in [M]} \in \rmP_T(\bfh_T) } \sum_{m=1}^M\bbE_{\bfdelta^m}\left[ r(S_T^m,A_T^m,S_{T+1}^m) | H_T^m=h_T^m \right] \\
		&= \max_{(\bftau^m,\bfdelta^m)_{m\in [M]} \in  \rmP_T(\bfh_T) } \sum_{m=1}^M\bbE_{\bfdelta^m}\left[ r(S_T^m,A_T^m,S_{T+1}^m) | S_T^m \sim p^m(\cdot|h_T^m) \right] = \bbE_{\bfdelta^{m,\rmIP}} \left[\sum_{m=1}^M r^m(S_T^m,A_T^m,S_{T+1}^m) | \bfH_{T}=\bfh_{T}\right]
	\end{align*}
	The first equality comes from the fact that the belief state is a sufficient statistic of the history.
	It proves the induction hypothesis for $t=T$.

	Suppose that the induction hypothesis holds from $t+1$.
	We compute the term in $t$:
	\begin{align*}
		&\max_{(\bftau^m,\bfdelta^m)_{m \in [M]} \in \rmP_t(\bfh_t)} \sum_{m=1}^M \bbE_{\bfdelta^{m}}\left[\sum_{t'=t}^T r^m(S_{t'}^m,A_{t'}^m,S_{t'+1}^m) | H_{t}^m=h_{t}^m\right] \\
		&= \max_{(\bftau^m,\bfdelta^m)_{m \in [M]} \in \rmP_t(\bfh_t)} \sum_{m=1}^M \bbE_{\bfdelta^{t,m}} \left[r^m(S_{t}^m,A_{t}^m,S_{t+1}^m) |H_t^m = h_t^m \right] \\
		&+ \sum_{m=1}^M \sum_{a_t^m,o_{t+1}^m} \Bigg(\bbP_{\bfdelta{t,m}}\left(O_{t+1}^m=o_{t+1}^m,A_{t}^m=a_t^m|H_t^m=h_t^m \right) \\
		&\times \overbrace{\bbE_{\bfdelta^m} \left[\sum_{t'=t+1}^T r^m(S_{t'}^m,A_{t'}^m,S_{t'+1}^m) |\underbrace{H_{t}^m = h_{t}^m, A_t^m=a_t^m,O_{t+1}^m=o_{t+1}^m}_{H_{t+1}^m=h_{t+1}^m} \right]}^{\text{does not depend on } \bfdelta^{t,m}} \Bigg)\\
		&\leq \bbE_{\bfdelta^{t,\rmIP}} \left[\sum_{m=1}^M r^m(S_{t}^m,A_{t}^m,S_{t+1}^m) |\bfH_t = \bfh_t \right] + \bigg(\sum_{\bfa_t,\bfo_{t+1}}\bbP_{\bfdelta{t,\rmIP}}\left(\bfO_{t+1}=\bfo_{t+1},\bfA_{t}=\bfa_t|\bfH_t=\bfh_t \right) \\
		&\times\overbrace{\max_{(\bftau^m,\bfdelta^m) \in \rmP_{t+1}(\bfh_{t+1})} \sum_{m=1}^M \bbE_{\bfdelta^m} \left[\sum_{t'=t+1}^T r^m(S_{t'}^m,A_{t'}^m,S_{t'+1}^m) |H_{t+1}^m = h_{t+1}^m \right]}^{\text{induction  hypothesis}} \bigg) \\
		&\leq \sum_{m=1}^M \bbE_{\bfdelta^{m,\rmIP}}\left[ r^m(S_t^m,A_t^m,S_{t+1}^m) | H_{t}^m=h_{t}^m\right]
	\end{align*}
	The first inequality above comes from the fact that there exists an optimal solution where $\bfdelta^{t,\rmIP}$ is the policy at time $t$ by definition of $\bfdelta^{\rmIP}$ and by decomposing the maximum operator in the sum of the second term. This latter operation can be done since $\bbE_{\bfdelta^m} \left[\sum_{t'=t+1}^T r^m(S_{t'}^m,A_{t'}^m,S_{t'+1}^m) |H_{t+1}^m = h_{t+1}^m \right]$ does not depend on the policy $\bfdelta^{t',m}$ for $t' < t+1$.  
	It proves the backward induction.
	Finally, given an optimal feasible solution $(\bftau*^m,\bfdelta*^m)_{m\in [M]}$ of MILP~\eqref{pb:wkpomdp:decPOMDP_MILP}, we get that:
	\begin{align*}
	z_{\rmIP} = \sum_{m=1}^M \bbE_{\bfdelta*^m}\left[ \sum_{t=1} r^m(S_t^m,A_t^m,S_{t+1}^m) \right] &= \bbE\left[ \bbE_{\bfdelta*^m}\left[ \sum_{m=1}^M \sum_{t=1}^T r^m(S_t^m,A_t^m,S_{t+1}^m) |O_1^m=o_1^m \right] \right] \\
	&\leq \bbE\left[ \max_{(\bftau^m,\bfdelta^m)_{m \in [M]} \in \rmIP_1(\bfO_1)} \sum_{m=1}^M \bbE_{\bfdelta^m}\left[ \sum_{t=1} r^m(S_t^m,A_t^m,S_{t+1}^m) |O_1^m=o_1^m \right] \right]\\
	&\leq \sum_{m=1}^M \bbE_{\bfdelta^{m,\rmIP}}\left[ \sum_{t=1} r^m(S_t^m,A_t^m,S_{t+1}^m) \right] = \nu_{\rmIP}
	\end{align*}
	The first inequality comes from the inversion of the maximum operator and the expectation. It achieves the proof.
\qed

\section{Lagrangian relaxations and column generation approach}
\label{app:ColGen}

In this section, we detail the proof of Proposition~\ref{prop:wkpomdp:lagrangian_relax} and we explain how we compute the value of the Lagrangian relaxation $z_{\rmLR}$ using a column generation approach.

\subsection{Proof of Proposition~\ref{prop:wkpomdp:lagrangian_relax}}
\label{sub:app_ColGen:Lagrangian_relax}

We denote by $\bfbeta=(\beta^t)_{t \in [T]}$ the dual variables associated with constraints~\eqref{eq:wkpomdp:decPOMDP_NLP_linking_cons}.
If we relax constraints~\eqref{eq:wkpomdp:decPOMDP_NLP_linking_cons}, then we obtain the Lagrangian function
\begin{align}\label{eq:app_ColGen:lagrangian_function}
 	\calL\left( \bftau,\bfdelta,\bfbeta \right) = \sum_{t=1}^T\sum_{m=1}^M \sum_{\substack{s,s' \in \calX_S^m \\ a \in \calX_A^m}} r^m(s,a,s') \tau_{sas'}^{t,m} + \sum_{t=1}^T  (\beta^t)^{\bfT} \left(\bfb - \sum_{m=1}^M \sum_{a \in \calX_A^m}\bfD^m(a)\tau_a^{t,m} \right),
 \end{align}
for any $(\bftau^m,\bfdelta^m)_{m \in [M]}.$
Then, we introduce the dual function $\calG : \bbR_{+}^{T\times q} \rightarrow \bbR,$ with values
\begin{subequations}\label{pb:app_ColGen:dual_function}
 \begin{alignat}{2}
 \calG \left(\bfbeta \right):= \max_{\bftau,\bfdelta} \enskip & \sum_{t=1}^T \sum_{m=1}^M \sum_{\substack{s,s' \in \calX_S^m \\ a \in \calX_A^m}} r^m(s,a,s')\tau_{sas'}^{t,m} + \sum_{t=1}^T  (\beta^t)^{\bfT} \left(\bfb - \sum_{m=1}^M \sum_{a \in \calX_A^m}\bfD^m(a)\tau_a^{t,m} \right) & \\
\mathrm{s.t.} \enskip 
 & \left(\bftau^m,\bfdelta^m\right) \in \calQ\left(T,\calX_S^m,\calX_O^m,\calX_A^m,\pfrak^m \right) \quad \forall m \in [M]
 \end{alignat}
\end{subequations}
By weak duality, for any $\bfbeta,$ the dual function~\eqref{pb:app_ColGen:dual_function} provides an upper bound obtained by using Approximation (A).

We now explain how to compute the dual function. 
As it is usually the case for Lagrangian relaxation, for every $\bfbeta \in \bbR_{+}^{T\times q},$ the maximum in the computation of $\calG(\bfbeta)$ decomposes over the sum of the maximum over each component.
However, the formulations obtained for each component are still nonlinear. 
Fortunately, the following proposition ensures that we can linearize the formulation without changing the value of the dual function.
\begin{prop}\label{prop:wkpomdp:dual_function}
	For all $\bfbeta \in \bbR_{+}^{T\times q},$ the dual function can be written as
	\begin{align}\label{pb:wkpomdp:MILP_dual_function}
		\calG\left( \bfbeta \right) = \sum_{t=1}^T (\beta^t)^{\bfT} \bfb + \sum_{m=1}^M \calG^m(\bfbeta) 
	\end{align}
	where $\calG^m(\bfbeta)$ is the quantity
	\begin{subequations}
		 \begin{alignat*}{2}
		 \calG^m \left(\bfbeta \right):= \max_{\bftau^m,\bfdelta^m} \enskip & \sum_{t=1}^T \sum_{\substack{s,s' \in \calX_S^m \\ a \in \calX_A^m}} \left(r^m(s,a,s') - (\bfbeta^t)^{\bfT}\bfD^m(a) \right)\tau_{sas'}^{t,m} & \quad &\\
		\mathrm{s.t.} \enskip 
		 & \left(\bftau^m,\bfdelta^m\right) \in \calQ^{\mathrm{d}}\left(T,\calX_S^m,\calX_O^m,\calX_A^m,\pfrak^m \right) &
		 \end{alignat*}
	\end{subequations}
\end{prop}

It follows from Proposition~\ref{prop:wkpomdp:dual_function} that the dual function can be computed by solving MILP~\eqref{pb:pomdp:MILP_pomdp} on each component of the system, which is in general easier than solving Problem~\eqref{pb:wkpomdp:decPOMDP_NLP}.

\proof[Proof of Proposition~\ref{prop:wkpomdp:dual_function}]
	Let $\bfbeta \in \bbR_{+}^{Tq}.$ Then, the value function $\calG$ in $\bfbeta$ can be written:
	\begin{subequations}
		 \begin{alignat*}{2}
		 \calG \left(\bfbeta \right)= \max_{\left(\bftau^m,\bfdelta^m\right)_{m\in [M]}} \enskip & \sum_{t=1}^T \sum_{m=1}^M \sum_{\substack{s,s' \in \calX_S^m \\ a \in \calX_A^m}} \left(r^m(s,a,s') - (\bfbeta^t)^{\bfT}\bfD^m(a) \right)\tau_{sas'}^{t,m} +  \sum_{t=1}^T (\beta^t)^{\bfT} \bfb & \quad &\\
		\mathrm{s.t.} \enskip 
		 & \left(\bftau^m,\bfdelta^m\right) \in \calQ\left(T,\calX_S^m,\calX_O^m,\calX_A^m,\pfrak^m \right) & \quad \forall m \in [M] 
		 \end{alignat*}
	\end{subequations}
	Since the second term does not depend on $(\bftau^m,\bfdelta^m)_{m\in [M]},$ we only consider the maximization on the first term. 
	In such a problem, there are no linking constraints between the components, which enables to decompose the maximization operator along the components as follows.
	\begin{subequations}
		 \begin{alignat*}{2}
		 \calG \left(\bfbeta \right)=  \sum_{m=1}^M\max_{\bftau^m,\bfdelta^m} \enskip & \sum_{t=1}^T \sum_{\substack{s,s' \in \calX_S^m \\ a \in \calX_A^m}} \left(r^m(s,a,s') - (\bfbeta^t)^{\bfT}\bfD^m(a) \right)\tau_{sas'}^{t,m} +  \sum_{t=1}^T (\beta^t)^{\bfT} \bfb & \quad &\\
		\mathrm{s.t.} \enskip 
		 & \left(\bftau^m,\bfdelta^m\right) \in \calQ\left(T,\calX_S^m,\calX_O^m,\calX_A^m,\pfrak^m \right) & 
		 \end{alignat*}
	\end{subequations}
	Theorem~\ref{theo:pomdp:NLP_optimal_solution} ensures that the optimization subproblem above on component $m$ corresponds to a POMDP problem with memoryless policies of POMDP $\left(\calX_S^m,\calX_O^m,\calX_A^m,\pfrak^m,\tilde{\bfr}\right)$ where $\tilde{r}^m(s,a,s') = r^m(s,a,s')- (\bfbeta^t)^{\bfT}\bfD^m(a)$ for any $s,s' \in \calX_S^m$ and $a \in \calX_A^m.$ 
	Thanks to Proposition~\ref{prop:pomdp:det_policies}, the subproblem on component $m$ can be solved using deterministic policies. Therefore, we can replace $\calQ\left(T,\calX_S^m,\calX_O^m,\calX_A^m,\pfrak^m \right)$ by $\calQ^{\mathrm{d}}\left(T,\calX_S^m,\calX_O^m,\calX_A^m,\pfrak^m \right)$ for any component $m$ and we obtain
	\begin{align*}
		\calG\left( \bfbeta \right) = \sum_{t=1}^T (\beta^t)^{\bfT} \bfb + \sum_{m=1}^M \calG^m(\bfbeta),
	\end{align*}
	which achieves the proof.
\qed


By definition of the Lagrangian relaxation, we have $z_{\rmLR}=\min_{\bfbeta \in \bbR_{+}^{Tq}} \calG(\bfbeta)$. We are now able to prove Proposition~\ref{prop:wkpomdp:lagrangian_relax}, which we recall here:

\lagrangian*

\proof
Thanks to Proposition~\ref{prop:wkpomdp:dual_function}, the dual functions of MILP~\eqref{pb:wkpomdp:decPOMDP_MILP} and NLP~\eqref{pb:wkpomdp:decPOMDP_NLP} are equal. It follows that the value of the Lagrangian relaxations are equal.

Now we prove inequality~\eqref{eq:wkpomdp:weak_duality}.
First, the inequality $z_{\rmUB} \leq z_{\rmLR}$ comes from weak duality (see e.g.~\citet[Proposition 5.1.3]{Bertsekas99}).
Second, to show the second inequality $z_{\rmLR} \leq z_{\rmR},$ it suffices to observe that the dual function $\calG(\bfbeta)$ of Problem~\eqref{pb:wkpomdp:decPOMDP_NLP} is also the dual function of Problem~\eqref{pb:wkpomdp:decPOMDP_MILP}. Indeed, in the expression of $\calG^m(\bfbeta)$ we can replace $\calQ\left(T,\calX_S^m,\calX_O^m,\calX_A^m,\pfrak^m \right)$ by $\calQ^{\mathrm{d}}\left(T,\calX_S^m,\calX_O^m,\calX_A^m,\pfrak^m \right)$  because the re always exists an optimal policy that is deterministic on the POMDP $\left(\calX_S^m,\calX_O^m,\calX_A^m,\pfrak^m,\bfr^m-\bfbeta^T\bfD^m \right)$.
A classical result in operations research (see e.g. \citet[Theorem 1]{Geoffrion1974}) states that the bound of the Lagrangian relaxation of an integer program is not worse than the bound of the linear relaxation. It shows that $z_{\rmLR} \leq z_{\rmR}.$

It remains to prove that $z_{\rmLR} \leq z_{\rmRc}$ and $z_{\rmRc} \leq z_{\rmR}.$ The second one comes from the fact that we have a smaller feasible set in the linear relaxation by adding the valid inequalities.
The first one comes by adding valid inequalities~\eqref{eq:wkpomdp:dec_Valid_cuts} in the expression of $\calG^m(\bfbeta)$, which is possible since the inequalities are valid, and by using the same arguments (weak duality and Geoffrion's Theorem) we conclude that $z_{\rmLR} \leq z_{\rmRc}.$   
\qed



\subsection{Column generation approach to compute $z_{\rmLR}$}
\label{sub:app_ColGen:col_gen}

In this section, we explain how we apply a column generation algorithm combined with a Dantzig-Wolfe decomposition to compute $z_{\rmLR}$. Proposition~\ref{prop:wkpomdp:lagrangian_relax} ensures that $z_{\rmLR}$ is also the value of the Lagrangian relaxation of MILP~\eqref{pb:wkpomdp:decPOMDP_MILP}.
Thanks to Geoffrion's Theorem \citep[Theorem 8.2]{Conforti:2014:IP:2765770}, the value of the Lagrangian relaxation MILP~\eqref{pb:wkpomdp:decPOMDP_MILP} satisfies:
\begin{equation}
\begin{aligned}\label{pb:app_ColGen:decPOMDP_MILP_DW}
z_{\rmLR} = \max_{(\bftau^m,\bfdelta^m)_{m \in [M]}} \enskip & \sum_{m=1}^M \sum_{t=1}^T \sum_{\substack{s,s' \in \calX_S^m \\ a \in \calX_A^m}} \tau_{sas'}^{t,m} r^m(s,a,s')  &\quad &\\
\mathrm{s.t.} \enskip 
  &(\bftau^m,\bfdelta^m) \in \rmConv\left( \calQ^{\rmd}(T,\calX_S^m,\calX_O^m,\calX_A^m, \pfrak^m) \right), & \forall m \in [M] \\
  & \sum_{m=1}^M \sum_{a \in \calX_A^m} \bfD^m(a)\tau_{a}^{t,m} \leq \bfb, & \forall t \in [T], 
 \end{aligned}
\end{equation}
where we included the constraints $\sum_{s \in \calX_S^m, o \in \calX_A^m}\tau_{soa}^{t,m} = \tau_a^{t,m}$ in the set $\calQ^{\rmd}(T,\calX_S^m,\calX_O^m,\calX_A^m, \pfrak^m)$ and $\rmConv(X)$ denotes the convex hull of a set $X$.

Using the definition of the convex hull, we can reformulate Problem~\eqref{pb:app_ColGen:decPOMDP_MILP_DW} as the following master problem:
\begin{subequations}\label{pb:app_ColGen:master}
 \begin{alignat}{2}
z_{\rmLR}=\max_{(\lambda^m,\bftau^m,\bfdelta^m)_{m\in [M]}} \enskip & \sum_{m=1}^M \sum_{\bftau \in \calM^m} \lambda_{\bftau}^m \bbE_{\bftau} \left[ \sum_{t=1}^T r^m(S_t^m,A_t^m,S_{t+1}^m) \right] &\quad &\\
\mathrm{s.t.} \enskip 
  &(\bftau^m,\bfdelta^m)= \sum_{(\bftau,\bfdelta)\in \calQ^{\rmd,m}} \lambda_{\bftau,\bfdelta}^m  (\bftau,\bfdelta) & \quad \forall m \in [M] \\
  &\sum_{(\bftau,\bfdelta) \in \calQ'^m} \lambda_{\bftau,\bfdelta}^m = 1  & \quad\forall m \in [M] \label{eq:app_ColGen:master_select_lambda} \\
  & \sum_{m=1}^M \bbE_{\bftau^m} \left[\bfD^m(A_t^m) \right] \leq \bfb  &  \forall t \in [T] \label{eq:app_ColGen:master_linking_contraints} \\
  & \lambda_{\bftau,\bfdelta}^m \geq 0  & \forall (\bftau,\bfdelta) \in \calQ^{\rmd,m}, \forall m \in [M]
 \end{alignat}
\end{subequations}
where $\calQ^{\rmd,m} =  \calQ^{\rmd}(T,\calX_S^m,\calX_O^m,\calX_A^m, \pfrak^m)$ for every $m$ in $[M]$.
It follows that the pricing subproblem on component $m$ writes down:
\begin{equation}\label{pb:app_ColGen:Subproblem}
 \begin{aligned}
z^m:=\max_{(\bftau,\bfdelta)} \enskip & \sum_{t=1}^T\sum_{\substack{s,s' \in\calX_S^m \\ a \in \calX_A^m}} \tau_{sas'}^{t}\left(r^m(s,a,s') - \beta_t^{\bfT}\bfD^m(a) \right) &\quad &\\
\mathrm{s.t.} \enskip 
  & (\bftau,\bfdelta) \in \calQ^{\rmd}(T,\calX_S^m,\calX_O^m,\calX_A^m,\pfrak^m), &  \forall t \in [T],
 \end{aligned}
\end{equation}
where $\beta_t \in \bbR_{+}^q$ is the vector dual variables of linking constraint~\eqref{eq:app_ColGen:master_linking_contraints}.  
It follows that the reduced cost can be written $c=\sum_{m=1}^M z^m + \pi^m$, where $(\pi^m)_{m\in[M]}$ is the vector of dual variables of constraints~\eqref{eq:app_ColGen:master_select_lambda}.

Now we are able to derive the column generation algorithm. We assume that $\calX_A \neq \emptyset$ (otherwise the decision maker cannot choose any action). Hence, there exists at least one element $\bfa \in \calX_A^1 \times \cdots \times \calX_A^M$ such that $\sum_{m=1}^M \bfD^m(a^m) \leq \bfb$.
Let $\bfa_{\rme}$ be such an element in $\calX_A$.

\begin{algorithm}[H]
\caption{Column generation to compute $z_{\rmLR}$.}
\label{alg:app_ColGen:ColGen}
\begin{algorithmic}[1]
\STATE \textbf{Input}: $T$, weakly coupled POMDP $\big((\calX_S^m,\calX_O^m,\calX_A^m,\pfrak^m,\bfr^m,\bfD^m)_{m\in[M]},\bfb\big)$
\STATE \textbf{Output}: The optimal value $z_{\rmLR}$ and 
\STATE \textbf{Initialize} $u_{\rmLB}\leftarrow - \infty$, $\calQ'^m \leftarrow \emptyset$
\FOR{$m=1,\ldots,M$}
\STATE Define $\bfdelta^m$ such that $\bfdelta^{t,m}_{a|o}= \mathds{1}_{a_{\rme}^m}(a)$ for every $a \in \calX_A^m$, $o \in \calX_O^m$ and $t \in [T]$.
\STATE Compute $\bftau^m$ such that $(\bftau^m,\bfdelta^m) \in \calQ\left(T,\calX_S^m,\calX_O^m,\calX_A^m,\pfrak^m,\bfr^m \right)$.
\STATE $z^m \leftarrow \infty$ and $\pi^m\leftarrow \infty$
\ENDFOR
\WHILE{$\sum_{m=1}^M z^m + \pi^m > 0$}
\STATE Add column: $\calQ'^m \leftarrow \calQ'^m \cup \{ (\bftau^m,\bfdelta^m) \}$
\STATE Solve Problem~\eqref{pb:app_ColGen:master} restricted to $(\calM'^m)_{m\in [M]}$ to obtain dual variables $(\bfbeta,\bfpi)$.
\STATE $u_{\rmLB} \leftarrow \text{Optimal value of the restricted master problem}$
\FOR{$m=1,\ldots,M$}
\STATE Set reward $\tilde{r}_t^m(s,a,s') :=r^m(s,a,s') - \beta_t^{\bfT}\bfD^m(a)$ for every $s,s' \in \calX_S^m$ and $a \in \calX_A^m$
\STATE Solve MILP~\eqref{pb:pomdp:MILP_pomdp} with valid inequalities~\eqref{eq:pomdp:Valid_cuts_pomdp} for POMDP $(\calX_S^m,\calX_O^m,\calX_A^m,\pfrak,\tilde{\bfr})$ to obtain $(\bftau^m,\bfdelta^m)$ and $z^m$
\ENDFOR
\ENDWHILE
\STATE Set $\bftau^m := \sum_{(\bftau,\bfdelta) \in \calQ'^m} \lambda_{\bftau'}^m \bftau'^m$ for every $m \in [M]$ 
\STATE \textbf{return} $(\bftau^m)_{m\in [M]}$ and $u_{\rmLB}$
\end{algorithmic}
\end{algorithm}

Algorithm~\ref{alg:app_ColGen:ColGen} computes the value of the Lagrangian relaxation $z_{\rmLR}$. We omit the proof in this paper.

\section{Numerical experiments}
\label{app:nums}

In this appendix, we provide further details on the instances of Section~\ref{sec:num} and additional numerical experiments.

\subsection{Generic POMDP: Random instances}
\label{sub:app_nums:random_instances}



All the instances are generated by first choosing $\vert\calX_S\vert$, $\vert\calX_O\vert$ and $\vert\calX_A\vert.$
We then randomly generate the initial probability $\big(p(s)\big)_{s \in \calX_S}$, the transition probability $\big(p(s'|s,a)\big)_{\substack{s,s' \in \calX_S\\a \in \calX_A}}$, the emission probability $\big(p(o|s)\big)_{\substack{s \in \calX_S \\o \in \calX_O}}$ and the immediate reward function $\big(r(s,a,s')\big)_{\substack{s,s' \in \calX_S\\a \in \calX_A}}$.
An instance is a tuple $(\calX_S,\calX_O,\calX_A,\pfrak,\bfr).$
A way to measure the difficulty of solving a MILP~\eqref{pb:pomdp:MILP_pomdp} for POMDP $(\calX_S,\calX_O,\calX_A,\pfrak,\bfr)$ with horizon $T$ can be characterized by the size of the set of deterministic policies $\vert \Delta_{\mathrm{ml}}^{\mathrm{d}} \vert = \vert\calX_A\vert^{T\vert\calX_O\vert}$, which only depends on the size of the observation space $\calX_O$ and the action space $\calX_A$. 
Since $\vert \Delta_{\mathrm{ml}}^{\mathrm{d}} \vert$ only depends on $\vert\calX_O\vert$ and $\vert\calX_A\vert$, we generate instances for different values of the pair $(\mathrm{k}_s,\mathrm{k}_a),$ where $\mathrm{k}_a = \vert\calX_O\vert = \vert\calX_A\vert$ and $\mathrm{k}_s=\vert\calX_S\vert.$

\subsection{Generic POMDP: Benchmark instances}
\label{sub:app_nums:benchmark}

All the instances can be found at the link~\url{http://pomdp.org/examples/} and further descriptions of each instance are available in the indicated literature on the same website.
In particular, it contains two instances \texttt{bridge-repair} and \texttt{machine} that model maintenance problems.
The first one, introduced by~\citet{Ellis1995}, consists of the maintenance of a bridge.
The modeling is almost similar to the one described in the introduction except that there are more available actions and they consider only one machine. Instead of just choosing whether or not to maintain the bridge, the decision maker chooses whether or not to inspect the bridge and, if so, whether or not to maintain it.
The second one, introduced by~\citet[Appendix H.3]{Cassandra1998}, consists of planning the maintenance of a machine with $4$ deteriorating components. Again, the decision maker can choose to inspect before performing a maintenance of the machine. In addition, the action ``maintenance'' is distinguished in two different actions: repair, which consists in maintaining internal components, and replace, which consists in replacing the machine by a new one. 
It leads to the set of available actions $ \calX_A = \{\mathrm{operate}, \mathrm{inspect}, \mathrm{repair}, \mathrm{replace}\}.$

\paragraph{Metrics.}
We give two metrics to evaluate MILP~\eqref{pb:pomdp:MILP_pomdp} against the SARSOP policy.
We want to compare the optimal value $z^*$ of MILP~\eqref{pb:pomdp:MILP_pomdp} with the value $z_{\mathrm{SARSOP}}$ obtained by using the SARSOP policy.
In addition, Theorem~\ref{theo:pomdp:MDP_approx_equivalence} says that $z^*$ and $z_{\mathrm{SARSOP}}$ are lower bounds of $v_{\mathrm{his}}^*.$
We also compare these values with $z_{\rmR^\mathrm{c}}^*$, the optimal value of the linear relaxation of MILP~\eqref{pb:pomdp:MILP_pomdp} with valid inequalities~\eqref{eq:pomdp:Valid_cuts_pomdp}. 
By Theorem~\ref{theo:pomdp:MDP_approx_equivalence}, the value of $z_{\rmR^\mathrm{c}}^*$ is an upper bound of $z^*$ and $v_{\mathrm{his}}^*$, and consequently an upper bound of $z_{\mathrm{SARSOP}}.$ 
It leads to the relative gap $g(z) = \frac{z_{\rmR^\mathrm{c}}^* - z}{z_{\rmR^\mathrm{c}}^*}$ for any $z$ belonging to $\{z^*, z_{\mathrm{SARSOP}}\}.$

\subsection{Weakly coupled POMDP: Multi-armed bandits with partial observations}
\label{sub:app_nums:bandits}

In this section, we provide numerical experiments on the partially observable multi-armed bandits problems that are introduced in Appendix~\ref{app:examples}.
We show the quality of the approximation~\eqref{pb:wkpomdp:decPOMDP_MILP} by comparing the values of $z_{\rm{LB}}$,$z_{\rmIP},$ $z_{\rmUB},$ $z_{\rmLR},$ $z_{\rmRc}$ and $z_{\rmR}.$


\paragraph{Instances.}
We consider instances where the state space and observation space of each bandit have the same cardinality $n$, i.e., $n:=|\calX_S^m| = |\calX_O^m|$ for any $m$ in $[M]$. 
The resulting system's state space and system's observation space have the size $|\calX_S| = |\calX_O| = n^M$. In each bandit state space, the states and observations are numbered from $1$ to $n$, i.e., $\calX_S^m=\calX_O^m = \{1,\ldots,n\}.$ 
Like \citet{bertsimas2016decomposable}, we generate different instances regular (REG.SAR), restless (RSTLS.SAR and RSTLS.SBR) or deterministic (RSTLS.DET.SBR) multi-armed bandits. For each set of instance, the emission probability vector $(p^m(o|s))_{o \in \calX_O^m,s \in \calX_S^m}$ is uniformly drawn from $[0,1]$ and renormalized. 
These sets of instances differ from each other in the structures of transition probabilities and reward functions.
We generate small-scale instances with $M \in \{2,3\}$ arms and $n=4$ states, and medium-scale instances with $M=5$ arms and $n=4$ states. 
The sets of instances are generated as follows. 
$T \in \{5,10, 20\}$.
\begin{itemize}
	\item REG.SAR consists of regular partially observable multi-armed bandits. 
	The reward function is defined by $r^m(s,1,s'):= (10/n) \cdot s$ and $r^m(s,0,s'):= 0$ for every state $s,s' \in \calX_S^m$ and every arm $m \in [M]$. 
	Each active transition probability vector $(p^m(s'|s,1))_{s,s' \in \calX_S^m}$ is drawn uniformly from $[0,1]$ and renormalized, for every arm $m \in [M]$. Each passive arm $m$ stays in the same state, i.e.,  $p^m(s'|s,0) = \mathds{1}_{s}(s')$ for every $s,s' \in \calX_S^m$.

	\item RSTLS.SAR consists of restless partially observable multi-armed bandits. 
	The reward function is the same as REG.SAR.
	Each active and passive transition probability vector $(p^m(s'|s,a))_{\substack{s,s' \in \calX_S^m \\ a \in \{0,1\}}}$ is drawn uniformly from $[0,1]$ and renormalized, for every arm $m \in [M]$. 

	\item RSTLS.SBR consists of restless partially observable multi-armed bandits. 
	The reward function is defined by $r^m(s,1,s'):= (10/n)\cdot s$ and $r^m(s,0,s'):= (1/M)\cdot(10/n)\cdot s$ for every state $s,s' \in \calX_S^m$ and every arm $m \in [M]$. 
	The transition probability is randomly drawn as RSTLS.SAR.   

	\item RSTLS.DET.SBR consists of restless partially observable multi-armed bandits. 
	The reward function is the same as RSTLS.SBR. 
	Each active and passive transition probability vector $(p^m(s'|s,a))_{\substack{s,s' \in \calX_S^m \\ a \in \{0,1\}}}$ is randomly drawn and deterministic, for every arm $m \in [M]$. 
\end{itemize}

\paragraph{Metrics.}
For each instance, we compute the value $z_{\rmIP},$ the lower bound $z_{\rmLB}$ and the upper bounds $z_{\rmUB},$ $z_{\rmLR},$ $z_{\rmRc}$ and $z_{\rmR}.$
Given an instance, we define the relative gaps with the largest upper bound $z_{\rmR}$: $\rmg_{\rmLB}=\frac{z_{\rmR} - z^{\rmLB}}{z_{\rmRc}},$  $\rmg_{\rmIP}=\frac{z_{\rmR} - z_{\rmIP}}{z_{\rmRc}},$ $\rmg_{\rmUB}=\frac{z_{\rmRc} - z_{\rmUB}}{z_{\rmRc}}$ and $\rmg_{\rmLR}=\frac{z_{\rmRc} - z_{\rmLR}}{z_{\rmRc}}$. 
Then, we define respectively the metrics $G_{\rm{mean}}(\rmg),$ $G_{\rm{95}}(\rmg)$ and $G_{\rm{max}}(\rmg)$ as the mean, the 95-th percentile and the maximum over a set of instances, for each gap $g$ in $\left\{\rmg_{\rmLB}, \rmg_{\rmIP}, \rmg_{\rmUB}, \rmg_{\rmLR} \right\}.$
In general, the lower the values of the metrics, the closer the bound is to the upper bound $z_{\rmRc}$.
In particular, thanks to Theorem~\ref{theo:wkpomdp:upper_bound_NLP} and Proposition~\ref{prop:wkpomdp:lagrangian_relax}the metrics $\rmg_{\rmLB}$ and $\rmg_{\rmLR}$ tell how close are the values of $z_{\rmIP}$ and $v_{\rmml}^*.$
Since the computation of $z_{\rmUB}$ becomes quickly difficult when the sizes of the instance increase, we only compute the values of $\rmg_{\rmUB}$ on small instances.  

Tables~\ref{tab:app_nums:bandits_numerical_results} summarizes the results.
For all the mathematical programs, we set the computation time limit to $3600$ seconds. If the resolution has not terminated before this time limit, then we keep the best upper bound obtained during the resolution. We do not have any guarantee about this upper bound, but it is the best one found by the solver during the resolution.
It explains why for some instances we obtain a smaller gap with the Lagrangian relaxation than with MILP~\eqref{pb:wkpomdp:decPOMDP_MILP}. 

\begin{table}
\centering
\resizebox{!}{7.6cm}
{\begin{tabular}{|c|cc|cccc|cccc|cccc|}
    	\hline
        Instance set & $T$ & $\rmg$ & \multicolumn{4}{c|}{$M=2$}  &  \multicolumn{4}{c|}{$M=3$}  &  \multicolumn{4}{c|}{$M=5$}  \\
        			 &     &     &  $G_{\rm{mean}}(\rmg)$ & $G_{\rm{95}}(\rmg)$ & $G_{\rm{max}}(\rmg)$ & Time(s) &$G_{\rm{mean}}(\rmg)$ & $G_{\rm{95}}(\rmg)$ & $G_{\rm{max}}(\rmg)$ & Time(s) &  $G_{\rm{mean}}(\rmg)$ & $G_{\rm{95}}(\rmg)$ & $G_{\rm{max}}(\rmg)$ & Time(s)\\
        \hline
        REG.SAR &  2  & $\rmg_{\rmLB}$  & 9.91 & 15.13 & 15.62 & 0.03 & 15.72 & 21.51 & 22.87 & 0.08 & 17.42 & 25.08 & 26.53 & 0.21 \\
			   	&	& $\rmg_{\rmIP}$  & 9.91 & 15.13 & 15.62 & 0.14 & 15.72 & 21.51 & 22.87 & 0.64 & 17.42 & 25.08 & 26.53 & 2.57 \\
			   	&	& $\rmg_{\rmUB}$ & 7.54 & 10.54 & 10.70 & 0.16 & 11.18 & 17.89 & 19.40 & 0.36 & $-$ & $-$ & $-$ & $-$ \\
			   	&	&  $\rmg_{\rmLR}$  & 7.02 & 10.21 & 10.33 & 9.55  & 10.86 & 16.97 & 18.73 & 14.08  & 13.00 & 18.41 & 18.64 & 15.52 \\ \cline{2-15}
				& 5  & $\rmg_{\rmLB}$  & 6.06 & 9.23 & 9.57 & 0.39 & 10.34 & 13.01 & 13.26 & 1.32 & 14.51 & 18.45 & 18.48 & 3.43\\
				&	 & $\rmg_{\rmIP}$  & 6.06 & 9.23 & 9.57 & 17.14 & 10.34 & 13.01 & 13.26 & 1525.63 & 17.10 & 27.03 & 27.04 & 2907.11 \\
				&	 &  $\rmg_{\rmUB}$  & 4.59 & 6.07 & 6.19 & 1260.58 & 7.18 & 9.89 & 9.99 & 3247.98 & $-$ & $-$ & $-$ &  $-$ \\
				&	 &  $\rmg_{\rmLR}$ & 4.04 & 5.89 & 5.95 & 43.00  & 6.79 & 9.36 & 9.42 & 54.66  & 10.46 & 13.15 & 13.61 & 51.17 \\ \cline{2-15}
				&	10  & $\rmg_{\rmLB}$  & 3.36 & 5.45 & 5.85 & 2.86 & 6.27 & 8.91 & 9.15 & 39.76 & 8.04 & 11.14 & 11.55 & 349.00 \\
				&	 & $\rmg_{\rmIP}$ & 3.38 & 5.55 & 6.04 & 1283.79 & 10.08 & 18.86 & 19.93 & 3205.30 & 15.53 & 23.63 & 24.08 & $>3600$ \\
				&	 &  $\rmg_{\rmUB}$  & 5.08 & 9.41 & 9.51 & 3248.01 & $-$ & $-$ & $-$ &  $-$ & $-$ & $-$ & $-$ &  $-$ \\
				&	 &  $\rmg_{\rmLR}$  & 2.13 & 3.27 & 3.46 & 946.01  & 3.90 & 5.74 & 5.84 & 1536.79  & 5.98 & 7.52 & 7.67 & 661.04 \\ 
		\hline
		RSTLS.SAR & 2  & $\rmg_{\rmLB}$  & 12.73 & 16.84 & 17.04 & 0.03 & 17.96 & 22.67 & 22.77 & 0.07 & 15.13 & 16.83 & 17.03 & 0.17 \\
				  &	 & $\rmg_{\rmIP}$  & 12.73 & 16.84 & 17.04 & 0.14 & 17.96 & 22.67 & 22.77 & 0.62 & 15.13 & 16.83 & 17.03 & 2.28 \\
				&	 &  $\rmg_{\rmUB}$  & 8.33 & 13.09 & 15.56 & 0.07 & 12.01 & 18.18 & 18.93 & 0.06 & $-$ & $-$ & $-$ &  $-$ \\
				&	 &  $\rmg_{\rmLR} $ & 8.16 & 12.71 & 14.99 & 9.22  & 11.99 & 18.15 & 18.91 & 14.44  & 9.94 & 11.20 & 11.28 & 15.79 \\ \cline{2-15}
				&	5  & $\rmg_{\rmLB}$  & 10.77 & 13.61 & 14.22 & 0.44 & 14.54 & 18.64 & 18.86 & 1.40 & 13.85 & 16.46 & 17.00 & 3.91 \\
				&	 & $\rmg_{\rmIP}$  & 10.77 & 13.61 & 14.22 & 12.50 & 14.54 & 18.64 & 18.86 & 358.36 & 18.31 & 21.51 & 21.66 & 3030.09 \\
				&	 &  $\rmg_{\rmUB} $ & 6.55 & 8.51 & 8.73 & 171.29 & 8.14 & 12.43 & 13.86 & 1716.52 & $-$ & $-$ & $-$ &  $-$ \\
				&	 &  $\rmg_{\rmLR}$  & 6.32 & 8.16 & 8.35 & 41.99  & 7.92 & 12.01 & 13.67 & 36.89  & 7.92 & 9.86 & 10.19 & 59.57 \\ \cline{2-15}
				&	10  & $\rmg_{\rmLB}$  & 10.39 & 13.73 & 14.32 & 3.20 & 13.18 & 16.36 & 16.83 & 29.58 & 14.55 & 17.85 & 18.18 & 312.76 \\
				&	 & $\rmg_{\rmIP}$ & 10.86 & 15.35 & 17.28 & 1896.65 & 14.42 & 17.97 & 18.07 & $>3600$ & 18.83 & 26.05 & 27.44 & $>3600$ \\
				&	 &  $\rmg_{\rmUB}$  & 5.99 & 8.44 & 9.28 & $>3600$ & $-$ & $-$ & $-$ &  $-$& $-$ & $-$ & $-$ &  $-$ \\
				&	 &  $\rmg_{\rmLR}$  & 5.44 & 7.29 & 7.83 & $>3600$  & 5.74 & 6.93 & 7.05 & 1669.80  & 7.34 & 9.42 & 9.68 & $>3600$ \\
		\hline
		RSTLS.SBR 	&  2 & 	$\rmg_{\rmLB}$  & 6.14 &  8.70 & 9.06 & 0.02 & 8.50 & 11.51 & 11.67 & 0.03 & 9.07 & 10.98 & 11.54 & 0.07 \\
		&  	 & 	$\rmg_{\rmIP}$  & 6.14 & 8.70 & 9.06 & 0.04 & 8.50 & 11.51 & 11.67 & 0.14 & 9.07 & 10.98 & 11.54 & 0.72 \\
		&  	 &  $\rmg_{\rmUB}$  & 3.54 & 5.63 & 6.00 & 0.02 & 5.39 & 7.60 & 7.75 & 0.02 & $-$ & $-$ & $-$ & $-$ \\
		&  	 &  $\rmg_{\rmLR}$  & 3.47 & 5.49 & 5.75 & 6.72  & 5.32 & 7.56 & 7.71 & 7.74  & 6.39 & 7.63 & 7.99 & 7.90 \\ \cline{2-15}
		&  	5  & $\rmg_{\rmLB}$  & 4.65 & 7.58 & 7.84 & 0.18 & 6.77 & 9.13 & 9.24 & 0.66 & 8.61 & 10.38 & 10.51 & 1.72 \\
		& 	 & $\rmg_{\rmIP}$  & 4.65 & 7.58 & 7.84 & 4.27 & 6.77 & 9.13 & 9.24 & 311.29 & 10.84 & 14.81 & 16.37 & 3291.89 \\
		&  	 &  $\rmg_{\rmUB}$  & 2.61 & 4.80 & 5.78 & 18.55 & 3.61 & 5.26 & 5.43 & 500.43 & $-$ & $-$ & $-$ & $-$ \\
		 & 	 & $\rmg_{\rmLR}$  & 2.33 & 4.16 & 4.85 & 14.43  & 3.43 & 5.14 & 5.27 & 21.78  & 5.10 & 5.89 & 5.91 & 19.59 \\ \cline{2-15}
		&  	10  & $\rmg_{\rmLB}$  & 4.27 & 7.82 & 8.52 & 1.29 & 6.00 & 8.86 & 9.20 & 16.71 & 7.99 & 10.75 & 11.28 & 127.66 \\
		&  	 & $\rmg_{\rmIP}$ & 4.37 & 7.90 & 8.52 & 890.65 & 9.87 & 15.37 & 15.70 & 3265.16 & 15.84 & 17.97 & 18.01 & $>3600$ \\
		&  	 &  $\rmg_{\rmUB}$  & 2.31 & 4.96 & 6.48 & 1836.06 & $-$ & $-$ & $-$ & $-$ & $-$ & $-$ & $-$ & $-$ \\
		&  	 &  $\rmg_{\rmLR}$  & 1.86 & 3.42 & 3.86 & 3168.63  & 2.77 & 4.69 & 4.76 & 3494.95  & 4.73 & 5.65 & 5.80 & 1443.62 \\
		\hline
		RSTLS.DET.SBR &2  & $\rmg_{\rmLB}$  & 6.08 & 11.58 & 13.78 & 0.02 & 8.52 & 13.94 & 14.43 & 0.03 & 9.74 & 13.54 & 13.57 & 0.11  \\
		&   & $\rmg_{\rmIP}$ & 6.08 & 11.58 & 13.78 & 0.05 & 8.52 & 13.94 & 14.43 & 0.26 & 9.74 & 13.54 & 13.57 & 1.76  \\
		&   &  $\rmg_{\rmUB}$  & 4.70 & 8.33 & 9.45 & 0.02 & 6.54 & 9.24 & 9.59 & 0.03 & $-$ & $-$ & $-$ & $-$ \\
		&   &  $\rmg_{\rmLR}$  & 4.67 & 8.31 & 9.42 & 12.44  & 6.35 & 9.22 & 9.56 & 14.00 & 7.43 & 10.17 & 10.60 & 13.76  \\ \cline{2-15}
		&  5  & $\rmg_{\rmLB}$  & 2.99 & 7.26 & 8.87 & 0.12 & 6.15 & 10.06 & 10.98 & 0.68 & 5.88 & 8.21 & 8.95 & 1.76  \\
		&   & $\rmg_{\rmIP}$ & 2.99 & 7.26 & 8.87 & 0.17 & 6.24 & 10.06 & 10.98 & 6.47 & 5.88 & 8.21 & 8.95 & 32.56  \\
		&   &  $\rmg_{\rmUB}$  & 2.72 & 6.86 & 8.66 & 2.97 & 5.01 & 8.11 & 9.52 & 77.43 & $-$ & $-$ & $-$ & $-$ \\
		&   &  $\rmg_{\rmLR}$  & 2.45 & 6.31 & 8.15 & 17.40  & 4.58 & 7.26 & 8.62 & 21.58  & 4.95 & 7.41 & 7.58 & 26.48  \\ \cline{2-15}
		&  10  & $\rmg_{\rmLB}$  & 2.39 & 6.94 & 7.50 & 0.59 & 4.29 & 8.36 & 10.01 & 3.90 & 4.70 & 7.54 & 8.25 & 43.95  \\
		&   & $\rmg_{\rmIP}$  & 2.39 & 6.94 & 7.50 & 4.79 & 4.29 & 8.36 & 10.01 & 514.32 & 6.38 & 13.10 & 13.87 & 1148.11  \\
		&   &  $\rmg_{\rmUB}$  & 2.23 & 7.04 & 7.60 & 816.32 & $-$ & $-$ & $-$ & $-$ & $-$ & $-$ & $-$ & $-$  \\
		&   &  $\rmg_{\rmLR}$  & 1.89 & 6.12 & 7.22 & 76.55  & 3.10 & 5.92 & 6.70 & 197.00  & 3.98 & 6.68 & 7.02 & 116.91  \\
		   \hline
  \end{tabular}}
\caption{The values of $G_{\rm{mean}}(\rmg),$ $G_{\rm{95}}(\rmg),$ and $G_{\rm{max}}(\rmg)$ obtained on the small-scale instances with $M \in \{2,3, 5\},$ $n=4$ and solved with different finite horizon $T \in \{2,5,10\}.$}
\label{tab:app_nums:bandits_numerical_results}
\end{table}

One can observe in Table~\ref{tab:app_nums:bandits_numerical_results} that for almost a large part of the instances, the values of $z_{\rm{LB}},$ $z_{\rmIP},$ $z_{\rmUB},$ and $z_{\rmLR}$ are close in general. 
It shows that our formulations have optimal values that are close to the optimal value $v_{\mathrm{ml}}^*$ of \ref{pb:decPOMDP_wc}.
In addition, the best bound obtained on the value of $z_{\rmIP}$ is very close to the value of the lower bound $z_{\rmLB}.$ Thanks to Theorem~\ref{theo:wkpomdp:lower_bound_MILP}, it means that most of the multi-armed bandit instances admit optimal policies that are ``decomposable'' (see Section~\ref{sub:wkpomdp:lb_ub}).

\subsection{Weakly coupled POMDP: Simulations of the implicit policy}
\label{sub:app_nums:implicit_policy}

In this section, detail how we build the instances of Section~\ref{sub:nums:implicit_policy}, and we provide the numerical results of matheuristic~\ref{alg:wkpomdp:matheuristic} involving the different policies in $\{\rmLB, \rmLR, \rmRc, \rmR \}$.

\paragraph{Instances} 
Like \citet{Walraven2018}, we build our instances of weakly coupled POMDP from the \texttt{bridge-repair} instance of \citet{Ellis1995} in which the decision maker has to perform maintenance on a bridge. 
In our problem, there are only two actions available on each bridge: either structural repair or keep. 
For each bridge $m$, the sizes of state space, observation space and action space are respectively $\vert \calX_S^m \vert=5$, $\vert \calX_O^m \vert=5$ and $\vert \calX_A^m \vert = 2$.
Each bridge starts almost surely in its most healthy state.
We add noises to the transition probabilities and emission probabilities of the bridges to ensure that they have slightly different parameters $\pfrak^m$ for all $m$ in $[M]$.
For every bridge $m$ in $[M]$, we set $C_F^m = 1000$ and $C_R^m=100$.
The bridges are inspected every months and evolve until the horizon of $H=24$ months. 
One instance consists in the value of the tuple $(M,K,(\pfrak^m)_{m\in [M]})$.
We build an instance as follows: first we choose a value of $M$ in $\{3,4,5,10,15, 20\}$, second we build the probabilities $\pfrak^m$ by adding a random real in $[0,0.1]$ to each non-zero value of the probabilities of \citet{Ellis1995}, and finally we choose $K = \max(\floor*{\gamma \times M},1)$, where $\gamma$ is a scalar belonging $\left\{ 0.2, 0.4, 0.6, 0.8\right\}$ (when $M=3$, then $K$ belongs to $\{1,2\}$). 
The range of values of $K$ is chosen in such a way that it goes from highly restrictive constraints (smallest values of $\gamma$) to more flexible constraints (largest values of $\gamma$), with respect to the value of $M$. When $K\geq M$, the decision maker can consider the subproblems separately, which is much easier. We enforce $K$ to be non smaller than $1$ because if $K=0$, the authority cannot maintain the bridges.

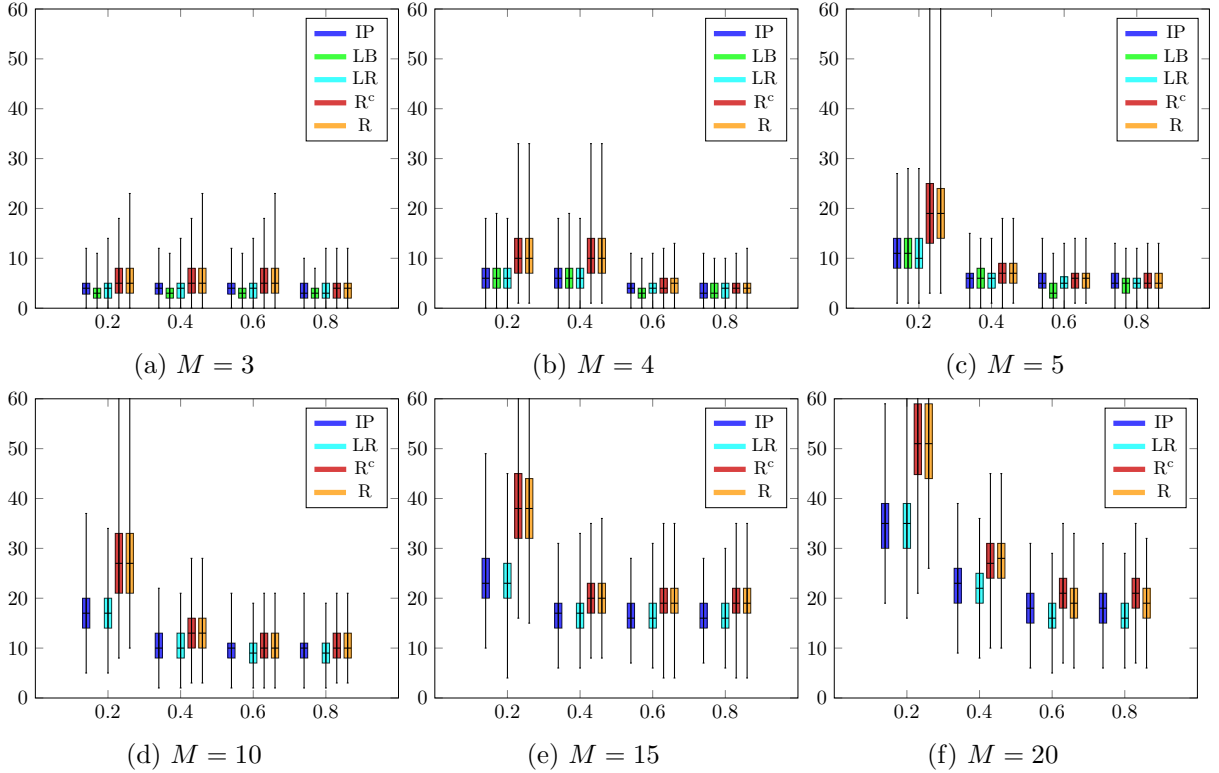
\begin{figure}
\begin{subfigure}[b]{0.33\textwidth}
  \centering
  \resizebox{!}{4.5cm}{
  \begin{tikzpicture}
    \begin{axis}[xmin=0,xmax=10,ymin=-.01, ymax=60,xtick={2,4,6,8}, ytick={0,10,20,30,40,50,60},xticklabels={$0.2$,$0.4$, $0.6$, $0.8$},legend pos=north east,legend style={font=\small}]

      \boxplot{1.4}{4.0}{5.0}{2.8}{0.0}{12.0}{blue!100!black}
	  \boxplot{1.7}{3.0}{4.0}{2.0}{0.0}{11.0}{green!100!black}
	  \boxplot{2.0}{4.0}{5.0}{2.0}{0.0}{14.0}{aqua}
	  \boxplot{2.3}{5.0}{8.0}{3.0}{0.0}{18.0}{red!80!black}
	  \boxplot{2.6}{5.0}{8.0}{3.0}{0.0}{23.0}{orange!120!black}
	  
	  \boxplot{3.4}{4.0}{5.0}{2.8}{0.0}{12.0}{blue!100!black}
	  \boxplot{3.7}{3.0}{4.0}{2.0}{0.0}{11.0}{green!100!black}
	  \boxplot{4.0}{4.0}{5.0}{2.0}{0.0}{14.0}{aqua}
	  \boxplot{4.3}{5.0}{8.0}{3.0}{0.0}{18.0}{red!80!black}
	  \boxplot{4.6}{5.0}{8.0}{3.0}{0.0}{23.0}{orange!120!black}
	  
	  \boxplot{5.4}{4.0}{5.0}{2.8}{0.0}{12.0}{blue!100!black}
	  \boxplot{5.7}{3.0}{4.0}{2.0}{0.0}{11.0}{green!100!black}
	  \boxplot{6.0}{4.0}{5.0}{2.0}{0.0}{14.0}{aqua}
	  \boxplot{6.3}{5.0}{8.0}{3.0}{0.0}{18.0}{red!80!black}
	  \boxplot{6.6}{5.0}{8.0}{3.0}{0.0}{23.0}{orange!120!black}
	  
	  \boxplot{7.4}{3.0}{5.0}{2.0}{0.0}{10.0}{blue!100!black}
	  \boxplot{7.7}{3.0}{4.0}{2.0}{0.0}{8.0}{green!100!black}
	  \boxplot{8.0}{3.0}{5.0}{2.0}{0.0}{12.0}{aqua}
	  \boxplot{8.3}{4.0}{5.0}{2.0}{0.0}{12.0}{red!80!black}
	  \boxplot{8.6}{4.0}{5.0}{2.0}{0.0}{12.0}{orange!120!black}


      \addlegendentry{$\rmIP$}          
      \addlegendimage{line width=1mm,color=blue!100!black,opacity=0.75}
      \addlegendentry{$\rmLB$}          
      \addlegendimage{line width=1mm,color=green!100!black,opacity=0.75}
      \addlegendentry{$\rmLR$}          
      \addlegendimage{line width=1mm,color=aqua, opacity=0.75}
      \addlegendentry{$\rmR^\rmc$}          
      \addlegendimage{line width=1mm,color=red!80!black,opacity=0.75}
      \addlegendentry{$\rmR$}          
      \addlegendimage{line width=1mm,color=orange!120!black,opacity=0.75}            
    \end{axis}
  \end{tikzpicture}}
  \subcaption{$M=3$}
  \label{fig:num:number_of_failures_3}
\end{subfigure}%
\begin{subfigure}[b]{0.33\textwidth}
  \centering
  \resizebox{!}{4.5cm}{
  \begin{tikzpicture}
    \begin{axis}[xmin=0,xmax=10,ymin=-.01, ymax=60,xtick={2,4,6,8}, ytick={0,10,20,30,40,50,60},xticklabels={$0.2$,$0.4$, $0.6$, $0.8$},legend pos=north east,legend style={font=\small}]

      \boxplot{1.4}{6.0}{8.0}{4.0}{0.0}{18.0}{blue!100!black}
	  \boxplot{1.7}{6.0}{8.0}{4.0}{0.0}{19.0}{green!100!black}
	  \boxplot{2.0}{6.0}{8.0}{4.0}{0.0}{18.0}{aqua}
	  \boxplot{2.3}{10.0}{14.0}{7.0}{1.0}{33.0}{red!80!black}
	  \boxplot{2.6}{10.0}{14.0}{7.0}{1.0}{33.0}{orange!120!black}
	  
	  \boxplot{3.4}{6.0}{8.0}{4.0}{0.0}{18.0}{blue!100!black}
	  \boxplot{3.7}{6.0}{8.0}{4.0}{0.0}{19.0}{green!100!black}
	  \boxplot{4.0}{6.0}{8.0}{4.0}{0.0}{18.0}{aqua}
	  \boxplot{4.3}{10.0}{14.0}{7.0}{1.0}{33.0}{red!80!black}
	  \boxplot{4.6}{10.0}{14.0}{7.0}{1.0}{33.0}{orange!120!black}
	  
	  \boxplot{5.4}{4.0}{5.0}{3.0}{0.0}{11.0}{blue!100!black}
	  \boxplot{5.7}{3.0}{4.0}{2.0}{0.0}{10.0}{green!100!black}
	  \boxplot{6.0}{4.0}{5.0}{3.0}{0.0}{11.0}{aqua}
	  \boxplot{6.3}{4.0}{6.0}{3.0}{0.0}{12.0}{red!80!black}
	  \boxplot{6.6}{5.0}{6.0}{3.0}{0.0}{13.0}{orange!120!black}
	  
	  \boxplot{7.4}{3.0}{5.0}{2.0}{0.0}{11.0}{blue!100!black}
	  \boxplot{7.7}{3.0}{5.0}{2.0}{0.0}{10.0}{green!100!black}
	  \boxplot{8.0}{4.0}{5.0}{2.0}{0.0}{10.0}{aqua}
	  \boxplot{8.3}{4.0}{5.0}{3.0}{0.0}{11.0}{red!80!black}
	  \boxplot{8.6}{4.0}{5.0}{3.0}{0.0}{12.0}{orange!120!black}

      \addlegendentry{$\rmIP$}          
      \addlegendimage{line width=1mm,color=blue!100!black, opacity=0.75}
      \addlegendentry{$\rmLB$}          
      \addlegendimage{line width=1mm,color=green!100!black, opacity=0.75}
      \addlegendentry{$\rmLR$}          
      \addlegendimage{line width=1mm,color=aqua, opacity=0.75}
      \addlegendentry{$\rmR^\rmc$}          
      \addlegendimage{line width=1mm,color=red!80!black, opacity=0.75}
      \addlegendentry{$\rmR$}          
      \addlegendimage{line width=1mm,color=orange!120!black, opacity=0.75}          
    \end{axis}
  \end{tikzpicture}}
  \subcaption{$M=4$}
  \label{fig:subfig2}
  \end{subfigure}
\begin{subfigure}[b]{0.33\textwidth}
  \centering
  \resizebox{!}{4.5cm}{
  \begin{tikzpicture}
    \begin{axis}[xmin=0,xmax=10,ymin=-.01, ymax=60,xtick={2,4,6,8}, ytick={0,10,20,30,40,50,60},xticklabels={$0.2$,$0.4$, $0.6$, $0.8$},legend pos=north east,legend style={font=\small}]

      \boxplot{1.4}{11.0}{14.0}{8.0}{1.0}{27.0}{blue!100!black}
	  \boxplot{1.7}{11.0}{14.0}{8.0}{1.0}{28.0}{green!100!black}
	  \boxplot{2.0}{10.0}{14.0}{8.0}{1.0}{28.0}{aqua}
	  \boxplot{2.3}{19.0}{25.0}{13.0}{3.0}{62.0}{red!80!black}
	  \boxplot{2.6}{19.0}{24.0}{14.0}{3.0}{62.0}{orange!120!black}
	  
	  \boxplot{3.4}{6.0}{7.0}{4.0}{0.0}{15.0}{blue!100!black}
	  \boxplot{3.7}{6.0}{8.0}{4.0}{0.0}{14.0}{green!100!black}
	  \boxplot{4.0}{6.0}{7.0}{4.0}{1.0}{14.0}{aqua}
	  \boxplot{4.3}{7.0}{9.0}{5.0}{0.0}{18.0}{red!80!black}
	  \boxplot{4.6}{7.0}{9.0}{5.0}{1.0}{18.0}{orange!120!black}
	  
	  \boxplot{5.4}{5.0}{7.0}{4.0}{0.0}{14.0}{blue!100!black}
	  \boxplot{5.7}{3.0}{5.0}{2.0}{0.0}{11.0}{green!100!black}
	  \boxplot{6.0}{5.0}{6.3}{4.0}{0.0}{13.0}{aqua}
	  \boxplot{6.3}{6.0}{7.0}{4.0}{1.0}{14.0}{red!80!black}
	  \boxplot{6.6}{6.0}{7.0}{4.0}{1.0}{14.0}{orange!120!black}
	  
	  \boxplot{7.4}{5.0}{7.0}{4.0}{0.0}{13.0}{blue!100!black}
	  \boxplot{7.7}{5.0}{6.0}{3.0}{0.0}{12.0}{green!100!black}
	  \boxplot{8.0}{5.0}{6.0}{4.0}{0.0}{12.0}{aqua}
	  \boxplot{8.3}{5.0}{7.0}{4.0}{0.0}{13.0}{red!80!black}
	  \boxplot{8.6}{5.0}{7.0}{4.0}{0.0}{13.0}{orange!120!black}

      \addlegendentry{$\rmIP$}          
      \addlegendimage{line width=1mm,color=blue!100!black, opacity=0.75}
      \addlegendentry{$\rmLB$}          
      \addlegendimage{line width=1mm,color=green!100!black, opacity=0.75}
      \addlegendentry{$\rmLR$}          
      \addlegendimage{line width=1mm,color=aqua, opacity=0.75}
      \addlegendentry{$\rmR^\rmc$}          
      \addlegendimage{line width=1mm,color=red!80!black, opacity=0.75}
      \addlegendentry{$\rmR$}          
      \addlegendimage{line width=1mm,color=orange!120!black, opacity=0.75}          
    \end{axis}
  \end{tikzpicture}}
  \subcaption{$M=5$}
  \label{fig:subfig2}
  \end{subfigure}
  \vspace{0.5cm}
  \begin{subfigure}[b]{0.33\textwidth}
  \centering
  \resizebox{!}{4.5cm}{
  \begin{tikzpicture}
    \begin{axis}[xmin=0,xmax=10,ymin=-.01, ymax=60,xtick={2,4,6,8}, ytick={0,10,20,30,40,50,60},xticklabels={$0.2$,$0.4$, $0.6$, $0.8$},legend pos=north east,legend style={font=\small}]

      \boxplot{1.4}{17.0}{20.0}{14.0}{5.0}{37.0}{blue!100!black}
	  \boxplot{2.0}{17.0}{20.0}{14.0}{5.0}{34.0}{aqua}
	  \boxplot{2.3}{27.0}{33.0}{21.0}{8.0}{71.0}{red!80!black}
	  \boxplot{2.6}{27.0}{33.0}{21.0}{10.0}{61.0}{orange!120!black}
	  
	  \boxplot{3.4}{10.0}{13.0}{8.0}{2.0}{22.0}{blue!100!black}
	  \boxplot{4.0}{10.0}{13.0}{8.0}{2.0}{21.0}{aqua}
	  \boxplot{4.3}{13.0}{16.0}{10.0}{3.0}{28.0}{red!80!black}
	  \boxplot{4.6}{13.0}{16.0}{10.0}{3.0}{28.0}{orange!120!black}
	  
	  \boxplot{5.4}{10.0}{11.0}{8.0}{2.0}{21.0}{blue!100!black}
	  \boxplot{6.0}{9.0}{11.0}{7.0}{2.0}{19.0}{aqua}
	  \boxplot{6.3}{10.0}{13.0}{8.0}{2.0}{21.0}{red!80!black}
	  \boxplot{6.6}{10.0}{13.0}{8.0}{2.0}{21.0}{orange!120!black}
	  
	  \boxplot{7.4}{10.0}{11.0}{8.0}{2.0}{21.0}{blue!100!black}
	  \boxplot{8.0}{9.0}{11.0}{7.0}{2.0}{19.0}{aqua}
	  \boxplot{8.3}{10.0}{13.0}{8.0}{3.0}{21.0}{red!80!black}
	  \boxplot{8.6}{10.0}{13.0}{8.0}{3.0}{21.0}{orange!120!black}

      \addlegendentry{$\rmIP$}          
      \addlegendimage{line width=1mm,color=blue!100!black,opacity=0.75}
      \addlegendentry{$\rmLR$}          
      \addlegendimage{line width=1mm,color=aqua,opacity=0.75}
      \addlegendentry{$\rmR^\rmc$}          
      \addlegendimage{line width=1mm,color=red!80!black,opacity=0.75}
      \addlegendentry{$\rmR$}          
      \addlegendimage{line width=1mm,color=orange!120!black,opacity=0.75}         
    \end{axis}
  \end{tikzpicture}}
  \subcaption{$M=10$}
  \label{fig:num:number_of_failures_10}
  \end{subfigure}%
  \begin{subfigure}[b]{0.33\textwidth}
  \centering
  \resizebox{!}{4.5cm}{
  \begin{tikzpicture}
    \begin{axis}[xmin=0,xmax=10,ymin=-.01, ymax=60,xtick={2,4,6,8}, ytick={0,10,20,30,40,50,60},xticklabels={$0.2$,$0.4$, $0.6$, $0.8$},legend pos=north east,legend style={font=\small}]
      \boxplot{1.4}{23.0}{28.0}{20.0}{10.0}{49.0}{blue!100!black}
	  \boxplot{2.0}{23.0}{27.0}{20.0}{4.0}{45.0}{aqua}
	  \boxplot{2.3}{38.0}{45.0}{32.0}{16.0}{88.0}{red!80!black}
	  \boxplot{2.6}{38.0}{44.0}{32.0}{15.0}{88.0}{orange!120!black}
	  
	  \boxplot{3.4}{17.0}{19.0}{14.0}{6.0}{31.0}{blue!100!black}
	  \boxplot{4.0}{17.0}{19.0}{14.0}{6.0}{33.0}{aqua}
	  \boxplot{4.3}{20.0}{23.0}{17.0}{8.0}{35.0}{red!80!black}
	  \boxplot{4.6}{20.0}{23.0}{17.0}{8.0}{36.0}{orange!120!black}
	  
	  \boxplot{5.4}{16.0}{19.0}{14.0}{7.0}{28.0}{blue!100!black}
	  \boxplot{6.0}{16.0}{19.0}{14.0}{6.0}{31.0}{aqua}
	  \boxplot{6.3}{19.0}{22.0}{17.0}{4.0}{35.0}{red!80!black}
	  \boxplot{6.6}{19.0}{22.0}{17.0}{4.0}{35.0}{orange!120!black}
	  
	  \boxplot{7.4}{16.0}{19.0}{14.0}{7.0}{28.0}{blue!100!black}
	  \boxplot{8.0}{16.0}{19.0}{14.0}{6.0}{30.0}{aqua}
	  \boxplot{8.3}{19.0}{22.0}{17.0}{4.0}{35.0}{red!80!black}
	  \boxplot{8.6}{19.0}{22.0}{17.0}{4.0}{35.0}{orange!120!black}


      \addlegendentry{$\rmIP$}          
      \addlegendimage{line width=1mm,color=blue!100!black,opacity=0.75}
      \addlegendentry{$\rmLR$}          
      \addlegendimage{line width=1mm,color=aqua,opacity=0.75}
      \addlegendentry{$\rmR^\rmc$}          
      \addlegendimage{line width=1mm,color=red!80!black,opacity=0.75}
      \addlegendentry{$\rmR$}          
      \addlegendimage{line width=1mm,color=orange!120!black,opacity=0.75}         
    \end{axis}
  \end{tikzpicture}}
  \subcaption{$M=15$}
  \label{fig:num:number_of_failures_15}
  \end{subfigure}%
  \begin{subfigure}[b]{0.33\textwidth}
  \centering
  \resizebox{!}{4.5cm}{
  \begin{tikzpicture}
    \begin{axis}[xmin=0,xmax=10,ymin=-.01, ymax=60,xtick={2,4,6,8}, ytick={0,10,20,30,40,50,60},xticklabels={$0.2$,$0.4$, $0.6$, $0.8$},legend pos=north east,legend style={font=\small}]
	  \boxplot{1.4}{35.0}{39.0}{30.0}{19.0}{59.0}{blue!100!black}
	  \boxplot{2.0}{35.0}{39.0}{30.0}{16.0}{65.0}{aqua}
	  \boxplot{2.3}{51.0}{59.0}{44.8}{21.0}{105.0}{red!80!black}
	  \boxplot{2.6}{51.0}{59.0}{44.0}{26.0}{108.0}{orange!120!black}
	  
	  \boxplot{3.4}{23.0}{26.0}{19.0}{9.0}{39.0}{blue!100!black}
	  \boxplot{4.0}{22.0}{25.0}{19.0}{8.0}{36.0}{aqua}
	  \boxplot{4.3}{27.0}{31.0}{24.0}{10.0}{45.0}{red!80!black}
	  \boxplot{4.6}{28.0}{31.0}{24.0}{10.0}{45.0}{orange!120!black}
	  
	  \boxplot{5.4}{18.0}{21.0}{15.0}{6.0}{31.0}{blue!100!black}
	  \boxplot{6.0}{16.0}{19.0}{14.0}{5.0}{29.0}{aqua}
	  \boxplot{6.3}{21.0}{24.0}{18.0}{7.0}{35.0}{red!80!black}
	  \boxplot{6.6}{19.0}{22.0}{16.0}{6.0}{33.0}{orange!120!black}
	  
	  \boxplot{7.4}{18.0}{21.0}{15.0}{6.0}{31.0}{blue!100!black}
	  \boxplot{8.0}{16.0}{19.0}{14.0}{6.0}{29.0}{aqua}
	  \boxplot{8.3}{21.0}{24.0}{18.0}{7.0}{35.0}{red!80!black}
	  \boxplot{8.6}{19.0}{22.0}{16.0}{6.0}{32.0}{orange!120!black}

      
      \addlegendentry{$\rmIP$}          
      \addlegendimage{line width=1mm,color=blue!100!black,opacity=0.75}
      \addlegendentry{$\rmLR$}          
      \addlegendimage{line width=1mm,color=aqua,opacity=0.75}
      \addlegendentry{$\rmR^\rmc$}          
      \addlegendimage{line width=1mm,color=red!80!black,opacity=0.75}
      \addlegendentry{$\rmR$}          
      \addlegendimage{line width=1mm,color=orange!120!black,opacity=0.75}
    \end{axis}

  \end{tikzpicture}}
  \subcaption{$M=20$}
  \label{fig:num:number_of_failures_20}
  \end{subfigure}%
\caption{Boxplots on the total number of failures counted during the simulations of the policies $\rmIP$ (blue boxplot), $\rmLB$ (green boxplot), $\rmLR$ (light blue) , $\rmR^\rmc$ (red boxplot) and $\rmR$ (orange boxplot) for different values of $M$ and $K$ and a rolling horizon $T_{\rmr}=5$.}
\label{fig:num:nb_failures}
\end{figure}

Tables~\ref{tab:wkpomdp:maintenance_1},~\ref{tab:wkpomdp:maintenance_2},~\ref{tab:wkpomdp:maintenance_3} summarize the results. In addition, Figure~\ref{fig:num:nb_failures} displays the statistics about the total number of failures counted during simulations of the policies.  
For all the mathematical programs, we set the computation time limit to $3600$ seconds and a final gap tolerance (\texttt{MIPGap} parameter in \texttt{Gurobi}) of $1 \%$, which is enough for the use of our matheuristic.
If the resolution has not terminated before this time limit, then we keep the best feasible solution at the end of the resolution.

\begin{table}
  \centering
  \resizebox{!}{7.2cm}{
  \begin{tabular}{|c|cc|ccccccc|ccccccc|}
      \hline
        \multirow{3}{*}{$\gamma$} & \multirow{3}{*}{$\rmp$} & \multirow{3}{*}{$T_\rmr$} &  \multicolumn{7}{c|}{$M=3$} & \multicolumn{7}{c|}{$M=4$}  \\
            & & & \multicolumn{2}{c}{$\vert \nu_{\rmp} \vert$ ($\times 10^3$)} & \multicolumn{2}{c}{$\rmF_{\rmp}$} & $\rmG_{\rmp}^{\rmLR}$ & $\rmG_{\rmp}^{\rmRc}$ & Time & \multicolumn{2}{c}{$\vert \nu_{\rmp} \vert$ ($\times 10^3$)} & \multicolumn{2}{c}{$\rmF_{\rmp}$} & $\rmG_{\rmp}^{\rmLR}$ & $\rmG_{\rmp}^{\rmRc}$ & Time \\
             & & & {Mean} & {Std. err.} & {Mean} & {Std. err.} & $(\%)$ & $(\%)$ & ($\si{s}$) & {Mean} & {Std. err.} & {Mean} & {Std. err.} & $(\%)$ & $(\%)$ & ($\si{s}$) \\
            \hline
         0.2 & $\rmIP$ & 2 & 5.71 & 2.32 & 4.6 & 2.2 & 11.60 & 14.09 & 0.004 & 8.81 & 2.87 & 6.5 & 2.9 & -1.77 & 13.79 & 0.007 \\ 
		     &      & 5 & 5.45 & 2.16 & 4.0 & 2.1 & 6.51 & 8.88 & 0.055 & \textbf{8.64} & \textbf{2.92} & \textbf{6.3} & \textbf{2.9} & \textbf{-3.64} & \textbf{11.62} & \textbf{0.116} \\ \cline{2-17}
		     & $\rmLB$ & 2 & 5.71 & 2.31 & 4.5 & 2.2 & 11.61 & 14.10 & 0.005 & 8.81 & 2.88 & 6.5 & 2.9 & -1.75 & 13.81 & 0.013 \\
		     &      & 5 & 5.47 & 1.89 & 3.2 & 1.9 & 6.79 & 9.18 & 0.096 & 8.86 & 2.98 & 6.6 & 3.0 & -1.18 & 14.47 & 0.165 \\ \cline{2-17}
		     & $\rmLR$ & 2 & 5.70 & 2.42 & 4.6 & 2.3 & 11.45 & 13.93 & 0.162 & 8.91 & 3.04 & 6.6 & 3.0 & -0.63 & 15.10 & 0.266 \\
		     &      & 5 & \textbf{5.42} & \textbf{2.28} & \textbf{4.0} & \textbf{2.2} & \textbf{5.95} & \textbf{8.31} & \textbf{0.615} & 8.69 & 3.07 & 6.4 & 3.1 & -3.09 & 12.26 & 2.699 \\ \cline{2-17}
		     & $\rmRc$ & 2 & 7.66 & 3.53 & 6.6 & 3.5 & 49.73 & 53.07 & 0.002 & 12.06 & 4.83 & 10.2 & 4.9 & 34.55 & 55.85 & 0.002 \\
		     &      & 5 & 6.95 & 3.14 & 5.7 & 3.1 & 35.85 & 38.88 & 0.006 & 12.51 & 5.33 & 10.7 & 5.4 & 39.60 & 61.71 & 0.010 \\ \cline{2-17}
		     & $\rmR$ & 2 & 7.67 & 3.52 & 6.6 & 3.5 & 49.85 & 53.19 & 0.004 & 11.87 & 4.62 & 10.0 & 4.7 & 32.44 & 53.41 & 0.006 \\ 
		     &      & 5 & 6.96 & 3.20 & 5.7 & 3.2 & 35.94 & 38.97 & 0.020 & 12.51 & 5.37 & 10.7 & 5.4 & 39.59 & 61.69 & 0.050 \\
		\hline
		0.4 & $\rmIP$  & 2 & 5.71 & 2.32 & 4.6 & 2.2 & 11.60 & 14.09 & 0.005  & 8.81 & 2.87 & 6.5 & 2.9 & -1.77 & 13.79 & 0.007 \\
		     &      & 5 & 5.45 & 2.16 & 4.0 & 2.1 & 6.51 & 8.88 & 0.062 & \textbf{8.64} & \textbf{2.92} & \textbf{6.3} & \textbf{2.9} & \textbf{-3.64} & \textbf{11.62} & \textbf{0.118} \\ \cline{2-17}
		     & $\rmLB$ & 2 & 5.71 & 2.31 & 4.5 & 2.2 & 11.61 & 14.10 & 0.006 & 8.81 & 2.88 & 6.5 & 2.9 & -1.75 & 13.81 & 0.013 \\
		     &      & 5 & 5.47 & 1.89 & 3.2 & 1.9 & 6.79 & 9.18 & 0.106 & 8.86 & 2.98 & 6.6 & 3.0 & -1.18 & 14.47 & 0.167 \\ \cline{2-17}
		     & $\rmLR$ & 2 & 5.70 & 2.42 & 4.6 & 2.3 & 11.45 & 13.93 & 0.163 & 8.91 & 3.04 & 6.6 & 3.0 & -0.63 & 15.10 & 0.267 \\
		     &      & 5 & \textbf{5.42} & \textbf{2.28} & \textbf{4.0} & \textbf{2.2} & \textbf{5.95} & \textbf{8.31} & \textbf{0.618} & 8.69 & 3.07 & 6.4 & 3.1 & -3.09 & 12.26 & 2.701 \\ \cline{2-17}
		     & $\rmRc$& 2 & 7.66 & 3.53 & 6.6 & 3.5 & 49.73 & 53.07 & 0.002 & 12.06 & 4.83 & 10.2 & 4.9 & 34.55 & 55.85 & 0.002 \\
		     &      & 5 & 6.95 & 3.14 & 5.7 & 3.1 & 35.85 & 38.88 & 0.006 & 12.51 & 5.33 & 10.7 & 5.4 & 39.60 & 61.71 & 0.010 \\ \cline{2-17}
		     & $\rmR$ & 2 & 7.67 & 3.52 & 6.6 & 3.5 & 49.85 & 53.19 & 0.004 & 11.87 & 4.62 & 10.0 & 4.7 & 32.44 & 53.41 & 0.006 \\
		     &      & 5 & 6.96 & 3.20 & 5.7 & 3.2 & 35.94 & 38.97 & 0.020 & 12.51 & 5.37 & 10.7 & 5.4 & 39.59 & 61.69 & 0.051 \\
		\hline
		0.6 & $\rmIP$ & 2 & 5.71 & 2.32 & 4.6 & 2.2 & 11.60 & 14.09 & 0.006 & 7.70 & 1.96& 4.3 & 1.9 & -0.94 & 0.46 & 0.006 \\
		     &      & 5 & 5.45 & 2.16 & 4.0 & 2.1 & 6.51 & 8.88 & 0.064 & \textbf{7.62} & \textbf{1.89} & \textbf{4.0} & \textbf{1.9} & \textbf{-2.05} & \textbf{-0.66} & \textbf{0.067} \\ \cline{2-17}
		     & $\rmLB$ & 2 & 5.71 & 2.31 & 4.5 & 2.2 & 11.61 & 14.10 & 0.006 & 7.71 & 1.96 & 4.3 & 1.9 & -0.85 & 0.55 & 0.010 \\
		     &      & 5 & 5.47 & 1.89 & 3.2 & 1.9 & 6.79 & 9.18 & 0.110 & 7.80 & 1.75 & 3.3 & 1.7 & 0.35 & 1.77 & 0.103 \\ \cline{2-17}
		     & $\rmLR$ & 2 & 5.70 & 2.42 & 4.6 & 2.3 & 11.45 & 13.93 & 0.164 & 7.88 & 2.11 & 4.5 & 2.0 & 1.37 & 2.80 & 0.161 \\
		     &      & 5 & \textbf{5.42} & \textbf{2.28} & \textbf{4.0} & \textbf{2.2} & \textbf{5.95} & \textbf{8.31} & \textbf{0.620} & 7.70 & 1.95 & 4.0 & 1.9 & -0.97 & 0.43 & 0.723 \\ \cline{2-17}
		     & $\rmRc$ & 2.0 & 7.66 & 3.53 & 6.6 & 3.5 & 49.73 & 53.07 & 0.002 & 8.48 & 2.35 & 5.1 & 2.3 & 9.01 & 10.55 & 0.002 \\
		     &      & 5 & 6.95 & 3.14 & 5.7 & 3.1 & 35.85 & 38.88 & 0.007 & 7.83 & 1.93 & 3.9 & 1.9 & 0.71 & 2.12 & 0.008 \\ \cline{2-17}
		     & $\rmR$ & 2 & 7.67 & 3.52 & 6.6 & 3.5& 49.85 & 53.19 & 0.005 & 8.29 & 2.24 & 4.5 & 2.2 & 6.63 & 8.13 & 0.002 \\
		     &      & 5 & 6.96 & 3.20 & 5.7 & 3.2 & 35.94 & 38.97 & 0.021 & 8.28 & 2.27& 4.8 & 2.3 & 6.46 & 7.97 & 0.046 \\
		    \hline
		0.8  & $\rmIP$ & 2 & 5.17 & 1.90 & 4.0 & 1.8 & 0.95 & 3.21 & 0.005 & 7.70 & 1.91 & 4.1 & 1.8 & -0.94 & 0.45 & 0.006 \\
		     &      & 5 & 5.11 & 1.86 & 3.5 &1.8 & -0.12 & 2.12 & 0.056 & \textbf{7.61} & \textbf{1.81} & \textbf{3.6} & \textbf{1.7} & \textbf{-2.16} & \textbf{-0.79} & \textbf{0.060} \\ \cline{2-17}
		     & $\rmLB$ & 2 & 5.17 & 1.89 & 3.9 & 1.8 & 1.01 & 3.27 & 0.005 & 7.70 & 1.90 & 4.1 & 1.8 & -0.97 & 0.42 & 0.009 \\
		     &      & 5 & 5.32 & 1.62 & 3.0 & 1.6 & 3.91 & 6.23 & 0.092 & 7.65 & 1.79 & 3.6 & 1.7 & -1.66 & -0.28 & 0.098 \\ \cline{2-17}
		     & $\rmLR$ & 2 & 5.18 & 2.03 & 4.0 & 1.9 & 1.23 & 3.50 & 0.124& 7.82 & 2.03 & 4.2 & 1.9 & 0.59 & 2.01 & 0.162 \\
		     &      & 5 & \textbf{5.10} & \textbf{1.91} & \textbf{3.5} & \textbf{1.8} & \textbf{-0.26} & \textbf{1.98} & \textbf{0.364} & 7.69 & 1.86 & 3.7 & 1.8 & -1.14 & 0.24 & 0.466 \\ \cline{2-17}
		     & $\rmRc$ & 2 & 6.34 & 2.58 & 5.1 & 2.5 & 23.79 & 26.56 & 0.002 & 8.29 & 2.24 & 4.5 & 2.2 & 6.63 & 8.13 & 0.002 \\
		     &      & 5 & 5.26 & 1.96 & 3.7 & 1.9 & 2.75 & 5.05 & 0.006 & 7.83 & 1.93 & 3.9 & 1.9 & 0.71 & 2.12 & 0.008 \\ \cline{2-17}
		     & $\rmR$ & 2 & 6.33 & 2.57 & 5.1 & 2.5 & 23.77 & 26.54 & 0.005 & 8.27 & 2.22& 4.5 & 2.2 & 6.35 & 7.84 & 0.006 \\
		     &      & 5 & 5.26 & 1.96 & 3.7 & 1.9 & 2.75 & 5.05 & 0.020 & 7.87 & 1.98 & 3.9 & 1.9 & 1.16 & 2.58 & 0.044 \\
		\hline
      \end{tabular}}
\caption{Numerical values of $\vert \nu_{\rmp} \vert$, $\rmF_{\rmp}$ (mean and standard errors), $\rmG_{\rmp}^{\rmLR}$ and $\rmG_{\rmp}^{\rmRc}$ obtained on an instance $(M,\gamma)$ with $M \in \{ 15, 20\}$ and $\gamma \in \{ 0.2, 0.4, 0.6, 0.8\}$. The values written in bold indicate the best performances regarding the objective values.}
\label{tab:wkpomdp:maintenance_1}
\end{table}

\begin{table}
  \centering
  \resizebox{!}{7.0cm}{
  \begin{tabular}{|c|cc|ccccccc|ccccccc|}
      \hline
        \multirow{3}{*}{$\gamma$} & \multirow{3}{*}{$\rmp$} & \multirow{3}{*}{$T_\rmr$} &  \multicolumn{7}{c|}{$M=5$} & \multicolumn{7}{c|}{$M=10$}  \\
            & & & \multicolumn{2}{c}{$\vert \nu_{\rmp} \vert$ ($\times 10^3$)} & \multicolumn{2}{c}{$\rmF_{\rmp}$} & $\rmG_{\rmp}^{\rmLR}$ & $\rmG_{\rmp}^{\rmRc}$ & Time & \multicolumn{2}{c}{$\vert \nu_{\rmp} \vert$ ($\times 10^3$)} & \multicolumn{2}{c}{$\rmF_{\rmp}$} & $\rmG_{\rmp}^{\rmLR}$ & $\rmG_{\rmp}^{\rmRc}$ & Time \\
             & & & {Mean} & {Std. err.} & {Mean} & {Std. err.} & $(\%)$ & $(\%)$ & ($\si{s}$) & {Mean} & {Std. err.} & {Mean} & {Std. err.} & $(\%)$ & $(\%)$ & ($\si{s}$) \\
        \hline
        0.2  & $\rmIP$ & 2 & 13.78 & 4.77 & 11.5 & 4.8 & -41.47 & 30.24 & 0.008 & 22.63 & (5.45) & 18.0 & (5.5) & -34.20 & 17.71 & 0.012 \\
		     &      & 5 & 13.38 &4.58 & 11.1 & 4.6 & -43.18 &26.44 & 0.250 & 22.06 & (5.24) & 17.5 & (5.2) & -35.85 & 14.74 & 0.384 \\ \cline{2-17}
		     & $\rmLB$ & 2 & 13.79 & 4.77 & 11.5 & 4.8 & -41.46 & 30.26 & 0.065 & $-$ & $-$ & $-$ & $-$ & $-$ & $-$ & $-$ \\
		     &      & 5 & 13.64 & 4.52 & 11.3 & 4.5 & -42.07 & 28.89 & 0.378 & $-$ & $-$ & $-$ & $-$ & $-$ & $-$ & $-$ \\ \cline{2-17}
		     & $\rmLR$ & 2 & 13.69 & 4.42 & 11.4 & 4.4 & -41.88 & 29.34 & 0.321 & 22.43 & (5.20) & 17.8 & (5.2) & -34.78 & 16.67 & 0.467 \\
		     &      & 5 & \textbf{13.27} & \textbf{4.41} & \textbf{11.0} & \textbf{4.4} & \textbf{-43.63} & \textbf{25.43} & \textbf{2.678} & \textbf{21.77} & \textbf{4.88} & \textbf{17.2} & \textbf{4.9} & \textbf{-36.70} & \textbf{13.23} & \textbf{5.839} \\ \cline{2-17}
		     & $\rmRc$ & 2 & 20.38 & 7.16 & 18.6 & 7.2 & -13.47 & 92.54 & 0.003 & 32.31 & 8.34 & 28.6 & 8.4 & -6.06 & 68.03 & 0.005 \\
		     &      & 5 & 21.09 & 7.83 & 19.3 & 7.9 & -10.46 & 99.24 & 0.015 & 31.74 & 9.13 & 28.0 & 9.2 & -7.72 & 65.07 & 0.027 \\ \cline{2-17}
		     & $\rmR$ & 2 & 20.05 & 6.85 & 18.2 & 6.9 & -14.88 & 89.41 & 0.007 & 31.92 & 8.03 & 28.2 & 8.1 & -7.18 & 66.04 & 0.015 \\
		     &      & 5 & 21.03 & 7.77 & 19.3 & 7.8 & -10.70 & 98.71 & 0.070 & 31.64 & 9.02 & 27.9 & 9.1 & -8.00 & 64.56 & 0.136 \\
		\hline
		0.4 & $\rmIP$ & 2 & 10.43 & 2.50 & 6.4 & 2.5 & 0.76 & 2.13 & 0.007 & 19.19 & 3.38 & 11.5 & 3.3 & 1.32 & 3.07 & 0.012 \\
		     &      & 5 & \textbf{10.27} & \textbf{2.39} & \textbf{5.7} & \textbf{2.4} & \textbf{-0.77} & \textbf{0.58} & \textbf{0.124} & 18.91 & 3.26 & 10.6 & 3.2 & -0.16 & 1.57 & 0.196 \\ \cline{2-17}
		     & $\rmLB$ & 2 & 10.44 & 2.49 & 6.4 & 2.5 & 0.90 & 2.27 & 0.030 & $-$ & $-$ & $-$ & $-$ & $-$ & $-$ & $-$ \\
		     &      & 5 & 10.51 & 2.38 & 6.0 & 2.4 & 1.56 & 2.94 & 0.776 & $-$ & $-$ & $-$& $-$ & $-$ & $-$ & $-$ \\ \cline{2-17}
		     & $\rmLR$ & 2 & 10.38 & 2.38 & 6.3 & 2.4 & 0.32 & 1.68 & 0.220 & 19.09 & 3.32 & 11.4 & 3.2 & 0.78 & 2.53 & 0.418 \\
		     &      & 5 & 10.28 & 2.31 & 5.8 & 2.3 & -0.62 & 0.72 & 2.451 & \textbf{18.85} & \textbf{3.19} & \textbf{10.6} & \textbf{3.1} & \textbf{-0.51} & \textbf{1.21} & \textbf{4.920} \\ \cline{2-17}
		     & $\rmRc$ & 2 & 13.17 & 3.59 & 9.5 & 3.6 & 27.26 & 28.98 & 0.003 & 22.13 & 4.32 & 14.6 & 4.3 & 16.83 & 18.85 & 0.005 \\
		     &      & 5 & 11.35 & 3.01 & 7.6 & 3.0 & 9.65 & 11.14 & 0.011 & 20.53 & 3.86 & 13.0 & 3.9 & 8.37 & 10.24 & 0.020 \\ \cline{2-17}
		     & $\rmR$ & 2 & 13.11 & 3.54 & 9.4 & 3.6 & 26.67 & 28.38 & 0.007 & 22.13 & 4.32 & 14.6 & 4.3 & 16.83 & 18.85 & 0.014 \\
		     &      & 5 & 11.16 & 2.97 & 7.3 & 3.0 & 7.86 & 9.32 & 0.056 & 20.62 & 3.97 & 13.1 & 4.0& 8.86 & 10.74 & 0.110 \\
		     \hline
		0.6 & $\rmIP$ & 2 & 10.26 & 2.31 & 5.7 & 2.2 & 0.00 & 0.80 & 0.006 & 19.10 & 3.28 & 10.8 & 3.1 & 0.92 & 2.60 & 0.011 \\
		     &      & 5 &10.22 & 2.26 & 5.2 & 2.2 & -0.44 & 0.35 & 0.071 & 18.81 & 3.06 & 9.7 & 2.9 & -0.62 & 1.03 & 0.138 \\ \cline{2-17}
		     & $\rmLB$ & 2 & 10.27 & 2.31 & 5.7 & 2.2 & 0.05 & 0.85 & 0.025 & $-$ & $-$ & $-$ & $-$ & $-$ & $-$ & $-$ \\
		     &      & 5 & 10.26 & 1.88 & 3.6 & 1.9 & -0.04 & 0.76 & 0.209 & $-$ & $-$ & $-$ & $-$ & $-$ & $-$ & $-$ \\  \cline{2-17}
		     & $\rmLR$ & 2 & 10.20 & 2.21 & 5.6 & 2.1 & -0.65 & 0.15 & 0.201 & 18.93 & 3.26 & 10.7 & 3.1 & 0.02 & 1.68 & 0.408 \\
		     &      & 5 & \textbf{10.17} & \textbf{2.13} & \textbf{5.1} & \textbf{2.1} & \textbf{-0.93} & \textbf{-0.14} & \textbf{1.220} & \textbf{18.68} & \textbf{3.05} & \textbf{9.5} & \textbf{2.9} & \textbf{-1.31} & \textbf{0.33} & \textbf{1.170} \\ \cline{2-17}
		     & $\rmRc$ & 2 & 11.48 & 2.75 & 7.1 & 2.7 & 11.86 & 12.76 & 0.002 & 20.67 & 3.63 & 12.1 & 3.5 & 9.20 & 11.01 & 0.004 \\
		     &      & 5 & 10.47 & 2.25 & 5.9 & 2.2 & 2.08 & 2.89 & 0.010 & 19.49 & 3.26 & 10.6 & 3.1 & 2.96 & 4.67 & 0.019 \\ \cline{2-17}
		     & $\rmR$ & 2 & 11.47 & 2.75 & 7.1 & 2.7 & 11.80 & 12.69 & 0.007  & 20.67 & 3.63 & 12.1 & 3.5 & 9.20 & 11.01 & 0.014 \\ 
		     &      & 5 & 10.47 & 2.25 & 5.9 & 2.2 & 1.99 & 2.81 & 0.051 & 19.49 & 3.26 & 10.6 & 3.1 & 2.96 & 4.67 & 0.019 \\
		     \hline
		0.8 & $\rmIP$ & 2 & 10.24 & 2.28 & 5.6 & 2.2 & -0.27 & 0.58 & 0.006 & 19.09 & 3.27 & 10.8 & 3.1 & 0.86 & 2.53 & 0.011 \\
		     &      & 5 & \textbf{10.20} & \textbf{2.22} & \textbf{5.1} & \textbf{2.2} & \textbf{-0.65} & \textbf{0.20} & \textbf{0.070} & 18.82 & 3.08 & 9.6 & 2.9 &-0.56 & 1.09 & 0.137 \\ \cline{2-17}
		     & $\rmLB$ & 2 & 10.24 & 2.28 & 5.6 & 2.2 & -0.23 & 0.62 & 0.021 & $-$& $-$ &$-$& $-$ & $-$ & $-$ & $-$ \\
		     &      & 5 & 10.17 & 2.17 & 4.8 & 2.1 & -0.97 & -0.12 & 0.164 & $-$ & $-$ & $-$ & $-$ & $-$ & $-$ & $-$ \\ \cline{2-17}
		     & $\rmLR$ & 2 & 10.19 & 2.20 & 5.5 & 2.1 & -0.75 & 0.10 & 0.187 & 18.93 & 3.26 & 10.6 & 3.1& -0.01 & 1.65 & 0.416 \\
		     &      & 5 & 10.20 & 2.11 & 5.0 & 2.1 & -0.67 & 0.18 & 0.561 & \textbf{18.65} & \textbf{3.02} & \textbf{9.4} & \textbf{2.9} & \textbf{-1.46} & \textbf{0.17} & \textbf{1.186} \\ \cline{2-17}
		     & $\rmRc$ & 2 & 11.17 & 2.62 & 6.7 & 2.6 & 8.82 & 9.74 & 0.002 & 20.54 & 3.59 & 11.9 & 3.5 & 8.53 & 10.33 & 0.004 \\
		     &      & 5 & 10.34 & 2.22 & 5.6 & 2.1 & 0.73 & 1.59 & 0.010  & 19.49 & 3.20 & 10.5 & 3.1 & 2.95 & 4.66 & 0.019 \\ \cline{2-17}
		     & $\rmR$ & 2 & 11.17 & 2.62 & 6.7 & 2.6 & 8.81 & 9.74 & 0.007 & 20.54 & 3.59 & 11.9 & 3.5 & 8.53 & 10.33 & 0.014 \\
		     &      & 5 & 10.34 & 2.22 & 5.6 & 2.1 & 0.74 & 1.60 & 0.051 & 19.49 & 3.20 & 10.5 & 3.1 & 2.95 & 4.66 & 0.101\\
        \hline
      \end{tabular}}
\caption{Numerical values of $\vert \nu_{\rmp} \vert$, $\rmF_{\rmp}$ (mean and standard errors), $\rmG_{\rmp}^{\rmLR}$ and $\rmG_{\rmp}^{\rmRc}$ obtained on an instance $(M,\gamma)$ with $M \in \{ 15, 20\}$ and $\gamma \in \{ 0.2, 0.4, 0.6, 0.8\}$. The values written in bold indicate the best performances regarding the objective values.}
\label{tab:wkpomdp:maintenance_2}
\end{table}

\begin{table}
  \centering
  \resizebox{!}{6.0cm}{
  \begin{tabular}{|c|cc|ccccccc|ccccccc|}
      \hline
        \multirow{3}{*}{$\gamma$} & \multirow{3}{*}{$\rmp$} & \multirow{3}{*}{$T_\rmr$} &  \multicolumn{7}{c|}{$M=15$} & \multicolumn{7}{c|}{$M=20$}  \\
            & & & \multicolumn{2}{c}{$\vert \nu_{\rmp} \vert$ ($\times 10^3$)} & \multicolumn{2}{c}{$\rmF_{\rmp}$} & $\rmG_{\rmp}^{\rmLR}$ & $\rmG_{\rmp}^{\rmRc}$ & Time & \multicolumn{2}{c}{$\vert \nu_{\rmp} \vert$ ($\times 10^3$)} & \multicolumn{2}{c}{$\rmF_{\rmp}$} & $\rmG_{\rmp}^{\rmLR}$ & $\rmG_{\rmp}^{\rmRc}$ & Time \\
             & & & {Mean} & {Std. err.} & {Mean} & {Std. err.} & $(\%)$ & $(\%)$ & ($\si{s}$) & {Mean} & {Std. err.} & {Mean} & {Std. err.} & $(\%)$ & $(\%)$ & ($\si{s}$) \\
        \hline
        0.2  & $\rmIP$ & 2 & 31.54 & 6.03 & 25.0 & 6.0 & -22.56 & 12.73 & 0.017 & 45.06 & 7.08 & 35.9 & 7.1 & -20.37 & 8.67 & 0.022 \\
		     &      & 5 & 30.89 & 5.88 & 24.0 & 5.9 & -24.14 & 10.43 & 0.591 & 44.28 & 6.90 & 35.1 & 6.9 & -21.74 & 6.80 & 0.660 \\ \cline{2-17}
		     & $\rmLR$ & 2 & 31.19 & 5.85 & 24.6 & 5.8 & -23.40 & 11.50 & 0.714 & 45.02 & 7.02 & 35.8 & 7.0 & -20.44 & 8.57 & 0.944 \\
		     &      & 5  & \textbf{30.68} & \textbf{5.65} & \textbf{23.8} & \textbf{5.7} & \textbf{-24.65} & \textbf{9.68} & \textbf{10.189} & \textbf{44.16} & \textbf{6.77} & \textbf{35.0} & \textbf{6.8} & \textbf{-21.95} & \textbf{6.52} & \textbf{13.309} \\ \cline{2-17}
		     & $\rmRc$ & 2 & 44.74 & 9.41 & 39.3 & 9.5 & 9.87 & 59.94 & 0.007 & 61.27 & 10.59 & 53.4 & 10.7 & 8.28 & 47.77 & 0.009 \\
		     &      & 5 & 44.41 & 9.74 & 38.8 & 9.8 & 9.06 & 58.75 & 0.058 & 59.82 & 10.47 & 52.0 & 10.6 & 5.73 & 44.28 & 0.084 \\ \cline{2-17}
		     & $\rmR$ & 2 & 44.71 & 9.53 & 39.3 & 9.6 & 9.78 & 59.81 & 0.041 & 61.64 & 10.75 & 53.8 & 10.8 & 8.94 & 48.66 & 0.061 \\
		     &      & 5 & 44.18 & 9.46 & 38.6 & 9.5 & 8.50 & 57.94 & 0.212 & 59.58 & 10.57 & 51.8 & 10.7 & 5.31 & 43.71 & 0.294 \\ 
		     \hline
		0.4  & $\rmIP$ & 2 & 28.18 & 4.22 & 19.1 & 4.1 & 0.45 & 1.93 & 0.016 & 41.18 & 4.72 & 23.0 & 4.7 & 0.35 & 1.50 & 0.020 \\
		     &      & 5 & \textbf{27.67} & \textbf{3.97} & \textbf{16.8} & \textbf{3.9} & \textbf{-1.39} & \textbf{0.06} & \textbf{0.232} & 40.83 & 4.66& 22.7 & 4.7 & -0.49 & 0.66 & 0.469 \\ \cline{2-17}
		     & $\rmLR$ & 2 & 28.10 & 4.13 & 19.0 & 4.0 & 0.15 & 1.63 & 0.646 & 41.22 & 4.71 & 23.0 & 4.7 & 0.44 & 1.60 & 0.939 \\
		     &      & 5 & 27.72 & 3.94 & 16.9 & 3.9 & -1.19 & 0.26 & 6.975 & \textbf{40.70} & \textbf{4.64} & \textbf{22.4} & \textbf{4.6} & \textbf{-0.82} & \textbf{0.32} & \textbf{26.353} \\ \cline{2-17}
		     & $\rmRc$ & 2 & 30.48 & 4.64 & 21.0 & 4.5 & 8.62 & 10.22 & 0.006 & 46.28 & 5.59 & 29.4 & 5.6 & 12.78 & 14.07 & 0.009 \\
		     &      & 5 & 29.69 & 4.58 & 20.1 & 4.5 & 5.82 & 7.37 & 0.046  &44.37 & 5.46 & 27.6 & 5.5 & 8.13 & 9.37 & 0.065 \\ \cline{2-17}
		     & $\rmR$ & 2 & 30.48 & 4.64 & 21.0 & 4.5 & 8.62 & 10.22 & 0.039 & 46.27 & 5.58 & 29.4 & 5.6 & 12.76 & 14.06 & 0.056 \\
		     &      & 5 & 29.68 & 4.57 & 20.1 & 4.5 & 5.76 & 7.32 & 0.167 & 44.40 & 5.43 & 27.5 & 5.4 & 8.21 & 9.45 & 0.281 \\
		     \hline
		0.6  & $\rmIP$ & 2 & 28.12 & 4.16 & 18.8 & 4.0 & 0.10 & 1.70 & 0.016 & 40.96 & 4.25 & 18.4 & 4.1 & 0.32 & 1.22 & 0.020 \\
		     &      & 5 & \textbf{27.67} & \textbf{3.84} & \textbf{16.3} & \textbf{3.7} & \textbf{-1.51} & \textbf{0.05} & \textbf{0.225} & 40.72 & 4.18 & 17.9 & 4.0 & -0.28 & 0.62 & 0.316 \\ \cline{2-17}
		     & $\rmLR$ & 2 & 28.10 & 4.09 & 18.8 & 3.9 & 0.04 & 1.63 & 0.620 & 40.83 & 4.12 & 18.3 & 4.0 & -0.02 & 0.89 & 0.831 \\
		     &      & 5 & 27.70 & 3.84 & 16.4 & 3.7 & -1.38 & 0.19 & 1.872 & \textbf{40.57} & \textbf{3.96} & \textbf{16.7} & \textbf{3.8} & \textbf{-0.64} & \textbf{0.25} & \textbf{11.835} \\ \cline{2-17}
		     & $\rmRc$ & 2 & 30.16 & 4.60 & 20.3 & 4.5 & 7.36 & 9.06 & 0.006 & 43.95 & 4.81 & 21.6 & 4.7 & 7.63 & 8.60 & 0.008 \\
		     &      & 5 & 29.55 & 4.38 & 19.5 & 4.2 & 5.19 & 6.87 & 0.045  & 42.58 & 4.63 & 21.3 & 4.5 & 4.26 & 5.20 & 0.059 \\ \cline{2-17}
		     & $\rmR$ & 2 & 30.16 & 4.60 & 20.3 & 4.5 & 7.36 & 9.06 & 0.038  & 43.95 & 4.81 & 21.6 & 4.7 & 7.62 & 8.60 & 0.052 \\
		     &      & 5 & 29.55 & 4.38 & 19.5 & 4.2 & 5.19 & 6.87 & 0.163 & 42.23 & 4.35 & 19.3 & 4.2 & 3.42 & 4.36 & 0.225 \\
		     \hline
		0.8  & $\rmIP$ & 2 & 28.12 & 4.16 & 18.8 & 4.0 & 0.11 & 1.71 & 0.015  & 40.96 & 4.24 & 18.3 & 4.0 & 0.29 & 1.21 & 0.019 \\
		     &      & 5 & \textbf{27.65} & \textbf{3.86} & \textbf{16.2} & \textbf{3.7} & \textbf{-1.57} & \textbf{-0.01} & \textbf{0.226} & 40.76 & 4.09 & 17.8 & 3.9 & -0.19 & 0.73 & 0.313 \\ \cline{2-17}
		     & $\rmLR$ & 2 & 28.10 & 4.09 & 18.8 & 3.9 & 0.03 & 1.62 & 0.637  & 40.83 & 4.09 & 18.1 & 3.9 & -0.04 & 0.88 & 0.829 \\
		     &      & 5  & 27.71 & 3.84 & 16.3 & 3.7 & -1.37 & 0.20 & 1.893 & \textbf{40.56} & \textbf{3.86} & \textbf{16.3} & \textbf{3.7} & \textbf{-0.69} & \textbf{0.22} & \textbf{2.238} \\ \cline{2-17}
		     & $\rmRc$ & 2 & 30.13 & 4.59 & 20.3 & 4.4 & 7.26 & 8.97 & 0.006 & 43.75 & 4.71 & 21.2 & 4.5 & 7.12 & 8.11 & 0.008 \\
		     &      & 5 & 29.52 & 4.37 & 19.5 & 4.2 & 5.10 & 6.77 & 0.044 & 42.40 & 4.57 & 21.0 & 4.4 & 3.81 & 4.77 & 0.057 \\ \cline{2-17}
		     & $\rmR$ & 2 & 30.13 & 4.59 & 20.3 & 4.4 & 7.26 & 8.97 & 0.038  & 43.75 & 4.71 & 21.2 & 4.5 & 7.12 & 8.11 & 0.050 \\
		     &      & 5 & 29.52 & 4.37 & 19.5 & 4.2 & 5.10 & 6.77 & 0.163 & 42.04 & 4.25 & 18.9 & 4.1 & 2.93 & 3.88 & 0.219 \\
        \hline
      \end{tabular}}
\caption{Numerical values of $\vert \nu_{\rmp} \vert$, $\rmF_{\rmp}$ (mean and standard errors), $\rmG_{\rmp}^{\rmLR}$ and $\rmG_{\rmp}^{\rmRc}$ obtained on an instance $(M,\gamma)$ with $M \in \{ 15, 20\}$ and $\gamma \in \{ 0.2, 0.4, 0.6, 0.8\}$. The values written in bold indicate the best performances regarding the objective values.}
\label{tab:wkpomdp:maintenance_3}
\end{table}

One may observe that for all instances, the matheuristic involving policy $\bfdelta^{\rmIP}$ significantly outperforms the matheuristic involving formulations $\{\rmLB,\rmR^\rmc,\rmR\}$ and delivers promising results even in the most challenging instances ($M=20$).
As mentionned in the introduction, it gives better results than the policy $\bfdelta^{\rmR}$, which does not consider the partially observable aspect of the components.
We also observe that the policy involving the Lagrangian relaxation gives competitive results with $\bfdelta^{\rmIP}$ in terms of rewards and failures. However, the standard errors of the Lagrangian relaxation policy are in general higher than the ones of $\rmIP$ and $\rmLR$. This may come from the fact that the action taken in $\rmAct_{T_\rmr}^{\rmLR,t}(\bfh)$ is sampled. 
In Tables~\ref{tab:wkpomdp:maintenance_1},~\ref{tab:wkpomdp:maintenance_2} and~\ref{tab:wkpomdp:maintenance_3}, the negative values of $\rmG_{\rmp}^{\rmRc}$ result from error approximations due to the Monte-Carlo simulations.  
From a solution time perspective, running $\rmAct_{T_\rmr}^{\rmIP,t}(\bfh)$ takes less time than $\rmAct_{T_\rmr}^{\rmLB,t}(\bfh)$, $\rmAct_{T_\rmr}^{\rmLR,t}(\bfh)$, $\rmAct_{T_\rmr}^{\rmR^\rmc,t}(\bfh)$ and $\rmAct_{T_\rmr}^{\rmR,t}(\bfh)$.

Even for the largest instances ($M = 15$ or $M=20$) and for $T=5$, the average time per action of $\rmAct{T_\rmr}^{\rmIP,t}(\bfh)$ is on the order of $1.0$ second; this amount of time is still feasible even if the $24$ decision times are close together.  
Figure~\ref{fig:num:nb_failures} shows that the decomposable policy $\rmLB$ may lead to a lower number of failures and a higher maintenance cost, which means that the policy $\rmLB$ is more conservative. One may observe that all the figures reduce when the maintenance capacity $K$ increases, as expected.

\section{Examples where $z_{\rmIP} < v_{\rmml}^*$ or $z_{\rmIP} > v_{\rmml}^*$}
\label{app:counter_example}

In this section, we describe two instances showing respectively that our MILP~\eqref{pb:wkpomdp:decPOMDP_MILP} is neither an upper bound or a lower bound. 

\subsection{The inequality $z_{\rmIP} \leq v_{\rmml}^*$ does not hold in general.}

Consider a weakly coupled POMDP with $M=2$, $K=1$, $\calX_S^1 = \calX_S^2 = \{1,2,3 \}$, and $\calX_O^1= \calX_O^2 = \{1,2\}$. We set the following initial probability data,
\begin{align*}
	p^1(\cdot) = \begin{bmatrix}
    0.0286 & 0.4429 & 0.5284
  \end{bmatrix}, & \quad & 
  	p^2(\cdot) = \begin{bmatrix}
    0.5328 & 0.2202 & 0.2469
  \end{bmatrix},
\end{align*}
the following transition probability data,
\begin{align*}
	p^1(\cdot|\cdot,0) = \begin{bmatrix}
    0.3149 & 0.2598 & 0.4253 \\
 	0.2542 & 0.5195 & 0.2263 \\
 	0.2016 & 0.7551 & 0.0433
  \end{bmatrix}, & \quad & 
  	p^2(\cdot|\cdot,0) = \begin{bmatrix}
    0.6833 & 0.1797 & 0.1371  \\
 	0.0398 & 0.9207 & 0.0394 \\
 	0.1422 & 0.2202 & 0.6376 
  \end{bmatrix}, \\
  	p^1(\cdot|\cdot,0) = \begin{bmatrix}
    0.3849 & 0.2891 & 0.3260 \\
 	0.4462 & 0.1346 & 0.4192 \\
 	0.0418 & 0.5297 & 0.4285
  	\end{bmatrix}, & \quad & 
  	p^2(\cdot|\cdot,1) = \begin{bmatrix}
    0.4665 & 0.0956 & 0.4379 \\
 	0.4510 & 0.5168 & 0.0322 \\
 	0.5864 & 0.2903 & 0.1234 
  	\end{bmatrix},
\end{align*}
the following emission probability data,
\begin{align*}
	p^1(\cdot|\cdot) = \begin{bmatrix}
    0.6823 & 0.3177 \\
 	0.0806 & 0.9194 \\
 	0.5018 & 0.4982
  \end{bmatrix}, & \quad & 
  	p^2(\cdot|\cdot) = \begin{bmatrix}
    0.4389 & 0.5611 \\
 	0.6657 & 0.3343 \\
 	0.1207 & 0.8793
  \end{bmatrix},
\end{align*}
and the following reward data
\begin{align*}
	r^1(\cdot|\cdot,0) = \begin{bmatrix}
    3.3101  & 7.8198 & 6.9773 \\
 	2.0722 & 2.6782 & 3.5715 \\
 	8.4428 & 2.6010 & 3.2765 
  \end{bmatrix}, & \quad & 
  	r^2(\cdot|\cdot,0) = \begin{bmatrix}
    2.9600 & 8.1503 & 4.5911 \\
 	2.2638 & 6.0290 & 2.5511 \\
 	8.0789 & 7.9927 & 5.0259 \\
  \end{bmatrix}, \\
  r^1(\cdot|\cdot,1) = \begin{bmatrix}
    1.9315 & 9.3614 & 2.8927 \\
 	4.8769 & 5.3131 & 7.3626 \\
 	3.7944 & 4.5557 & 8.6462 
  \end{bmatrix}, & \quad & 
  	r^2(\cdot|\cdot,1) = \begin{bmatrix}
    6.2647 & 6.6832  & 1.1263 \\
 	9.9182 & 9.0278  & 5.9492 \\
 	9.8333 & 0.4466 & 4.3798 
  \end{bmatrix}.
\end{align*}
Solving~\ref{pb:decPOMDP_wc} with $T=4$ using MILP~\eqref{pb:pomdp:MILP_pomdp} on $\calX_S$, $\calX_O$ and $\calX_A$, we obtain an optimal value of $v_{\rm{ml}}^* = 44.7122$, while the optimal value of our MILP~\eqref{pb:wkpomdp:decPOMDP_MILP} is $z_{\rmIP} = 44.2834$. Hence, we obtain $z_{\rmIP} <  v_{\rmml}^*$. Therefore, $v_{\rmml}^* \leq z_{\rmIP}$ does not hold in general.

\subsection{The inequality $z_{\rmIP} \geq v_{\rm{ml}}^*$ does not hold in general.}

Consider a weakly coupled POMDP with $M=2$, $K=1$, $\calX_S^1 = \calX_S^2 = \{1,2,3 \}$, and $\calX_O^1= \calX_O^2 = \{1,2\}$. We set the following initial probability data,
\begin{align*}
	p^1(\cdot) = \begin{bmatrix}
    0.4311 & 0.5255 & 0.0434
  \end{bmatrix}, & \quad & 
  	p^2(\cdot) = \begin{bmatrix}
    0.4835 & 0.1745 & 0.3421
  \end{bmatrix},
\end{align*}
the following transition probability data,
\begin{align*}
	p^1(\cdot|\cdot,0) = \begin{bmatrix}
    0.1517 & 0.3481 & 0.5002 \\
 	0.1639 & 0.0922 & 0.7439 \\
 	0.3395 & 0.2385 & 0.4220 
  \end{bmatrix}, & \quad & 
  	p^2(\cdot|\cdot,0) = \begin{bmatrix}
    0.3435 & 0.3291 & 0.3274 \\
 	0.5964 & 0.1653 & 0.2383 \\
 	0.3968 & 0.2626 & 0.3406 
  \end{bmatrix}, \\
  	p^1(\cdot|\cdot,1) = \begin{bmatrix}
    0.3467 & 0.2733  & 0.3800 \\
 	0.5027 & 0.3548 & 0.1425 \\
 	0.2530 & 0.5466 & 0.2003 
  	\end{bmatrix}, & \quad & 
  	p^2(\cdot|\cdot,1) = \begin{bmatrix}
    0.3160 & 0.4210 & 0.2630  \\
 	0.3583 & 0.3882 & 0.2535   \\
 	0.3611 & 0.4308 & 0.2081 
  	\end{bmatrix},
\end{align*}
the following emission probability data,
\begin{align*}
	p^1(\cdot|\cdot) = \begin{bmatrix}
    0.2052 & 0.7948 \\
 	0.8296 & 0.1704 \\
 	0.5330 & 0.4670 
  \end{bmatrix}, & \quad & 
  	p^2(\cdot|\cdot) = \begin{bmatrix}
    0.6273 & 0.3727  \\
 	0.0392 & 0.9608 \\
 	0.4024 & 0.5976 
  \end{bmatrix},
\end{align*}
and the following reward data
\begin{align*}
	r^1(\cdot|\cdot,0) = \begin{bmatrix}
    7.0075 & 6.2135 & 8.4122 \\
 	9.7198 & 9.5152 & 2.6182 \\
 	1.8522 & 7.4390 & 4.9132 
  \end{bmatrix}, & \quad & 
  	r^2(\cdot|\cdot,0) = \begin{bmatrix}
    8.7418 & 2.6682 & 2.5227  \\
 	8.7673 & 6.1198 & 6.4814  \\
 	6.4971 & 3.8810 & 0.3476 
  \end{bmatrix}, \\
  r^1(\cdot|\cdot,1) = \begin{bmatrix}
    2.8154 & 7.0215 & 1.6752 \\
 	7.8149 & 0.7849 & 4.3722 \\
 	5.9378 & 9.1273 & 1.1657 
  \end{bmatrix}, & \quad & 
  	r^2(\cdot|\cdot,1) = \begin{bmatrix}
    7.4528 & 8.5013 & 9.1925 \\
 	4.3003 & 2.0946 & 4.2973 \\
 	4.2865 & 0.8470 & 9.5848 
  \end{bmatrix}.
\end{align*}
Solving~\ref{pb:decPOMDP_wc} with $T=4$ using MILP~\eqref{pb:pomdp:MILP_pomdp} on $\calX_S$, $\calX_O$ and $\calX_A$, we obtain an optimal value of $v_{\rmml}^* = 47.3693$, while the optimal value of our MILP~\eqref{pb:wkpomdp:decPOMDP_MILP} is $z_{\rmIP} = 47.7356$. Hence, we obtain $v_{\rmml}^* < z_{\rmIP}$. Therefore, $v_{\rmml}^* \geq z_{\rmIP}$ does not hold in general.

\section{Links between the linear relaxation of MILP~\eqref{pb:wkpomdp:decPOMDP_MILP} and the fluid formulation of~\citet{bertsimas2016decomposable}}
\label{app:linksBertsimas}

In this appendix, we show that when we consider decomposable POMDP (i.e., when the action space does not necessary decompose along the components), the linear relaxation of our MILP~\eqref{pb:wkpomdp:decPOMDP_MILP} is equivalent to the fluid formulation of \citet{bertsimas2016decomposable} over a finite horizon and without discounting applied on the MDP relaxation of the problem.


Consider a decomposable MDP $\left((\calX_S^m,\pfrak^m),\calX_A\right)$, where $\pfrak^m = \left( (p^m(s))_{s \in \calX_S^m}, (p^m(s'|s,a))_{\substack{s,s' \in \calX_S \\ a \in \calX_A}} \right)$. Using Remark~\ref{rem:problem:decomposablePOMDP}, we transform the decomposable MDP into a weakly coupled MDP, where the action space is $\tilde{\calX_A} = \left\{ \bfa \in \calX_A^1\times \cdots \times \calX_A^M \colon \mathds{1}_{a}(a^m) - \mathds{1}_{a}(a^{m+1}) = 0, \enskip \forall a \in \calX_A \right\}$.  
Now we write the fluid formulation of \citet[Problem (3)]{bertsimas2016decomposable} with our notation for a finite horizon, i.e., without considering constraints at time $t' \geq T+1$ and by setting a discount factor $\beta=1$:
\begin{subequations}\label{pb:app_linksBertsimas:fluid}
 \begin{alignat}{2}
\max_{\bfx,\bfA} \enskip & \sum_{t=1}^T \sum_{m=1}^M \sum_{\substack{s,s' \in \calX_S^m \\ a \in \calX_A}} r_m(s,a,s')x_{sas'}^{t,m} & \quad &\\
\mathrm{s.t.} \enskip 
 & x_{s}^{1,m} = \sum_{a' \in \calX_A, s' \in \calX_S^m} x_{sa's'}^{1,m} & \forall m \in  [M], s \in \calX_S \label{eq:app_linksBertsimas:LP_MDP_consistent_action_initial}\\
& \sum_{s'\in \calX_S^m,a' \in \calX_A} x_{s'a's}^{t,m} = \sum_{a' \in \calX_A, s' \in \calX_S^m} x_{sa's'}^{t+1,m} & \forall s \in \calX_S^m, m \in [M], t \in [T] \label{eq:app_linksBertsimas:LP_MDP_consistent_action}\\
& x_{s}^{1,m} = p^m(s) & \forall s \in \calX_S^m, m \in [M] \label{eq:app_linksBertsimas:LP_MDP_initial} \\
& x_{sas'}^{t,m} = p^m(s'|s,a) \sum_{s'' \in \calX_S^m}x_{sas''}^{t,m} &  \forall s \in \calX_S^m, a \in \calX_A, m \in [M], t \in [T] \label{eq:app_linksBertsimas:LP_MDP_consistent_state} \\
& \sum_{s,s' \in\calX_S^m} x_{sas'}^{t,m} = A_{a}^{t} &  \forall a \in \calX_A, m \in [M], t \in [T] \label{eq:app_linksBertsimas:same_action}
 \end{alignat}
\end{subequations}

Now, we relate this fluid formulation to the linear relaxation of our integer formulation~\eqref{pb:wkpomdp:decPOMDP_MILP}. 
Our MILP formulation written on the weakly coupled POMDP writes down:
\begin{subequations}\label{pb:app_linksBertsimas:decPOMDP_MILP}
 \begin{alignat}{2}
 \max_{\bftau,\bfdelta} \enskip & \sum_{t=1}^T \sum_{m=1}^M \sum_{\substack{s,s' \in \calX_S^m \\ a \in \calX_A^m}} r^m(s,a,s')\tau_{sas'}^{t,m} & \quad &\\
\mathrm{s.t.} \enskip 
 & \left(\bftau^m,\bfdelta^m\right) \in \calQ^{\mathrm{d}}\left(\calX_S^m, \calX_O^m, \calX_A^m,\pfrak^m\right)   & \forall m \in [M]  \label{eq:app_linksBertsimas:POMDP}\\
 &\sum_{s \in \calX_S^m, o \in \calX_A^m}\tau_{soa}^{t,m} = \tau_a^{t,m} & \forall a \in \calX_A^m, m \in [M], t \in [T]  \\
 &\tau_a^{t,m} = \tau_{a}^{t,m+1} & \forall a \in \calX_A, m \in [M-1], t \in [T]
 \end{alignat}
\end{subequations}
The following proposition states that the linear relaxation of our MILP~\eqref{pb:app_linksBertsimas:decPOMDP_MILP} is equivalent to the fluid formulation written on the decomposable MDP $((\calX_S^m,\pfrak^m,\bfr^m)_{m\in [M]}, \calX_A)$.
We denote by $z_{\rmF}$ the optimal value of fluid formulation~\eqref{pb:app_linksBertsimas:fluid}.

\begin{theo}\label{theo:app_linksBertsimas:linear_relaxation}
    The linear relaxation of MILP~\eqref{pb:app_linksBertsimas:decPOMDP_MILP} is equivalent to the fluid formulation~\eqref{pb:app_linksBertsimas:fluid}. In particular, $z_{\rmF} = z_{\rmR}$.
\end{theo}

The equivalence in Theorem~\ref{theo:app_linksBertsimas:linear_relaxation} should to be understood in the sense that there exists a solution of MILP~\eqref{pb:app_linksBertsimas:decPOMDP_MILP} if, and only if there exists a feasible solution of the fluid formulation~\eqref{pb:app_linksBertsimas:fluid} with the same objective value.

\proof[Proof of Theorem~\ref{theo:app_linksBertsimas:linear_relaxation}]
    Let $(\bfx,\bfA)$ be a feasible solution of linear program~\eqref{pb:app_linksBertsimas:fluid}. 
    We set $\tau_s^{1,m} = x_s^{1,m}$ and $\tau_{sas'}^{t,m} = x_{sas'}^{t,m}$ for every $s,s' \in \calX_S^m$, $a \in \calX_A$, $m \in [M]$ and $t \in [T].$
    We also define variables $(\tau_{soa}^{t,m})_{s,o,a,t,m}$ and $(\delta_{a|o}^{t,m})_{o,a,t,m}$ using the definitions~\eqref{eq:pomdp:proof_mu_soa} and~\eqref{eq:pomdp:proof_delta} for each component $m$.
    Then, $(\bftau^m,\bfdelta^m)$ satisfies all the constraints of $Q^{\rm{d}}(\calX_S^m,\calX_O^m,\calX_A,\pfrak^m)$ when the variable $\bfdelta^m$ are continuous. 
    It remains to show that constraints~\eqref{eq:app_linksBertsimas:same_action} are satisfied.
    To do so, we set $\tau_a^{t,m} = A_a^{t}$ for every $a \in \calX_A$, $m \in [M]$ and $t \in [T]$.
    We obtain that:
    \begin{align*}
        &\tau_a^{t,m}= A_a^t = \tau_a^{t,m+1}, \quad \text{and}, \quad
        &\tau_a^{t,m}= A_a^t =  \sum_{s,s' \in \calX_S^m} \tau_{sas'}^{t,m}= \sum_{\substack{s,s' \in \calX_S^m \\ o \in \calX_S^m}} \tau_{soa}^{t,m}p^m(s'|s,a) = \sum_{s \in \calX_S^m, o \in \calX_O^m} \tau_{soa}^{t,m},
    \end{align*}
    which ensure that $(\bftau^m,\bfdelta^m)_{m \in [M]}$ is a feasible solution of the linear relaxation of Problem~\eqref{pb:app_linksBertsimas:decPOMDP_MILP}. In addition, the objective values are equal.

    Let $(\bftau^m,\bfdelta^m)_{m\in [M]}$ be a feasible solution of the linear relaxation of MILP~\eqref{pb:app_linksBertsimas:decPOMDP_MILP}. 
    We set $\bfx^m = ((\tau_s^{1,m})_{s \in \calX_S^m},(\tau_{sas'}^{t,m})_{\substack{s,s' \in \calX_S^m \\ a \in \calX_A, t \in [T]}})$ and $\bfA = (\tau_a^{t,1})_{a \in \calX_A,t \in [T]}$.
    By definition of $\calQ^{\mathrm{d}}(\calX_S^m,\calX_O^m,\calX_A,\pfrak^m)$, constraints~\eqref{eq:app_linksBertsimas:LP_MDP_consistent_action_initial}-\eqref{eq:app_linksBertsimas:LP_MDP_consistent_action} are satisfied. It remains to show that constraints~\eqref{eq:app_linksBertsimas:same_action} are satisfied.
    First, since $\tau_a^{t,m} = \tau_{a}^{t,m+1}$ for every $ m \in [M-1]$, then $A_a^t = \tau_a^{t,m}$ for every $m$ in $[M]$.
    This comes from the following computation:
    \begin{align*}
        A_a^{t} = \sum{s,o} \tau_{soa}^{t,m} = \sum_{s} \underbrace{\sum_{o} \tau_{soa}^{t,m}}_{= \sum_{s''} \tau_{sas''}^{t,m}} = \sum_{s,s'} \tau_{sas''}^{t,m}
    \end{align*}
    In addition, the objective values are equal.
    It achieves the proof.
\qed

\section{Definition of an implicit policy based on the linear relaxation of MILP~\eqref{pb:wkpomdp:decPOMDP_MILP}}
\label{app:LPpolicy}

In this appendix, we introduce three algorithms $\rm{Act}_{T}^{\rmR,t}(\bfh)$, $\rm{Act}_{T}^{\rmR^\rmc,t}(\bfh)$ and $\rm{Act}_{T}^{\rmLR,t}(\bfh)$ that are slightly different to $\rm{Act}_{T}^{IP,t}(\bfh)$ and that respectively involve the linear relaxation of MILP~\eqref{pb:wkpomdp:decPOMDP_MILP} without and with valid inequalities~\eqref{eq:wkpomdp:dec_Valid_cuts}, and the Lagrangian relaxation.
Algorithm $\rm{Act}_{T}^{IP,t}(\bfh)$ needs to be modified because the linear relaxation of MILP~\eqref{pb:wkpomdp:decPOMDP_MILP} and the Lagrangian relaxation do not provide policies $\bfdelta^m$ such that the action taken at step~\ref{alg:wkpomdp:take_action} belongs to $\calX_A$. For this reason we slightly modify this step in the following algorithm:

\begin{algorithm}[H]
\caption{History-dependent policy $\rm{Act}_{T}^{\rmR,t}(\bfh)$}
\label{alg:app_LPpolicy:heuristic_individual}
\begin{algorithmic}[1]
\STATE \textbf{Input} An history of observations and actions $\bfh \in (\calX_O \times \calX_A)^{t-1}\times \calX_O$.
\STATE \textbf{Output} An action $\bfa \in \calX_A.$
\STATE Compute the belief state $p^m(s|h^m)$ according to the belief state update for every state $s$ in $\calX_S^m$ and every component $m$.
\STATE Remove constraints and variables indexed by $t'<t$ in MILP~\eqref{pb:wkpomdp:decPOMDP_MILP} and solve the \textbf{linear relaxation} of the resulting problem with horizon $T - t$, initial probability distributions $\left(p^m(s|h^m)\right)_{s \in \calX_S^m}$ for every component $m$ in $[M]$ and initial observation $\bfo$ (see Remark~\ref{rem:pomdp:with_observation}) to obtain an optimal solution $(\bftau^m,\bfdelta^m)_{m\in [M]}.$
\label{alg:app_LPpolicy:modify_constraints}
\STATE Choose $\bfa \notin \calX_A$
\WHILE{$\bfa \notin \calX_A$}\label{alg:app_LPpolicy:heuristic_individual:loop_while}
\STATE Sample $a^m$ according to the probability distribution $(\tau_a^{t,m})_{a \in \calX_A^m}$ for all $m$ in $[M].$ \label{alg:app_LPpolicy:take_action}
\ENDWHILE
\STATE Return $\bfa$
\end{algorithmic}
\end{algorithm}

Then we define the implicit policy $\bfdelta$ as follows:
\begin{align}\label{eq:app_LPpolicy:implicit_policy}
    \delta_{\bfa|\bfh}^{\rmR,t} = \begin{cases}
                            1, & \text{if}\ \bfa=\mathrm{Act}_{T}^{\rmR,t}(\bfh) \\
                            0, & \text{otherwise}
                        \end{cases}, & \quad \forall \bfh \in \calX_H^t, \enskip \bfa \in \calX_A, \enskip  t \in [T],
\end{align}
Similarly, we define algorithm $\rm{Act}_{T}^{\rmR^\rmc,t}(\bfh)$ by adding valid inequalities~\eqref{eq:wkpomdp:dec_Valid_cuts} in the resolution of the linear formulation at step~\ref{alg:app_LPpolicy:modify_constraints}.
We also define $\rm{Act}_{T}^{\rmLR,t}(\bfh)$ by solving the master problem~\eqref{pb:app_ColGen:master} using the column generation approach (see Appendix~\ref{app:ColGen}) at step~\ref{alg:app_LPpolicy:modify_constraints}.

We assume that $\calX_A \subsetneq \calX_A^1 \times \cdots \times \calX_A^M$. Otherwise, solving~\ref{pb:POMDP_perfectRecall} can be solved by solving the subproblems independently. 
At first sight, there is no reason to believe that the loop starting at step~\ref{alg:app_LPpolicy:heuristic_individual:loop_while} ends in a finite number iterations. 
It turns out the theorem below ensures that algorithm $\rm{Act}_{T}^{\rmR,t}(\bfh)$ ends in a finite number of iteration.

\begin{theo}\label{theo:app_LPpolicy:alg_finite}
    Algorithm $\rm{Act}_{T}^{\rmR,t}(\bfh)$ (resp. $\rm{Act}_{T}^{\rmR^\rmc,t}(\bfh)$ and $\rm{Act}_{T}^{\rmLR,t}(\bfh)$) ends in a finite number of iterations.
\end{theo}

In fact, unlike \citet{bertsimas2016decomposable}, we cannot choose the action $a^m \in \argmax_{a} \tau_a^{t,m}$ for every $m$ in $[M]$ because there is no guarantee that the induced action $\bfa$ will belong to $\calX_A$.
For this reason, we choose to sample according to $\tau_a^{t,m}$.

\proof[Proof of Theorem~\ref{theo:app_LPpolicy:alg_finite}]
    To prove this result, we denote by $Z$ the stopping time that belongs to $\bbN$ and that represents the number of iterations of the loop starting at step~\ref{alg:app_LPpolicy:heuristic_individual:loop_while}. 
    We want to prove that $\bbP( Z < \infty ) = 1$.  
    To do so we introduce the random variable $\bfA_i$ that represents the sample drawn at iteration $i$.
    Since $\bftau$ is a feasible solution of the linear relaxation of MILP~\eqref{pb:wkpomdp:decPOMDP_MILP}, constraint~\eqref{eq:wkpomdp:decPOMDP_MILP_linking_cons} at time $t$ ensures that $\bbE_{\tilde{\bbP}} \left[\sum_{m=1}^M \bfD^m(A_{i}^m) \right] \leq \bfb$ where $\tilde{\bbP}(\bfA_i = \bfa) = \prod_{m=1}^M \tau_{a^m}^{t,m}$ for every $\bfa \in \calX_A^1 \times \cdots \times \calX_A^M$.
    The random variable $Z$ follows a geometric law with a probability of success $\tilde{\bbP}\left(\bfA_i \in \calX_A \right)$. To prove that $\bbP( Z < \infty ) = 1$, it suffices to prove that $\tilde{\bbP}\left(\bfA_i \in \calX_A \right) > 0$.
    We introduce the quantity $E = \bbE_{\tilde{\bbP}} \left[\sum{m=1}^M \bfD^m(A_i^m) \right]$ and we have the following computation:
    \begin{align*}
        \tilde{\bbP}(\bfA_i \in \calX_A) = \tilde{\bbP}(\sum_{m=1}^M \bfD^m(A_i^m) \leq \bfb) &= 1 - \bbP(\sum_{m=1}^M \bfD^m(A_i^m) - E > \bfb - E)
    \end{align*}
    By contradiction, suppose that $\tilde{\bbP}(\bfA_i \in \calX_A)=0$.
    It follows that $\tilde{\bbP}(\sum_{m=1}^M \bfD^m(A_i^m) - E > \bfb - E) = 1$.
    Since $\bfb \geq E$, we deduce that $\sum_{m=1}^M \bfD^m(A_i^m) - E$ is strictly positive almost surely according to the probability distribution $\tilde{\bbP}$.
    On the other hand, we have $\bbE_{\tilde{\bbP}} \left[ \sum_{m=1}^M \bfD^m(A_i^m) - E \right] =0$.
    Therefore, the random variable $\sum_{m=1}^M \bfD^m(A_i^m) - E$ is a non-negative random variable with an expected value equals to $0$. We deduce that $ \sum_{m=1}^M \bfD^m(A_i^m) - E$ is equal to $0$ almost surely according to $\tilde{\bbP}$, which contradicts the fact that $\tilde{\bbP}(\sum_{m=1}^M \bfD^m(A_i^m) - E > 0)=1$.
    It achieves the proof.
    The proof also holds for Lagrangian relaxation because the feasible solution also satisfies the linking constraint in expectation.
\qed

\end{document}